\documentclass[reqno,12pt,a4paper]{amsart}

\usepackage{mathrsfs}
\usepackage{amssymb}
\usepackage{amsfonts}
\usepackage{latexsym}
\usepackage{amsthm}

\usepackage{graphicx}

\def\lb{\label}

\newcommand{\er}[1]{\textrm{(\ref{#1})}}

\begin{document}

%%%%%%%%%% Some definitions %%%%%%%%%%

%%%%%%%% Equations, theorems %%%%%%%%%
\renewcommand{\theequation}{\arabic{section}.\arabic{equation}}
\theoremstyle{plain}
\newtheorem{theorem}{\bf Theorem}[section]
\newtheorem{lemma}[theorem]{\bf Lemma}
\newtheorem{corollary}[theorem]{\bf Corollary}
\newtheorem{proposition}[theorem]{\bf Proposition}
\newtheorem{definition}[theorem]{\bf Definition}

\newtheorem{remark}[theorem]{\bf Remark}

%%%%% Alphabet %%%%%
\def\a{\alpha}  \def\cA{{\mathcal A}}     \def\bA{{\bf A}}  \def\mA{{\mathscr A}}
\def\b{\beta}   \def\cB{{\mathcal B}}     \def\bB{{\bf B}}  \def\mB{{\mathscr B}}
\def\g{\gamma}  \def\cC{{\mathcal C}}     \def\bC{{\bf C}}  \def\mC{{\mathscr C}}
\def\G{\Gamma}  \def\cD{{\mathcal D}}     \def\bD{{\bf D}}  \def\mD{{\mathscr D}}
\def\d{\delta}  \def\cE{{\mathcal E}}     \def\bE{{\bf E}}  \def\mE{{\mathscr E}}
\def\D{\Delta}  \def\cF{{\mathcal F}}     \def\bF{{\bf F}}  \def\mF{{\mathscr F}}
\def\c{\chi}    \def\cG{{\mathcal G}}     \def\bG{{\bf G}}  \def\mG{{\mathscr G}}
\def\z{\zeta}   \def\cH{{\mathcal H}}     \def\bH{{\bf H}}  \def\mH{{\mathscr H}}
\def\e{\eta}    \def\cI{{\mathcal I}}     \def\bI{{\bf I}}  \def\mI{{\mathscr I}}
\def\p{\psi}    \def\cJ{{\mathcal J}}     \def\bJ{{\bf J}}  \def\mJ{{\mathscr J}}
\def\vT{\Theta} \def\cK{{\mathcal K}}     \def\bK{{\bf K}}  \def\mK{{\mathscr K}}
\def\k{\kappa}  \def\cL{{\mathcal L}}     \def\bL{{\bf L}}  \def\mL{{\mathscr L}}
\def\l{\lambda} \def\cM{{\mathcal M}}     \def\bM{{\bf M}}  \def\mM{{\mathscr M}}
\def\L{\Lambda} \def\cN{{\mathcal N}}     \def\bN{{\bf N}}  \def\mN{{\mathscr N}}
\def\m{\mu}     \def\cO{{\mathcal O}}     \def\bO{{\bf O}}  \def\mO{{\mathscr O}}
\def\n{\nu}     \def\cP{{\mathcal P}}     \def\bP{{\bf P}}  \def\mP{{\mathscr P}}
\def\r{\rho}    \def\cQ{{\mathcal Q}}     \def\bQ{{\bf Q}}  \def\mQ{{\mathscr Q}}
\def\s{\sigma}  \def\cR{{\mathcal R}}     \def\bR{{\bf R}}  \def\mR{{\mathscr R}}
\def\S{\Sigma}  \def\cS{{\mathcal S}}     \def\bS{{\bf S}}  \def\mS{{\mathscr S}}
\def\t{\tau}    \def\cT{{\mathcal T}}     \def\bT{{\bf T}}  \def\mT{{\mathscr T}}
\def\f{\phi}    \def\cU{{\mathcal U}}     \def\bU{{\bf U}}  \def\mU{{\mathscr U}}
\def\F{\Phi}    \def\cV{{\mathcal V}}     \def\bV{{\bf V}}  \def\mV{{\mathscr V}}
\def\P{\Psi}    \def\cW{{\mathcal W}}     \def\bW{{\bf W}}  \def\mW{{\mathscr W}}
\def\o{\omega}  \def\cX{{\mathcal X}}     \def\bX{{\bf X}}  \def\mX{{\mathscr X}}
\def\x{\xi}     \def\cY{{\mathcal Y}}     \def\bY{{\bf Y}}  \def\mY{{\mathscr Y}}
\def\X{\Xi}     \def\cZ{{\mathcal Z}}     \def\bZ{{\bf Z}}  \def\mZ{{\mathscr Z}}
%*********************
\def\be{{\bf e}}
\def\bv{{\bf v}} \def\bu{{\bf u}}
\def\Om{\Omega}
%************************
\def\bbD{\pmb \Delta}
\def\mm{\mathrm m}
\def\mn{\mathrm n}
%*************************

\newcommand{\mc}{\mathscr {c}}

\newcommand{\gA}{\mathfrak{A}}          \newcommand{\ga}{\mathfrak{a}}
\newcommand{\gB}{\mathfrak{B}}          \newcommand{\gb}{\mathfrak{b}}
\newcommand{\gC}{\mathfrak{C}}          \newcommand{\gc}{\mathfrak{c}}
\newcommand{\gD}{\mathfrak{D}}          \newcommand{\gd}{\mathfrak{d}}
\newcommand{\gE}{\mathfrak{E}}
\newcommand{\gF}{\mathfrak{F}}           \newcommand{\gf}{\mathfrak{f}}
\newcommand{\gG}{\mathfrak{G}}           %\newcommand{\gg}{\mathfrak{g}}
\newcommand{\gH}{\mathfrak{H}}           \newcommand{\gh}{\mathfrak{h}}
\newcommand{\gI}{\mathfrak{I}}           \newcommand{\gi}{\mathfrak{i}}
\newcommand{\gJ}{\mathfrak{J}}           \newcommand{\gj}{\mathfrak{j}}
\newcommand{\gK}{\mathfrak{K}}            \newcommand{\gk}{\mathfrak{k}}
\newcommand{\gL}{\mathfrak{L}}            \newcommand{\gl}{\mathfrak{l}}
\newcommand{\gM}{\mathfrak{M}}            \newcommand{\gm}{\mathfrak{m}}
\newcommand{\gN}{\mathfrak{N}}            \newcommand{\gn}{\mathfrak{n}}
\newcommand{\gO}{\mathfrak{O}}
\newcommand{\gP}{\mathfrak{P}}             \newcommand{\gp}{\mathfrak{p}}
\newcommand{\gQ}{\mathfrak{Q}}             \newcommand{\gq}{\mathfrak{q}}
\newcommand{\gR}{\mathfrak{R}}             \newcommand{\gr}{\mathfrak{r}}
\newcommand{\gS}{\mathfrak{S}}              \newcommand{\gs}{\mathfrak{s}}
\newcommand{\gT}{\mathfrak{T}}             \newcommand{\gt}{\mathfrak{t}}
\newcommand{\gU}{\mathfrak{U}}             \newcommand{\gu}{\mathfrak{u}}
\newcommand{\gV}{\mathfrak{V}}             \newcommand{\gv}{\mathfrak{v}}
\newcommand{\gW}{\mathfrak{W}}             \newcommand{\gw}{\mathfrak{w}}
\newcommand{\gX}{\mathfrak{X}}               \newcommand{\gx}{\mathfrak{x}}
\newcommand{\gY}{\mathfrak{Y}}              \newcommand{\gy}{\mathfrak{y}}
\newcommand{\gZ}{\mathfrak{Z}}             \newcommand{\gz}{\mathfrak{z}}

\def\ve{\varepsilon}   \def\vt{\vartheta}    \def\vp{\varphi}    \def\vk{\varkappa}

\def\A{{\mathbb A}} \def\B{{\mathbb B}} \def\C{{\mathbb C}}
\def\dD{{\mathbb D}} \def\E{{\mathbb E}} \def\dF{{\mathbb F}} \def\dG{{\mathbb G}} \def\H{{\mathbb H}}\def\I{{\mathbb I}} \def\J{{\mathbb J}} \def\K{{\mathbb K}} \def\dL{{\mathbb L}}\def\M{{\mathbb M}} \def\N{{\mathbb N}} \def\O{{\mathbb O}} \def\dP{{\mathbb P}} \def\R{{\mathbb R}}\def\S{{\mathbb S}} \def\T{{\mathbb T}} \def\U{{\mathbb U}} \def\V{{\mathbb V}}\def\W{{\mathbb W}} \def\X{{\mathbb X}} \def\Y{{\mathbb Y}} \def\Z{{\mathbb Z}}

%%%%% Arrows %%%%%

\def\la{\leftarrow}              \def\ra{\rightarrow}            \def\Ra{\Rightarrow}
\def\ua{\uparrow}                \def\da{\downarrow}
\def\lra{\leftrightarrow}        \def\Lra{\Leftrightarrow}

%%%%% Typography %%%%%

\def\lt{\biggl}                  \def\rt{\biggr}
\def\ol{\overline}               \def\wt{\widetilde}
\def\no{\noindent}

%%%%% Math signs %%%%%

\let\ge\geqslant                 \let\le\leqslant
\def\lan{\langle}                \def\ran{\rangle}
\def\/{\over}                    \def\iy{\infty}
\def\sm{\setminus}               \def\es{\emptyset}
\def\ss{\subset}                 \def\ts{\times}
\def\pa{\partial}                \def\os{\oplus}
\def\om{\ominus}                 \def\ev{\equiv}
\def\iint{\int\!\!\!\int}        \def\iintt{\mathop{\int\!\!\int\!\!\dots\!\!\int}\limits}
\def\el2{\ell^{\,2}}             \def\1{1\!\!1}
\def\sh{\sharp}
\def\wh{\widehat}
\def\bs{\backslash}
\def\intl{\int\limits}
%%%%% Math operations %%%%%

\def\na{\mathop{\mathrm{\nabla}}\nolimits}
\def\sh{\mathop{\mathrm{sh}}\nolimits}
\def\ch{\mathop{\mathrm{ch}}\nolimits}
\def\where{\mathop{\mathrm{where}}\nolimits}
\def\all{\mathop{\mathrm{all}}\nolimits}
\def\as{\mathop{\mathrm{as}}\nolimits}
\def\Area{\mathop{\mathrm{Area}}\nolimits}
\def\arg{\mathop{\mathrm{arg}}\nolimits}
\def\const{\mathop{\mathrm{const}}\nolimits}
\def\det{\mathop{\mathrm{det}}\nolimits}
\def\diag{\mathop{\mathrm{diag}}\nolimits}
\def\diam{\mathop{\mathrm{diam}}\nolimits}
\def\dim{\mathop{\mathrm{dim}}\nolimits}
\def\dist{\mathop{\mathrm{dist}}\nolimits}
\def\Im{\mathop{\mathrm{Im}}\nolimits}
\def\Iso{\mathop{\mathrm{Iso}}\nolimits}
\def\Ker{\mathop{\mathrm{Ker}}\nolimits}
\def\Lip{\mathop{\mathrm{Lip}}\nolimits}
\def\rank{\mathop{\mathrm{rank}}\limits}
\def\Ran{\mathop{\mathrm{Ran}}\nolimits}
\def\Re{\mathop{\mathrm{Re}}\nolimits}
\def\Res{\mathop{\mathrm{Res}}\nolimits}
\def\res{\mathop{\mathrm{res}}\limits}
\def\sign{\mathop{\mathrm{sign}}\nolimits}
\def\span{\mathop{\mathrm{span}}\nolimits}
\def\supp{\mathop{\mathrm{supp}}\nolimits}
\def\Tr{\mathop{\mathrm{Tr}}\nolimits}
\def\BBox{\hspace{1mm}\vrule height6pt width5.5pt depth0pt \hspace{6pt}}

%%%%%%%%%%%%% specialities %%%%%%%%%%%%%%

\newcommand\nh[2]{\widehat{#1}\vphantom{#1}^{(#2)}}
%{{\mathop{#1}\limits^\wedge}\vphantom{#1}^{(#2)}}
\def\dia{\diamond}

\def\Oplus{\bigoplus\nolimits}

%%%%%%%%%%% End of definitions %%%%%%%%%%

%%%%% OLD OLD OLD

\def\qqq{\qquad}
\def\qq{\quad}
\let\ge\geqslant
\let\le\leqslant
\let\geq\geqslant
\let\leq\leqslant
\newcommand{\ca}{\begin{cases}}
\newcommand{\ac}{\end{cases}}
\newcommand{\ma}{\begin{pmatrix}}
\newcommand{\am}{\end{pmatrix}}
\renewcommand{\[}{\begin{equation}}
\renewcommand{\]}{\end{equation}}
\def\eq{\begin{equation}}
\def\qe{\end{equation}}
\def\[{\begin{equation}}
\def\bu{\bullet}

\title[{Schr\"odinger operators on periodic discrete graphs}]
{Schr\"odinger operators on periodic discrete graphs}

\date{\today}
\author[Evgeny Korotyaev]{Evgeny Korotyaev}
\address{Mathematical Physics Department, Faculty of Physics, Ulianovskaya 2,
St. Petersburg State University, St. Petersburg, 198904, Russia,
 \ korotyaev@gmail.com,}
\author[Natalia Saburova]{Natalia Saburova}
\address{Department of Mathematical Analysis, Algebra and Geometry, Institute of Mathematics,
Information and Space Technologies, Uritskogo St. 68, Northern (Arctic)
Federal University,
Arkhangelsk, 163002,
 \ n.saburova@gmail.com}


\begin{abstract}
We consider Schr\"odinger operators with periodic potentials on
periodic  discrete  graphs. The spectrum of the Schr\"odinger
operator consists of an absolutely continuous part (a~union of a
finite number of non-degenerated bands)
plus a finite number of flat bands, i.e., eigenvalues  of infinite
multiplicity. We obtain estimates of
the Lebesgue measure of the spectrum in terms of geometric
parameters of the graph and show that they become identities for some class of graphs.
Moreover, we obtain stability estimates and show 
the existence and positions of large number of  flat
bands for specific graphs.
The proof is based on the Floquet theory
and the precise representation of fiber Schr\"odinger operators,
constructed in the paper.

\end{abstract}

\maketitle

\begin{quotation}

\begin{center}
{\bf Table of Contents}
\end{center}

\vskip 6pt

{\footnotesize

1. Introduction \hfill \pageref{Sec1}\ \ \ \ \

2. Main results  \hfill \pageref{Sec2}\ \ \ \ \

3. Direct integrals for Schr\"odinger operators \hfill \pageref{Sec3}\ \ \ \ \

4. Spectral estimates in terms of geometric graph parameters
\hfill \pageref{Sec4}\ \ \ \ \

5.  Schr\"odinger operators on bipartite regular graphs \hfill
\pageref{Sec7b}\ \ \ \ \

6. Stability estimates  \hfill \pageref{Sec5}\ \ \ \ \

7. Various number  of
 flat bands of Laplacians on specific graphs \hfill
\pageref{Sec6}\ \ \ \ \

8. Crystal models \hfill
\pageref{Sec7}\ \ \ \ \

9. Appendix, well-known  properties of matrices \hfill
\pageref{Sec8}\ \ \ \ \ }
\end{quotation}

%%%%%%%%%%%%%%%%%%%%%%%%%%%%%%%%%%%%%%%%%%%%%%%%%%%%%

\vskip 0.25cm

\section {\lb{Sec1}Introduction}
\setcounter{equation}{0}

We  discuss the spectral properties of both Laplacians and
Schr\"odinger operators on $\Z^d$-periodic discrete graphs, $d\ge
2$.  Schr\"odinger operators on periodic graphs are of interest due
to their applications to problems of physics and chemistry. They are
used to study properties of different periodic media, e.g.
nanomedia, see \cite{Ha85}, \cite{NG04} and a nice survey
\cite{CGPNG09}.

There are a lot of papers, and even books,  on the spectrum of
discrete Laplacians  on finite and infinite graphs (see \cite{BK12},
\cite{Ch97}, \cite{CDS95}, \cite{CDGT88}, \cite{P12}  and references
therein). There are results about spectral properties of discrete
Schr\"odinger operators on specific $\Z^d$-periodic graphs.
Schr\"odinger operators with decreasing potentials on the lattice
$\Z^d$ are considered by Boutet de Monvel-Sahbani \cite{BS99},
Isozaki-Korotyaev \cite{IK12}, Rosenblum-Solomjak \cite{RoS09} and
see references therein. Ando \cite{A12} considers the inverse
spectral theory for the discrete Schr\"odinger operators with
finitely supported potentials on the hexagonal lattice.
Gieseker-Kn\"orrer-Trubowitz \cite{GKT93} consider
 Schr\"odinger operators with periodic potentials on the lattice $\Z^2$, the
simplest example of $\Z^2$-periodic graphs.  They study its Bloch
variety and its integrated density of states. Korotyaev-Kutsenko
\cite{KK10} -- \cite{KK10b} study the spectra of the discrete
Schr\"odinger operators on graphene nano-tubes and nano-ribbons in
external fields.

\medskip

\subsection{The definition of Schr\"odinger operators on periodic graphs.}
Let $\G=(V,\cE)$ be a connected graph, possibly  having loops and
multiple edges, where $V$ is the set of its vertices and  $\cE$ is
the set of its unoriented edges. The graphs under consideration are
embedded into $\R^d$. An edge connecting vertices $u$ and $v$ from
$V$ will be denoted as the unordered pair $(u,v)_e\in \cE$ and is
said to be \emph{incident} to the vertices. Vertices $u,v\in V$ will
be called \emph{adjacent} and denoted by $u\sim v$, if $(u,v)_e\in
\cE$. We define the degree ${\vk}_v=\deg v$ of the vertex $v\in V$ as the
number of all its incident edges from $\cE$ (here a loop is counted
twice). Below we consider locally finite
$\Z^d$-periodic graphs $\G$ (see examples
in Figures \ref{ff.0.3}\emph{a}, \ref{ff.FCC}\emph{a}), i.e., graphs satisfying the following conditions:

1) {\it the number of vertices from $V$ in any bounded domain $\ss\R^d$ is
finite;

2) the degree of each vertex is finite;

3) $\G$ has the periods (a basis) $a_1,\ldots,a_d$ in $\R^d$, such that $\G$ is invariant
under translations through the vectors $a_1,\ldots,a_d$:}
$$
\G+a_s=\G, \qqq  \forall \ s\in\N_d=\{1,\ldots,d\}.
$$

In the space $\R^d$ we consider a coordinate system with the origin at some point $O$. The coordinate axes of this system are directed along the vectors $a_1,\ldots,a_d$. Below the coordinates of all vertices of $\G$ will be expressed  in this coordinate system.
From the definition it follows that a $\Z^d$-periodic graph $\G$ is invariant under translations through any integer vector $m$:
$$
\G+m=\G,\qqq \forall\, m\in\Z^d.
$$

Let $\ell^2(V)$ be the
Hilbert space of all square summable functions $f:V\to \C$, equipped
with the norm
$$
\|f\|^2_{\ell^2(V)}=\sum_{v\in V}|f(v)|^2<\infty.
$$
We  define the self-adjoint Laplacian  (or the Laplace operator) $\D$ on
$f\in\ell^2(V)$ by
\[
\lb{DLO}
 \big(\D f\big)(v)= \sum\limits_{(v,\,u)_e\in\cE}\big(f(v)-f(u)\big), \qqq
 v\in V.
\]
We recall basic facts about the spectrum for both finite and
periodic graphs (see \cite{Me94}, \cite{M91}, \cite{M92},
\cite{MW89}): \emph{the point 0 belongs to the spectrum $\s(\D)$ containning in  $[0,2\vk_+]$, i.e., }
\[
\lb{bf}
\begin{aligned}
0\in\s(\D)\subset[0,2\vk_+],\qqq
\textrm{where}\qqq
\vk_+=\sup_{ v\in V}\deg v<\infty.
\end{aligned}
\]

We consider the Schr\"odinger operator $H$ acting on the Hilbert
space $\ell^2(V)$ and given by
\[
\lb{Sh}
H=\D+Q,
\]
\[
\lb{Pot}
\big(Q f\big)(v)=Q(v)f(v),\qqq \forall v\in V,
\]
where we assume that the potential $Q$ is real valued and satisfies
$$
Q(v+\wt a_s)=Q(v), \qqq  \forall\, (v,s)\in V\ts\N_d,
$$
for some linearly independent integer vectors  $\wt a_1,\ldots,\wt a_d\in\Z^d$ (in the basis $a_1,\ldots,a_d$).
The vectors $\wt a_1,\ldots,\wt a_d$ are called \emph{the periods of the potential} $Q$.
Since the periods $\wt a_1,\ldots,\wt a_d$ of the potential are also periods of the periodic graph, we may assume that the periods of the potential are the same as the periods of the graph.

\subsection{The definitions of fundamental
graphs and edge indices.} In order to define the Floquet-Bloch
decomposition \er{raz} of Schr\"odinger operators we need to
introduce the two oriented edges $(u,v)$ and $(v,u)$ for each
unoriented edge $(u,v)_e\in \cE$: the oriented edge starting at
$u\in V$ and ending at $v\in V$ will be denoted as the ordered pair
$(u,v)$. We denote the set  of all oriented edges by $\cA$.

We define \emph{the fundamental graph} $\G_f=(V_f,\cE_f)$ of the
periodic graph $\Gamma$ as a graph on the surface $\R^d/\Z^d$ by
\[
\lb{G0} \G_f=\G/{\Z}^d\ss \R^d/\Z^d.
\]
The fundamental graph $\G_f$ has the vertex set $V_f$, the set
$\cE_f$ of unoriented edges and the set $\cA_f$ of oriented edges, which are finite (see Proposition \ref{pro0}.i).
Denote by $v_1,\ldots,v_\nu$ the vertices of $V_f$, where $\nu<\iy$ is the number of vertices of $\G_f$.
We identify them
%$\G_f=(V_f,\cE_f)$
with the vertices of  $V$ from the set $[0,1)^d$ by
\[
\lb{V0} V_f=[0,1)^d\cap V=\{v_1, v_2,\ldots,v_\n\},
\]
see Fig.\ref{ff.0.11}.  Due to \er{V0} for any
$v\in V$ the following unique representation holds true:
\[
\lb{Dv} v=[v]+\tilde v, \qquad [v]\in\Z^d,\qquad \tilde v\in
V_f\subset[0,1)^d.
\]
In other words, each vertex $v$ can be represented uniquely as the
sum of an integer part $[v]\in \Z^d$ and a fractional part $\tilde
v$ that is a vertex of the fundamental graph $\G_f$.
We introduce {\it an edge index}, which is
important to study the spectrum of Schr\"odinger operators on periodic graphs.
For any
oriented edge $\be=(u,v)\in\cA$ we define {\bf the edge "index"}
$\t({\bf e})$ as the integer vector by
\[
\lb{in}
\t({\bf e})=[v]-[u]\in\Z^d,
\]
where due to \er{Dv} we have
$$
u=[u]+\tilde{u},\qquad v=[v]+\tilde{v}, \qquad [u],
[v]\in\Z^d,\qquad \tilde{u},\tilde{v}\in V_f.
$$
If $\be=(u,v)$ is an oriented edge of the graph $\G$, then by the
definition of the fundamental graph there is an oriented edge
$\tilde\be=(\tilde u,\tilde v\,)$ on $\G_f$. For the edge
$\tilde\be\in\cA_f$ we define the edge index $\t(\tilde{\bf e})$ by
\[
\lb{inf}
\t(\tilde{\bf e})=\t(\be).
\]
In other words, edge indices of the fundamental graph $\G_f$  are
induced by edge indices of the periodic graph $\G$.
In a fixed coordinate system the index of the
fundamental graph edge is uniquely determined by \er{inf}, since
due to Proposition \ref{pro0}.ii.c, we have
$$
\t(\be+m)=\t(\be),\qqq \forall\, (\be,m)\in\cA \ts \Z^d.
$$
But generally speaking, the edge indices  depend on the choice of the coordinate
origin $O$.
Edges with nonzero indices will be called {\bf bridges} (see
Fig.\ref{ff.0.11}). They  are important to describe the
spectrum of the Schr\"odinger operator. The bridges provide the
connectivity of the periodic graph and the removal of all bridges
disconnects the graph into infinitely many connected components. The
set of all bridges of the fundamental graph $\G_f$ we denote by
$\cB_f$.

\setlength{\unitlength}{1.0mm}
\begin{figure}[h]
\centering
\unitlength 1mm % = 2.845pt
\linethickness{0.4pt}
\ifx\plotpoint\undefined\newsavebox{\plotpoint}\fi % GNUPLOT compatibility

\begin{picture}(60,50)(0,0)

\multiput(10,10)(4,0){10}{\line(1,0){2}}
\multiput(10,30)(4,0){10}{\line(1,0){2}}
\multiput(10,50)(4,0){10}{\line(1,0){2}}

\multiput(10,10)(0,4){10}{\line(0,1){2}}
\multiput(30,10)(0,4){10}{\line(0,1){2}}
\multiput(50,10)(0,4){10}{\line(0,1){2}}
\put(10,10){\circle{1}}

\put(10,10){\vector(1,0){20.00}}
\put(10,10){\vector(0,1){20.00}}

\put(7.0,8.0){$\scriptstyle O$}

\put(20,8){$\scriptstyle a_1$}
\put(6,20){$\scriptstyle a_2$}

\put(15,15){\line(1,0){10.00}}
\put(25,25){\line(1,0){5.00}}
\put(25,25.1){\line(1,0){5.00}}
\put(25,24.9){\line(1,0){5.00}}
\put(25,24.8){\line(1,0){5.00}}

\put(25,25.1){\line(2,1){20.00}}
\put(25,25.2){\line(2,1){20.00}}
\put(25,24.9){\line(2,1){20.00}}
\put(25,24.8){\line(2,1){20.00}}
\put(25,25){\line(2,1){20.00}}

\put(15,15){\line(-1,2){5.00}}

\put(15,15.1){\line(3,2){15.00}}
\put(15,15.2){\line(3,2){15.00}}
\put(15,14.9){\line(3,2){15.00}}
\put(15,14.8){\line(3,2){15.00}}
\put(15,15){\line(3,2){15.00}}

\put(15,15.4){\line(1,3){10.00}}
\put(15,15.2){\line(1,3){10.00}}
\put(15,15){\line(1,3){10.00}}
\put(15,14.8){\line(1,3){10.00}}
\put(15,14.6){\line(1,3){10.00}}
\put(15,15){\line(1,3){10.00}}

\put(20,20){\line(1,1){5.00}}

\put(15,15){\circle{1}}
\put(10,25){\circle{1}}
\put(20,20){\circle{1}}
\put(25,25){\circle{1}}
\put(25,15){\circle{1}}
\put(12,13){$\scriptstyle v_1$}
\put(11,25){$\scriptstyle v_2$}\put(31,25){$\scriptstyle v_2+a_1$}
\put(26,13){$\scriptstyle v_4$}\put(46,33){$\scriptstyle v_4+a_1+a_2$}
\put(25,27){$\scriptstyle v_3$}\put(21,47){$\scriptstyle v_3+a_2$}
\put(18.5,22){$\scriptstyle v_5$}

%*************************************

\put(35,15){\circle{1}}
\put(30,25){\circle{1}}
\put(40,20){\circle{1}}
\put(45,25){\circle{1}}
\put(45,15){\circle{1}}

\put(15,35){\circle{1}}
\put(10,45){\circle{1}}
\put(20,40){\circle{1}}
\put(25,45){\circle{1}}
\put(25,35){\circle{1}}

%*************************************

\put(35,35){\circle{1}}
\put(30,45){\circle{1}}
\put(40,40){\circle{1}}
\put(45,45){\circle{1}}
\put(45,35){\circle{1}}
\end{picture}

\vspace{-0.5cm} \caption{ \footnotesize  A graph $\G$
with $\n=5$; only edges of the fundamental graph $\G_f$ are shown;
the bridges $(v_1,v_2+a_1)$, $(v_1,v_3+a_2)$, $(v_3,v_2+a_1)$,
$(v_3,v_4+a_1+a_2)$ of $\G_f$ are marked by bold.} \label{ff.0.11}
\end{figure}
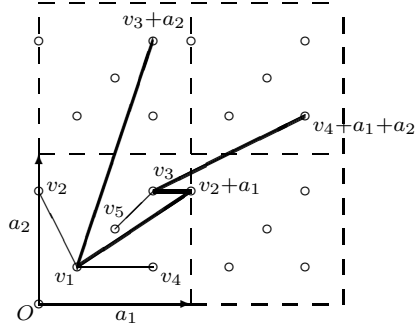

\subsection{Floquet decomposition of Schr\"odinger operators.}
Recall that the fundamental graph $\G_f=(V_f,\cE_f)$  has the finite vertex set
$V_f=\{v_1,\ldots,v_\nu\}\ss [0,1)^d$.
 Due to this
notation, we can denote the potential $Q$ on the fundamental
graph $\G_f$  by
\[\lb{pott}
Q(v_j)=q_j, \qqq j\in \N_\n=\{1,\ldots,\nu\}.
\]

The Schr\"odinger operator $H=\D+Q$ on $\ell^2(V)$ has  the standard
decomposition into a constant fiber direct integral
\[
\lb{raz}
\begin{aligned}
& \ell^2(V)={1\/(2\pi)^d}\int^\oplus_{\T^d}\ell^2(V_f)\,d\vt ,\qqq
UH U^{-1}={1\/(2\pi)^d}\int^\oplus_{\T^d}H(\vt)d\vt,
\end{aligned}
\]
     for some unitary operator $U$. Here
$\ell^2(V_f)=\C^\nu$ is the fiber space  and   $H(\vt)$ is the
Floquet $\nu\ts\nu$ (fiber) matrix  and $\vt\in \T^d=\R^d/(2\pi\Z)^d$ is the
quasimomentum.

Note that the decomposition of discrete Schr\"odinger operators on periodic
graphs into a constant fiber direct integral \er{raz} (without an
exact form of fiber operators)  was
discussed  by Higuchi-Shirai \cite{HS04}, Rabinovich-Roch
\cite{RR07}, Higuchi-Nomura \cite{HN09}. In particular, they prove that the spectrum of Schr\"odinger operators consists  of an
absolutely continuous part and a finite number of flat bands (i.e.,
eigenvalues with infinite multiplicity). The absolutely continuous
spectrum consists of a finite number of intervals (spectral bands)
separated by gaps.

\medskip

\begin{theorem}\label{pro2}
i) The Schr\"odinger operator $H=\D+Q$ acting on $\ell^2(V)$ has the
decomposition into a constant fiber direct integral \er{raz}, where
 the Floquet (fiber) matrix $H(\vt)$ is given by
\[
\lb{Hvt}
H(\vt)=\D(\vt)+q,\qqq q=\diag(q_1,\ldots,q_\n),\qqq \forall\,\vt\in \T^d.
\]
The Floquet matrix
$\D(\vt)=\{\D_{jk}(\vt)\}_{j,k=1}^\n$ for the Laplacian $\D$ is given by
\medskip
\[
\label{l2.15} \D_{jk}(\vt )=\vk_j\d_{jk}-\ca \sum\limits_{{\bf
e}=(v_j,\,v_k)\in{\cA}_f}e^{\,i\lan\t
({\bf e}),\,\vt\ran }, \qq &  {\rm if}\  \ (v_j,v_k)\in \cA_f \\
\qqq 0, &  {\rm if}\  \ (v_j,v_k)\notin \cA_f \ac,
\]
where $\vk_j$ is the degree of $v_j$, $\d_{jk}$ is the Kronecker delta and
$\lan\cdot\,,\cdot\ran$ denotes the standard inner product in
$\R^d$.

ii) Let $H^{(1)}(\vt )$ be a Floquet matrix for $H$ defined by (\ref{Hvt}), (\ref{l2.15}) in another coordinate system with an origin $O_1$. Then the matrices $H^{(1)}(\vt )$ and $H(\vt)$ are unitarily equivalent for all $\vt \in\T^d$.

iii) The entry $\D_{jk}(\cdot)$ of the Floquet matrix $\D(\cdot)
=\{\D_{jk}(\cdot)\}_{j,k=1}^\n$ is constant iff there is no bridge
$(v_j,v_k)\in\cA_f$.

iv) The Floquet matrix  $\D(\cdot)$ has at least one non-constant entry $\D_{jk}(\cdot)$ for some $j\le k$.

\end{theorem}

{\bf Remark.} 1) The identity \er{l2.15} for
the Floquet (fiber) operator is new.
It is important to study
spectral properties of Schr\"odinger operators acting on  graphs.

2) Badanin-Korotyaev-Saburova \cite{BKS13} derived different spectral
properties of normalized Laplacians on $\Z^2$-periodic graphs.
Their proof is based on an exact form of a fiber operator, which is used in our paper.

\section{\lb{Sec2}Main results}
\setcounter{equation}{0}

\medskip

\subsection{Estimates of bands.}
Theorem \ref{pro2} and standard arguments (see Theorem XIII.85 in
\cite{RS78})  describe the spectrum of the Schr\"odinger operator
$H=\D+Q$. Each Floquet ${\nu\ts\nu}$ matrix $H(\vt), \vt\in\T^d$,
has $\n$ eigenvalues $\l_n(\vt)$, $n\in\N_\n$, which are labeled
in increasing order (counting multiplicities) by
\[
\label{eq.3} \l_1(\vt )\leq\l_2(\vt )\leq\ldots\leq\l_{\nu}(\vt),
\qqq \forall\,\vt\in\T^d =\R^d/(2\pi\Z)^d.
\]
Since $H(\vt)$ is
self-adjoint and analytic in $\vt\in\T^d$, each $\l_n(\cdot)$, $n\in\N_\n$, is a real and piecewise analytic function on the torus $\T^d$ and
defines \emph{a dispersion relation}. Define the \emph{spectral bands} $\s_n(H)$ by
\[
\lb{ban.1}
\s_n(H)=[\l_n^-,\l_n^+]=\l_n(\T^d),\qqq n\in\N_{\n}.
\]
 Sy and Sunada
\cite{SS92} show that the lower point of the spectrum $\s(H)$ of the
operator $H$ is $\l_1(0)$, i.e., $\l_1(0)=\l_1^{-}$.
Thus, the spectrum of the operator $H$ on the  graph $\G$ is
given by
\[
\s(H)=\bigcup_{\vt\in\T^d}\s\big(H(\vt)\big)=\bigcup_{n=1}^{\nu}\s_n(H).
\]
Note that if $\l_n(\cdot)= C_n=\const$ on some set $\mB\ss\T^d$ of
positive Lebesgue measure, then  the operator $H$ on $\G$ has the
eigenvalue $C_n$ with infinite multiplicity. We call $C_n$ a
\emph{flat band}. Each
flat band is generated by finitely supported eigenfunction, see \cite{HN09}. Thus, the spectrum of the Schr\"odinger operator
$H$ on the periodic graph $\G$ has the form
\[
\lb{r0}
\s(H)=\s_{ac}(H)\cup \s_{fb}(H).
\]
Here $\s_{ac}(H)$ is the absolutely continuous spectrum, which is a
union  of non-degenerated intervals, and $\s_{fb}(H)$ is the set of
all flat bands (eigenvalues of infinite multiplicity). An open
interval between two neighboring non-degenerated spectral bands is
called a \emph{spectral gap}.

The eigenvalues of the Floquet matrix  $\D(\vt)$ for the Laplacian $\D$ will be denoted by
$\l^0_n(\vt)$, $n\in\N_\n$. The \emph{spectral bands} for the Laplacian $\s_n^0=\s_n(\D), n\in\N_{\n}$, have the form
\[
\lb{ban0} \s_n^0=\s_n(\D)=[\l_n^{0-},\l_n^{0+}]=\l_n^0(\T^d).
\]

\begin{theorem}
\lb{T1}
Let the Schr\"odinger operator $H=\D+Q$ act on $\ell^2(V)$. Then

i)  Let a coefficient $\D_{jk}(\cdot)$ for some $j,k\in\N_\n$
satisfy
\[
\lb{ner}
|\D_{jk}(\cdot)|\neq\const.
\]
Then the first spectral band $\s_1(H)=[\l_1^-,\l_1^+]$  is non-degenerated, i.e., $\l_1^-<\l_1^+$.

ii) The Lebesgue measure $|\s(H)|$ of the spectrum of  $H$ satisfies
\[
\lb{eq.7}
|\s(H)|\le \sum_{n=1}^{\n}|\s_n(H)|\le 2\b,
\]
where $\b$ is the number of fundamental graph bridges.
Moreover, if in the spectrum $\s(H)$ there exist $s$  spectral gaps $\g_1(H),\ldots,\g_s(H)$, then the following estimates hold true:
\[\lb{GEga}
\begin{aligned}
\sum_{n=1}^s|\g_n(H)|\ge
\l^+_\n-\l_1^--2\b\ge C_0-2\b,\\
C_0=\max\{\l^{0+}_\n-q_\bu\,,\, q_\bu-2\vk_+\}, \qq q_\bu=\max_n q_n-\min_n q_n.
\end{aligned}
\]
The estimates \er{eq.7} and the first estimate in \er{GEga} become identities for some classes of
graphs, see \er{es.111}.

\end{theorem}

\no \textbf{Remark.} 1) The total length of spectral  bands depends
essentially on the number of bridges on the fundamental graph
$\G_f$. If we remove the coordinate system, then the number of
bridges on $\G_f$ is changed in general. In order to get the best
estimate in \er{eq.7} we  have to choose a coordinate system in
which  the number $\b$  is minimal.

2) The condition \er{ner} holds true for very large
class of graphs. However, there exist graphs for which this condition does not hold true. An example of such graph is shown in Fig.\ref{f.10}\emph{a}.

3) Below we need the simple, basic fact of Laplacian on periodic graphs:

\no
{\it the first spectral band $\s_1(\D)=[0,\l_1^{0+}]$ of $\D$ is non-degenerated,
i.e., $\l_1^{0+}>\l_1^{0-}=0$}.

\no It can be reformulated:
\emph{the point 0 is never a flat band of $\D$}.

\no Unfortunately, we can not find a paper, proving this fact and we will
prove one in Proposition \ref{pp1}.

4) Sy and Sunada prove that $\l_1^-=\l_1(0)$ for a more general class of graphs including $\Z^d$-periodic graphs (see Theorem~1, p.143 in \cite{SS92}). But they do not discuss the question if the first spectral band is degenerated or not.

5) If $q_\bu>2\vk_+$ ($q_\bu$ is large enough), then $C_0=q_\bu-2\vk_+$.
If $q_\bu<\l^{0+}_\n$ ($q_\bu$ is small enough), then $C_0=\l^{0+}_\n-q_\bu$.

\medskip

We consider the  Schr\"odinger operator $H_t=\D+tQ$, where the
potential $Q$ is  "generic" and   $t\in \R$  is the coupling
constant. We discuss spectral bands of $H_t$ for $t$
 large enough.

\begin{theorem}
\lb{T17} Let the Schr\"odinger operator $H_t=\D+tQ$, where the
potential $Q$ satisfies $q_j\neq q_k$  for all $j,k\in\N_\n, \ j\neq
k$, and the real coupling constant $t$ is large enough. Without loss
of generality we assume that $q_1<q_2<\ldots<q_\n$.  Then each
eigenvalue $\l_n(\vt,t)$ of the corresponding Floquet matrix
$H_t(\vt)$ and each spectral band $\s_n(H_t)$, $n\in\N_\n$,  satisfy
\[
\lb{Qt}
\begin{aligned}
\l_n(\vt,t)=tq_n+\D_{nn}(\vt)-{1\/t}\sum_{j=1 \atop j\neq n}^\n\frac{|\D_{jn}(\vt)|^2}{q_j-q_n}+{O(1)\/t^2}\,,\\
|\s_n(H_t)|=|\D_{nn}(\T^d)|+O(1/t)
\end{aligned}
\]
as $t\to \iy$, uniformly in $\vt\in \T^d$. In particular, we have
\[\lb{Qt1}
|\s(H_t)|=C+O(1/t),\qq \textrm{ where } \qq C=\sum\limits_{n=1}^\n|\D_{nn}(\T^d)|
\]
and $C>0$, if there are bridge-loops on $\G_f$ and $C=0$
if there are no bridge-loops on $\G_f$.
\end{theorem}

\no \textbf{Remark.} Asymptotics  \er{Qt} yield that a small
change  of the potential gives  that all spectral bands of the
Schr\"odinger operator $H_t$ become open for $t$ large enough, i.e.,
the spectrum of $H_t$ is absolutely continuous.

\medskip

{\bf Definition of Loop Graphs.} {\it i) A periodic graph $\G$ is
called a loop  graph if all bridges of some fundamental graph $\G_f$
are loops. This graph $\G_f$ is called a loop fundamental graph.

ii) A loop graph $\G$ is called precise  if $\cos\lan\t
({\bf e}),\,\vt_0\ran=-1$ for all bridges $\be\in\cB_f$ and some
$\vt_0\in\T^d$, where $\t(\be)\in\Z^d$ is the index of a bridge
$\be$ of $\G_f$. This point $\vt_0$ is called a precise quasimomentum of the
loop  graph $\G$.}

\medskip

 A class of all precise loop graphs is large  enough. The simplest example of precise loop graphs is the lattice graph $\dL^d=(V,\cE)$, where the vertex set  and the edge set are given by
\[
\lb{dLg}
V=\Z^d,\qquad  \cE=\big\{(m,m+a_1),\ldots,(m,m+a_d), \quad
\forall\,m\in\Z^d\big\},
\]
and $a_1,\ldots,a_d$ is the standard orthonormal basis, see Fig.\ref{ff.0.1}\emph{a}. The
graph $\dL^d$ has an infinite number of fundamental graphs. The
"minimal"\, fundamental graph $\dL_f^d$ of the lattice $\dL^d$
consists of one vertex $v=0$ and $d$ unoriented edge-loops
$(v,v)$, see Fig.\ref{ff.0.1}\emph{b}. All
bridges of $\dL_f^d$ are loops and their indices have the form $\pm a_1,\pm a_2,\ldots,\pm a_d$. Thus, for the quasimomentum
$\vt_0=(\pi,\ldots,\pi)\in\T^d$  we have $\cos\lan\t ({\bf
e}),\,\vt_0\ran=-1$ for all bridges $\be\in \dL^d_f$ and the graph $\dL^d$ is a
precise loop  graph. It is known that the spectrum of the Laplacian
$\D$ on $\dL^d$ has the form $\s(\D)=\s_{ac}(\D)=[0,4d]$.

We consider perturbations of loop graphs and precise loop graphs.
 A simple example of a precise loop graph $\G_*$ obtained by
perturbations of the square lattice $\dL^2$ is given in
Fig.\ref{ffS'}\emph{a}.

\begin{proposition}\lb{TpLG}
i) There exists a loop graph, which is not precise.

ii) Let $\G=(\cE,V)$ be a loop graph and
let $\G'=(\cE',V')\ss \R^d$ be any connected finite graph such that its diameter is small enough.
We take some points $v\in V$ and $v'\in V'$.
We joint the graph $\G'$ with each point from the vertex set $v+\Z^d$, identifying the vertex  $v'$ with each vertex
of $v+\Z^d$. Then the obtained graph $\G_*$ is a loop graph. Moreover, if $\G$ is precise, then $\G_*$ is also precise and the precise quasimomentum $\vt_0$ of $\G$ is also a precise quasimomentum
of $\G_*$.
\end{proposition}

%*************************************************************
\setlength{\unitlength}{1.0mm}
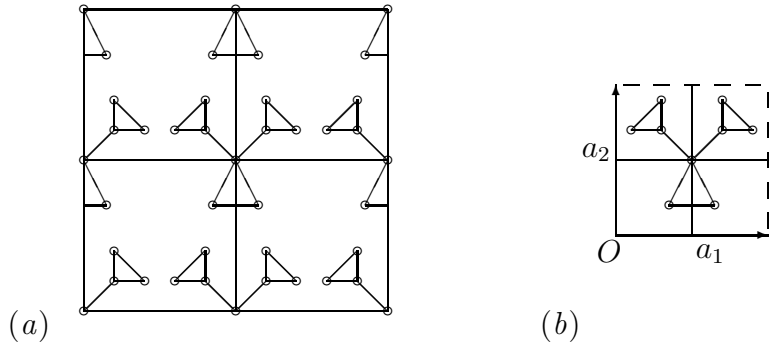
\begin{figure}[h]
\centering
\unitlength 1mm % = 2.845pt
\linethickness{0.4pt}
\ifx\plotpoint\undefined\newsavebox{\plotpoint}\fi % GNUPLOT compatibility
\begin{picture}(120,50)(0,0)

\put(90,30){\circle{1}}
\put(80,20){\vector(1,0){20.00}}
\put(80,20){\vector(0,1){20.00}}
\multiput(81,40)(4,0){5}{\line(1,0){2}}
\multiput(100,21)(0,4){5}{\line(0,1){2}}
\put(77.5,16.5){$O$}
\put(90.5,17.0){$a_1$}
\put(75.5,30.5){$a_2$}
\put(80,30){\line(1,0){20.00}}
\put(90,20){\line(0,1){20.00}}

\put(90,30){\line(1,1){4.00}}
\put(94,34){\line(1,0){4.00}}
\put(94,34){\line(0,1){4.00}}
\put(94,38){\line(1,-1){4.00}}
\put(94,34){\circle{1}}
\put(94,38){\circle{1}}
\put(98,34){\circle{1}}

\put(90,30){\line(-1,1){4.00}}
\put(86,34){\line(-1,0){4.00}}
\put(86,34){\line(0,1){4.00}}
\put(86,38){\line(-1,-1){4.00}}
\put(86,34){\circle{1}}
\put(86,38){\circle{1}}
\put(82,34){\circle{1}}

\put(90,30){\line(-1,-2){3.00}}
\put(90,30){\line(1,-2){3.00}}
\put(87,24){\line(1,0){6.00}}
\put(87,24){\circle{1}}
\put(93,24){\circle{1}}

\put(0,7){(\emph{a})}
%**********************************
\put(70,7){(\emph{b})}

\put(10,10){\line(1,0){40.00}}
\put(10,30){\line(1,0){40.00}}
\put(10,50){\line(1,0){40.00}}
\put(10,10){\line(0,1){40.00}}
\put(30,10){\line(0,1){40.00}}
\put(50,10){\line(0,1){40.00}}

\put(10,10){\circle{1}}
\put(30,10){\circle{1}}
\put(50,10){\circle{1}}

\put(10,30){\circle{1}}
\put(30,30){\circle{1}}
\put(50,30){\circle{1}}

\put(10,50){\circle{1}}
\put(30,50){\circle{1}}
\put(50,50){\circle{1}}
\put(30,30){\line(1,1){4.00}}
\put(34,34){\line(1,0){4.00}}
\put(34,34){\line(0,1){4.00}}
\put(34,38){\line(1,-1){4.00}}
\put(34,34){\circle{1}}
\put(34,38){\circle{1}}
\put(38,34){\circle{1}}

\put(30,30){\line(-1,1){4.00}}
\put(26,34){\line(-1,0){4.00}}
\put(26,34){\line(0,1){4.00}}
\put(26,38){\line(-1,-1){4.00}}
\put(26,34){\circle{1}}
\put(26,38){\circle{1}}
\put(22,34){\circle{1}}

\put(30,30){\line(-1,-2){3.00}}
\put(30,30){\line(1,-2){3.00}}
\put(27,24){\line(1,0){6.00}}
\put(27,24){\circle{1}}
\put(33,24){\circle{1}}
%************************

\put(50,30){\line(-1,1){4.00}}
\put(46,34){\line(-1,0){4.00}}
\put(46,34){\line(0,1){4.00}}
\put(46,38){\line(-1,-1){4.00}}
\put(46,34){\circle{1}}
\put(46,38){\circle{1}}
\put(42,34){\circle{1}}

\put(50,30){\line(-1,-2){3.00}}
\put(47,24){\line(1,0){3.00}}
\put(47,24){\circle{1}}
%************************
\put(10,30){\line(1,1){4.00}}
\put(14,34){\line(1,0){4.00}}
\put(14,34){\line(0,1){4.00}}
\put(14,38){\line(1,-1){4.00}}
\put(14,34){\circle{1}}
\put(14,38){\circle{1}}
\put(18,34){\circle{1}}

\put(10,30){\line(1,-2){3.00}}
\put(10,24){\line(1,0){3.00}}
\put(13,24){\circle{1}}
%***************************
\put(30,10){\line(1,1){4.00}}
\put(34,14){\line(1,0){4.00}}
\put(34,14){\line(0,1){4.00}}
\put(34,18){\line(1,-1){4.00}}
\put(34,14){\circle{1}}
\put(34,18){\circle{1}}
\put(38,14){\circle{1}}

\put(30,10){\line(-1,1){4.00}}
\put(26,14){\line(-1,0){4.00}}
\put(26,14){\line(0,1){4.00}}
\put(26,18){\line(-1,-1){4.00}}
\put(26,14){\circle{1}}
\put(26,18){\circle{1}}
\put(22,14){\circle{1}}

%************************
\put(50,10){\line(-1,1){4.00}}
\put(46,14){\line(-1,0){4.00}}
\put(46,14){\line(0,1){4.00}}
\put(46,18){\line(-1,-1){4.00}}
\put(46,14){\circle{1}}
\put(46,18){\circle{1}}
\put(42,14){\circle{1}}

%************************
\put(10,10){\line(1,1){4.00}}
\put(14,14){\line(1,0){4.00}}
\put(14,14){\line(0,1){4.00}}
\put(14,18){\line(1,-1){4.00}}
\put(14,14){\circle{1}}
\put(14,18){\circle{1}}
\put(18,14){\circle{1}}
%***********************
\put(30,50){\line(-1,-2){3.00}}
\put(30,50){\line(1,-2){3.00}}
\put(27,44){\line(1,0){6.00}}
\put(27,44){\circle{1}}
\put(33,44){\circle{1}}
%************************
\put(50,50){\line(-1,-2){3.00}}
\put(47,44){\line(1,0){3.00}}
\put(47,44){\circle{1}}
%************************
\put(10,50){\line(1,-2){3.00}}
\put(10,44){\line(1,0){3.00}}
\put(13,44){\circle{1}}
\end{picture}
\vspace{-0.5cm}
\caption{\footnotesize  \emph{a}) Precise loop graph $\G_*$;\quad \emph{b}) the fundamental graph $\G_{*f}$.} \label{ffS'}
\end{figure}

{\bf Remark.}
Applying this procedure to the obtained loop graph
$\G_*$  and to another connected finite graph $\G_1$ we obtain  a new
loop graph $\G_{**}$ and so on. Thus, from one loop graph we obtain
a whole class of loop graphs.

\medskip

We now describe bands for precise loop periodic graphs.

\begin{theorem}\lb{T100}
i) Let the Schr\"odinger operator $H=\D+Q$ act on a loop graph $\G$. Then spectral bands $\s_n=\s_n(H)=[\l_n^-,\l_n^+]$ satisfy
\[
\lb{eq.5} \l_n^-=\l_n(0),\qqq \forall \;n\in \N_\n.
\]

ii) Let, in addition, $\G$ be  precise with a  precise
quasimomentum $\vt_0\in\T^d$. Then
\[
\lb{es}
\s_n=[\l_n^-,\l_n^+]=[\l_n(0),\l_n(\vt_0)],\qqq \forall\,n\in \N_\n,
\]
\[
\lb{es.1}
\sum_{n=1}^\n|\s_n|=2\b,
\]
where $\b$ is the number of bridge-loops on the loop fundamental graph
$\G_f$. In particular, if all bridges of $\G_f$ have the form
$(v_k,v_k)$ for some vertex $v_k\in V_f$, then
\[
\lb{es.111}
|\s(H)|=\sum_{n=1}^{\n}|\s_n|=2\b.
\]

\end{theorem}

\no\textbf{Remark.} 1) Due to \er{es.1}, the total length of all spectral bands of the Schr\"odinger operators $H=\D+Q$ on precise loop graphs does not depend on the potential $Q$.

2) The number of the loop fundamental graph bridges can be any integer, then
due to \er{es.111} the Lebesgue measure $|\s(H)|$ of the spectrum of  $H$
(on the specific graphs) can be  greater than any number.

3) $\l_n^-$, $n\in\N_\n$, are the eigenvalues  of the Schr\"odinger  operator
$H(0)$ defined by  (\ref{Hvt}), (\ref{l2.15})  on the fundamental graph $\G_f$. The
identities \er{es} are similar to the case of $N$-periodic Jacobi
 matrices on the lattice $\Z$ (and for Hill operators).
The spectrum of these operators is absolutely continuous and is a
union of spectral bands, separated by gaps. The endpoints of the
bands are the so-called $2N$-periodic eigenvalues.

\medskip

\begin{proposition}\lb{Tfb}
Let $\n,d\geq2$. Then there exists a $\Z^d$-periodic graph $\G$, such that the spectrum of $\D$ on $\G$ has exactly 2 open separated spectral bands $\s_1(\D)$ and $\s_\n(\D)$ and between them, in the gap,
  $\n-2$ degenerated spectral bands (flat bands) $\s_2(\D)=\ldots=\s_{\n-1}(\D)$.
\end{proposition}

{\bf Remark.} There is an open problem: does there exist a $\Z^d$-periodic graph with any $\n\geq2$ vertices in the fundamental graph such that the spectrum of the Laplacian on $\G$ has only 1 spectral band and $\n-1$ flat bands, counting multiplicity?

%**********************************************

%**********************************************
\subsection{Crystal models.}
It is known that the majority of common metals have either a face center cubic (FCC) structure (Fig.\ref{ff.FCC}), a body centered cubic (BCC) structure (Fig.\ref{ff.11}) or a hexagonal close packed (HCP) structure (see \cite{BM80}). The differences between these structures lead to different physical properties of bulk metals. For example, FCC metals, Cu, Au, Ag, are usually soft and ductile, which means they can be
bent and shaped easily. BCC metals are less ductile but stronger, for example iron, while HCP metals are usually brittle. Zinc is HCP and is difficult to bend without breaking, unlike copper.
These structures are obtained from the cubic lattice $\dL^3$
by adding vertices and edges.
In Section \ref{Sec7}  we consider the Schr\"odinger operator
on the hexagonal lattice, on the face-centered cubic lattice and on the body-centered cubic lattice.

\

\section{\lb{Sec3} Direct integrals for Schr\"odinger operators}
\setcounter{equation}{0}

In Proposition \ref{pro0} we present properties of periodic
graphs, needed to prove main results.  We omit the proof, since the
proof  repeats the case of $\Z^2$-periodic graph from \cite{BKS13}.

\begin{proposition}\label{pro0}
i) Each fundamental graph of a $\Z^d$-periodic graph is finite.

ii) Let $(u,v)\in\cA$ and let $m,n\in\Z^d$. Then the following
statements  hold true:

a) $\t(u,v)=-\t(v,u)$.

b) If $(v+m,v+n)\in\cA$, then the index $\t(v+m,v+n)=n-m$.

c) The edge $(u+m,v+m)\in\cA$ and its index $\t(u+m,v+m)=\t(u,v)$.

iii) Let $\t^{(1)}(\tilde\be)$ be the index of an edge
$\tilde\be=(u,v)\in\cA_f$ in the coordinate system with an origin
$O_1$.  Let $[v]\in \Z^d$ denote the integer part of $v$. Then (see
Fig.\ref{gg.00})
\begin{equation}\label{ind'}
\t^{(1)}(\tilde\be)=\t(\tilde\be)+[v-b]-[u-b],\qqq\textrm{ where }\;
b=\overrightarrow{OO}_1.
\end{equation}
\end{proposition}

\setlength{\unitlength}{1.0mm}
\begin{figure}[h]
\centering
\unitlength 1mm % = 2.845pt
\linethickness{0.4pt}
\ifx\plotpoint\undefined\newsavebox{\plotpoint}\fi % GNUPLOT compatibility
\begin{picture}(60,30)(0,0)

\put(10,5){\circle*{1}}
\put(7.5,1.0){$O$}
\put(20.0,2.0){$a_1$}
\put(9.0,12.5){$a_2$}
\put(10,5){\vector(1,0){20.00}}
\put(10,5){\vector(1,2){7.00}}

\put(10,5){\vector(3,1){24.00}}

\put(34,13){\circle*{1}}
\put(32.0,8.5){$O_1$}
\put(44.0,10.0){$a_1$}
\put(33.0,20.5){$a_2$}
\put(34,13){\vector(1,0){20.00}}
\put(34,13){\vector(1,2){7.00}}

\put(23.0,11.0){$b$}

\put(70,8){\line(-1,1){18.00}}
\put(70,8){\circle{1}}

\put(52,26){\circle{1}}
\put(50,28.0){$v+n$ ($v+n-b$)}
\put(61.0,18.0){$\be$}

\put(69,4.0){$u+m$ ($u+m-b$)}
\end{picture}

\caption{\footnotesize  The edge $\be=(u+m,v+n)$ in the coordinate system with the origin~$O$; $\be=(u+m-b,v+n-b)$ in the coordinate system with the origin $O_1$; $a_1,a_2$ are the periods of the graph.} \label{gg.00}
\end{figure}
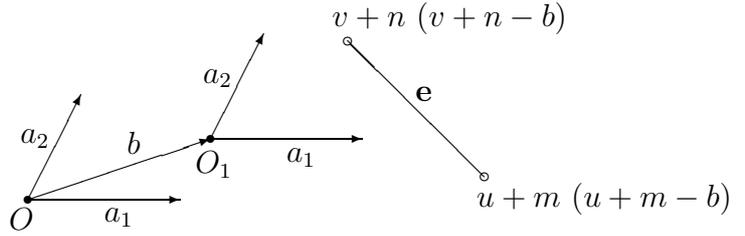

\no \textbf{Remark.}
1) From  iii) it follows that an edge index, generally
speaking, depends on the choice of the coordinate origin.

2) Item ii.b shows that the index of the edge
$(v+m,v+n)\in\cA$ does not depend on the choice of the coordinate origin $O$. It also means that the indices of all loops
on the fundamental graph do not depend on the choice of the point
$O$.

3) Under the group $\Z^d$ action the set $\cA$ of oriented edges of
the graph  $\G$ is divided into equivalence classes. Each
equivalence class is an oriented edge $\be\in\cA_f$ of the
fundamental graph $\G_f$. From item ii.c it follows that all edges
from one equivalence class $\mathbf{e}$ have the same index that is
also the index of the fundamental graph edge $\mathbf{e}$.

\

Now we discuss definitions of   other  Laplacians \cite{MW89}
which are used in the literature.

Firstly,  we consider the so-called \emph{adjacency operator} $\bA$
acting on $\ell^2(V)$ and given  by
\[
\lb{DAO}
 \big(\bA f\big)(v)= \sum\limits_{(v,\,u)_e\in\cE}f(u),\qqq
 v\in V,
\]
where all loops in the sum are counted twice. Note that  the
adjacency operator $\bA$ satisfies:
\[
\lb{Oid}
\D=\vk-\bA,\qqq
\s(\bA)\ss [-\vk_+,\vk_+],
\]
where $\vk$ is the  degree  operator defined by $\big(\vk
f\big)(v)=\vk_vf(v)$ for all $f\in \ell^2(V)$.

Secondly, there exists the \emph{normalized Laplace operator}
$\D_\ast$ acting on $\ell^2(V)$ and given  by
\[
\lb{DAvO}
\D_\ast=\1-\vk^{-{1\/2}}\bA\vk^{-{1\/2}},
\]
where $\1$ is the identity operator.

Thirdly, we introduce the Hilbert space  $\ell^2(\cA)$ of all square
summable functions $\phi:\cA\to \C$ such that $\phi(\be)=-\phi(\bar
\be)$ for all $\be\in\cA$, equipped with the norm
$$
\|\phi\|^2_{\ell^2(\cA)}=\textstyle\frac12\sum\limits_{\be\in\cA}|\phi(\be)|^2<\infty.
$$
Here $\bar\be$ is the inverse edge of $\be$.
We define the operator $\nabla_\G:\ell^2(V)\rightarrow\ell^2(\cA)$
by
$$
(\nabla_\G f)(\be)=f(v)-f(u),\qqq \forall f\in \ell^2(V),\qqq
\textrm{ where } \ \be=(u,v).
$$
The  conjugate operator $\nabla_\G^\ast:\ell^2(\cA)\to \ell^2(V)$
has the form
$$
(\nabla_\G^\ast \phi)(v)=-\sum_{\be=(v,u)\in\cA}\phi(\be),\qqq
\forall \phi\in \ell^2(\cA).
$$
Then the Laplacian $\D$ satisfies
$$
\D=\nabla_\G^\ast \nabla_\G.
$$

These operators are used in different applications. It is known
\cite{C97} that the investigation of the spectrum of Laplacians and
Schr\"odinger operators on an equilateral metric graph (i.e., a
graph consisting of identical segments) can be reduced to the study
of the spectrum of the discrete normalized Laplacian $\D_*$, see
\er{DAvO}. In \cite{KS} we describe spectral properties (including
the Bethe-Sommerfeld conjecture) of Laplace operators on
quantum graphs. In order to do this we need to study the spectrum of
the normalized Laplacian $\D_\ast$. The spectrum of the perturbed
normalized Laplacians \er{DAvO} on $\Z^d$-periodic graphs is studied
in our another paper \cite{KS1}.

%%%%%%%%%%%%%%%%%%%%%%

\

We begin to discuss the Laplacian $\D$ given by \er{l2.15}.
We rewrite the Floquet matrix $\D(\vt)$ for the Laplacian $\D$
 in terms of  a degree matrix $\vk_0$ and an adjacency  matrix $\bA(\vt)$ by
\[
\lb{MAdj}
\D(\vt)=\vk_0-\bA(\vt), \qqq   \vk_0=\diag(\vk_1,\ldots,\vk_\n),\qqq \vt\in\T^d.
\]
 From \er{MAdj} and \er{l2.15} it follows that the matrix $\bA(\vt)$ is given by
\[
\lb{Adj}
\bA(\vt)=\{\bA_{jk}(\vt)\}_{j,k=1}^\n,\qqq
\bA_{jk}(\vt)=\ca
\sum\limits_{{\bf e}=(v_j,\,v_k)\in{\cA}_f}e^{\,i\lan\t
({\bf e}),\,\vt\ran }, \qq &  {\rm if}\  \ (v_j,v_k)\in \cA_f \\
\qqq 0, &  {\rm if}\  \ (v_j,v_k)\notin \cA_f \ac.
\]
From \er{Hvt} and \er{MAdj} it follows that
\[
\lb{HHAA} H(\vt)=\vk_0-\bA(\vt)+q.
\]

\

\no {\bf Proof of Theorem \ref{pro2}.i -- iii.} i) In the proof we use some
arguments from \cite{BKS13} about a normalized Laplacian on $\Z^2$ -periodic graphs, but we need additional ones, since we consider Schr\"odinger operators $H=\D+Q$ and the Laplacian $\D$ is not normalized.

We introduce the Hilbert space (a constant fiber direct integral) $\mH=L^2\big(\T^d,{d\vt
\/(2\pi)^d}\,,\cH\big)=\int_{\T^d}^{\os}\cH\,{d\vt \/(2\pi)^d}$, where $\cH=\C^\nu$, equipped with the norm
$$
\|g\|^2_{\mH}=\int_{\T^d}\|g(\vt )\|_{\C^\nu}^2\frac{d\vt
}{(2\pi)^d}\,,
$$
where a function $g:\T^d\rightarrow\C^\nu$.

%\ell_{fin}^2(V)\ell^2(V)

Denote by $\ell_{fin}^2(V)$ the set of all finitely supported functions $f\in
\ell^2(V)$. Recall that the vertices of $\G_f$ are identified with
the vertices $v_1,\ldots,v_\nu$ of the periodic graph $\G$ from the set $[0,1)^d$. Let $U:\ell^2(V)\to\mH$ be the operator defined by
\[
 \lb{5001}
(Uf)_n(\vt )=\sum\limits_{m\in\mathbb{Z}^d}e^{-i\lan m,\vt\ran }
f(v_n+m), \qqq
(\vt,n)\in \T^d\ts\N_\nu.
\]
Standard arguments (see pp. 290--291 in \cite{RS78}) give that $U$ is well defined on $\ell_{fin}^2(V)$ and has an unique extension to a unitary operator.
For  $f\in \ell_{fin}^2(V)$ the sum \er{5001} is finite and using the identity $V=\big\{v_n+m: (n,m)\in
\N_\nu\ts\Z^d\big\}$  we have
$$
\begin{aligned}
&\|Uf\|^2_{\mH}=\int_{\T^d}\|(Uf)(\vt )\|_{\C^\n}^2{d\vt\/(2\pi)^d}\\
&=
 \int_{\T^d}\sum_{n=1}^\nu\bigg(\sum\limits_{m\in\Z^d}e^{-i\lan m,\vt\ran }
 f(v_n+m)\bigg)
 \bigg(\sum\limits_{m'\in\Z^d}e^{i\lan m',\vt\ran }\bar{f}(v_n+m')\bigg)
{d\vt \/(2\pi)^d}\\
&=\sum_{n=1}^\nu\sum_{m,m'\in\Z^d}\lt(
 f(v_n+m)\bar f(v_n+m')\int_{\T^d}e^{-i\lan m-m',\vt\ran }
 {d\vt \/(2\pi)^d}\rt)\\
 &=
 \sum_{n=1}^\nu\sum_{m\in\Z^d}\big|f(v_n+m)\big|^2
 =\sum_{v\in V}|f(v)|^2=\|f\|_{\ell^2(V)}^2.
\end{aligned}
$$
 Thus, $U$ is  well defined on $\ell_{fin}^2(V)$ and has a unique
isometric extension. In order to prove  that $U$ is onto $\mH$ we compute
$U^*$.   Let
$g=\big(g_n\big)_{n\in\N_\nu}\in\mH$, where $g_n
:\T^d\to \C$. We define
\begin{equation}\label{rep}
(U^*g)(v)=\int_{\T^d}e^{i\lan m,\vt\ran }g_n(\vt ){d\vt \/(2\pi)^d}\,,\qquad
 v=v_n+m\in V, \qq (n,m)\in \N_\nu\ts\Z^d,
\end{equation}
where $(n,m)\in \N_\nu\ts\Z^d$ are uniquely defined.  A direct computation
gives  that it is indeed the formula for the adjoint of $U$.
Moreover, the Parseval's identity for the Fourier series gives
$$
\|U^*g\|^2_{\ell^2(V)}=\sum_{v\in
V}\big|(U^*g)(v)\big|^2=\sum_{n=1}^\nu\sum_{m\in\Z^d}\big|(U^*g)(v_n+m)\big|^2=
\sum_{n=1}^\nu\sum_{m\in\Z^d}\bigg|\int_{\T^d}e^{i\lan m,\vt\ran
}g_n(\vt ){d\vt \/(2\pi)^d}\bigg|^2
$$
$$
=\sum_{n=1}^\nu\int_{\T^d}\big|g_n(\vt)\big|^2{d\vt \/(2\pi)^d}=
 \int_{\T^d}\sum_{n=1}^\nu\big|g_n(\vt)\big|^2{d\vt \/(2\pi)^d}=\|g\|_\mH^2,
$$
Recall that  $\bA_{jk}(\vt)$ is defined by \er{Adj}. Then for $f\in \ell_{fin}^2(V)$ and $j\in\N_\nu$ we obtain
\[
\begin{aligned}
\label{ext1}
(U\bA f)_j(\vt )=\sum_{m\in\Z^d}e^{-i\lan m,\vt\ran}(\bA
f)(v_j+m)=\sum_{m\in\Z^d}e^{-i\lan m,\vt\ran }
\sum\limits_{(v_j+m,\,u)_e\in\cE}f(u)\\=\sum_{m\in\Z^d}e^{-i\lan m,\vt\ran
}
\sum\limits_{k=1}^\nu\sum\limits_{\be=
(v_j,v_k)\in\cA_f}f\big(v_k+k+\t (\be)\big)\\
=\sum\limits_{k=1}^\nu\sum\limits_{\be=
(v_j,v_k)\in\cA_f}\,e^{i\lan\t(\be),\vt\ran }
\sum_{m\in\Z^d}e^{-i\lan m+\t(\be),\vt\ran }f\big(v_k+m+\t (\mathbf{e})\big)\\
=\sum\limits_{k=1}^\nu\sum\limits_{\be= (v_j,v_k)\in\cA_f}e^{i\lan\t
(\be),\vt\ran } (Uf)_k(\vt)=\sum\limits_{k=1}^\nu
\bA_{jk}(\vt)(Uf)_k(\vt),
\end{aligned}
\]
and the operator $\vk+Q$ satisfies
\[
\begin{aligned}\label{ext1'}
\big(U(\vk+Q)f\big)_j(\vt)=\sum_{m\in\Z^d}e^{-i\lan m,\vt\ran}\big((\vk+Q)
f\big)(v_j+m)\\=\sum_{m\in\Z^d}e^{-i\lan m,\vt\ran }(\vk_j+q_j)f(v_j+m)=(\vk_j+q_j)(Uf)_j(\vt).
\end{aligned}
\]
For the matrices $\bA(\vt)$ and $\vk_0$ defined by \er{MAdj}, the identities (\ref{ext1}), (\ref{ext1'}) yield
$$
(U\bA f)(\vt
)=\bA(\vt)(Uf)(\vt ),\qqq \big(U(\vk+Q)f\big)(\vt
)=(\vk_0+q)(Uf)(\vt).
$$
%where the matrices $\bA(\vt)$ and $\vk_0$ are defined by \er{MAdj}.
Thus, due to $\D(\vt)=\vk_0-\bA(\vt)$, we obtain
\begin{multline*}
UHU^{-1}=U(\D+Q)U^{-1}=U(\vk+Q-\bA)U^{-1}\\=
\frac1{(2\pi)^d}\int_{\T^d}^{\os}\big(\vk_0+q-\bA(\vt)\big)\,d\vt=
\frac1{(2\pi)^d}\int_{\T^d}^{\os}\big(\D(\vt)+q\big)\,d\vt=
\frac1{(2\pi)^d}\int_{\T^d}^{\os}H(\vt)\,d\vt,
\end{multline*}
which gives the proof of i).

ii) Recall that the vertices $v_1,\ldots,v_\nu$ of the  fundamental graph
$\G_f$ are
identified with the vertices of the periodic graph $\G$ from the set
$[0,1)^d$ in the coordinate system with the origin $O$. Due to Proposition
\ref{pro0}.iii, we have that for each $(v_j,v_k)\in\cA_f$
\begin{equation}\label{ind}
\t^{(1)}(v_j,v_k)=\t(v_j,v_k)+m_k-m_j,\qqq m_k=[v_k-b],\qqq m_j=[v_j-b],
\end{equation}
where $\t^{(1)}(v_j,v_k)$ is the index of the edge $(v_j,v_k)$ in
the coordinate system with the origin $O_1$,
$b=\overrightarrow{OO}_1$. Let $\bA^{(1)}(\vt)$ be the matrix
defined by (\ref{MAdj}), (\ref{Adj}) in the coordinate system with
the origin $O_1$. Applying (\ref{ind}), we obtain the following form for the entries of
$\bA^{(1)}(\vt)$
\[
\lb{AA1} \bA_{jk}^{(1)}(\vt
)=\hspace{-2mm}\sum\limits_{\be=(v_j,\,v_k)\in\cA_f}\hspace{-2mm}
e^{i\lan\t^{(1)} (\be),\vt\ran }=e^{i\lan m_k-m_j,\vt\ran
}\hspace{-2mm}
\sum\limits_{\be=(v_j,\,v_k)\in\cA_f}\hspace{-2mm}e^{i\lan\t
(\be),\,\vt\ran }=e^{i\lan m_k-m_j,\,\vt\ran }\bA_{jk}(\vt ).
\]
We introduce the diagonal $\nu\times\nu$ matrix
$$
\cU(\vt )=\mathrm{diag}  \left(
\begin{array}{ccc}
   e^{-i\lan m_1,\,\vt\ran}, &  \ldots,& e^{-i\lan m_\nu,\,\vt\ran}
\end{array}\right),\qquad \forall\, \vt \in\T^d.
$$
Using \er{AA1}, we obtain
\[
\lb{AAA2}
\cU(\vt )\,\bA(\vt )\,\cU^{-1}(\vt
)=\bA^{(1)}(\vt),\qq \forall\, \vt \in\T^d.
\]
Since the matrices $\cU(\vt)$ and $\vk_0+q$ are diagonal, from \er{AAA2} we deduce that for each $\vt\in\T^d$
the matrices $H(\vt)=\vk_0+q-\bA(\vt)$ and $H^{(1)}(\vt)=\vk_0+q-\bA^{(1)}(\vt)$ are unitarily equivalent.

iii) This statement is a direct consequence of (\ref{l2.15}) and  the
definition of a bridge.

\medskip

%*****************************************************************************
In order to prove  iv), we need to discuss the following properties.

\begin{proposition}
\label{pp1} i) The point 0 is never a flat band of the Laplacian $\D$.

ii) The matrix $H(0)$ is the Schr\"odinger operator on the
fundamental  graph $\G_f$ given by
\[
\label{ll6}
H(0)=\D(0)+q,\qqq
\D_{jk}(0)=\vk_j\d_{jk}-\vk_{jk}\,,\qq \forall\,(j,k)\in\N_\n^2,
\]
where $\vk_{jk}\ge 1$ is the multiplicity of the edge $(v_j,v_k)$,
if $(v_j,v_k)\in \cA_f$ and  $\vk_{jk}=0$, if $(v_j,v_k)\notin
\cA_f$. Moreover, they satisfy
\begin{equation}
\label{l.83}
\vk_j=\sum_{k=1}^\n\vk_{jk}\ge1,\qquad \forall\,j\in \N_\nu.
\end{equation}

\end{proposition}

\no \textbf{Proof.}
i) The proof is by contradiction. Let the point 0 be an eigenvalue
of the Laplacian $\D$ on a graph $\G$. Then there exists an
eigenfunction  $0\neq f\in\ell^2(V)$ with the eigenvalue 0 and with
a finite support $\mB\ss V$ (see Theorem 3.2 in \cite{HN09}). Let
$\max\limits_{v\in\mB}f(v)=f(\tilde v)$ for some $\tilde v\in\mB$.
Thus, we obtain
$$
0=\big(\D f\big)(\tilde v)=\vk_{\tilde v}\,f(\tilde v)-
\sum\limits_{(\tilde v,\,u)_e\in\cE}f(u)\geq\vk_{\tilde v}\,f(\tilde v)-
\sum\limits_{(\tilde v,\,u)_e\in\cE}f(\tilde v)=0.
$$
From this we deduce that the inequality has to be an equality, and therefore
$
f(u)=f(\tilde v), \ \forall\,u\sim
\tilde v.
$
Repeating this argument until we reach a vertex from
$V\setminus\mB$, we conclude that $f=0$. We obtain a contradiction.
Thus, the point 0 is never a flat band of $\D$.

ii)  A direct calculation of the
entries of the matrices
$H(0)$ and $\D(0)$ using the formulas \er{Hvt}, (\ref{l2.15}) gives the identities (\ref{ll6}). From the
definitions \er{DLO}, \er{Sh}, \er{Pot} it follows that the  Schr\"odinger
operator on $\G_f$ has the form \er{ll6}. Since the degree
${\vk}_j$ of the vertex $v_j$ is equal to the number of oriented
edges starting at $v_j$, we obtain (\ref{l.83}). From connectivity of the graph $\G$ it follows that $\vk_j\geq1$.
 \quad $\BBox$

\medskip

\no {\bf Proof of Theorem \ref{pro2}.iv.}
The proof is by contradiction. Assume that all entries $\D_{jk}(\cdot)$,
$1\leq j\leq k\leq\n$, are constant. Since the matrix $\D(\vt)$ is self-adjoint,  all its entries are constant and $\D(\cdot)=\D(0)$.
Using the fact that $\D(0)$ is the Laplacian on the fundamental graph $\G_f$,
which has the eigenvalue $0$, we deduce that
the point 0 is an eigenvalue of $\D$ with
infinite multiplicity. This contradicts Proposition \ref{pp1}.i.
\quad $\BBox$

\medskip

\section{\lb{Sec4} Spectral estimates in terms of geometric parameters of the graph}
\setcounter{equation}{0}

Below we need the following representation of the Floquet matrix
$H(\vt),\vt\in\T^d$:
\[
\label{eq.1}
H(\cdot)=H_0+\wt\D(\cdot),\qqq
H_0={1\/(2\pi)^d}\int_{\T^d}H(\vt )d\vt.
\]
From \er{eq.1}, \er{Hvt}, \er{l2.15} we deduce that the matrix
$\wt\D(\cdot)$ has the form
\[
\label{tl2.15} \wt\D(\vt)=\{\wt\D_{jk}(\vt)\}_{j,k=1}^\n,\qqq
\wt\D_{jk}(\vt )= -\sum\limits_{{\bf
e}=(v_j,v_k)\in{\cB}_f}e^{\,i\lan\t ({\bf e}),\,\vt\ran }.
\]

\begin{proposition}\lb{LP}
The lower point of the spectrum $\s(H)$ of the operator $H$ is
$\l_1(0)=\l_1^-$.
\end{proposition}
\no {\bf Proof.} Recall that Sy and Sunada prove this result for a more general class of graphs (see Theorem~1, p.143 in \cite{SS92}). Their proof is rather complicated. For readers' convenience we give here a simple proof, based on matrix properties. We define $\n\ts\n$ matrix
\[
\lb{MaK}
K(\vt)=\a \1_\n-H(\vt)=\bA(\vt)-(\vk_0+q)+\a \1_\n,\qqq
\a=\max_{j\in\N_\n} (\vk_j+q_j),
\]
where $\bA(\vt)$ is defined by \er{MAdj}, \er{Adj}. We arrange
the eigenvalues of the matrix $K(\vt)$ in increasing order
$\z_1(\vt)\leq\ldots\leq\z_\n(\vt)$.
Since $\bA(0)$ is an adjacency matrix of the connected graph $\G_f$,
due to Proposition \ref{MP}.vi the matrix $\bA(0)$ is irreducible.
Then the matrix $K(0)$, that differs from $\bA(0)$ in the diagonal
entries only, is also irreducible. From \er{MaK} it follows that all
entries of the matrix $K(0)$ are nonnegative. Then Proposition
\ref{MP}.vii  implies that the spectral radius $\r\big(K(0)\big)$ is
a simple eigenvalue of $K(0)$, i.e.,
$\r\big(K(0)\big)=\z_\n(0)$. From the formulas \er{MaK},
(\ref{Adj}) it follows that the entries of $K(\vt)$ satisfy
\[
\lb{vv1}
|K_{jk}(\vt)|\leq K_{jk}(0), \qquad \forall\,(j,k, \vt)\in\N_\nu^2\ts \T^d.
\]
Then, Proposition \ref{MP}.i  implies that the spectral radius  satisfies
$
\r\big(K(\vt)\big)\leq\r\big(K(0)\big)=\z_\n(0),
$
which yields
$$
\a-\l_{\n-n+1}(\vt)=\z_n(\vt)\leq\z_\n(0)
=\a-\l_1(0),\qqq \forall\,(n,\vt)\in\N_\nu\ts\T^d.
$$
Then
$
\l_1(0)\leq\l_n(\vt),\ \forall\,(n,\vt)\in\N_\nu\ts\T^d,
$
which gives that $\l_1(0)=\l_1^-$. \qq $\BBox$

\

%*****************************************
\no {\bf Proof of Theorem \ref{T1}.} i) The proof for the case $Q=0$  has been given in Proposition \ref{pp1}.i.

Let now $|\D_{jk}(\cdot)|\neq\const$ for some $j,k\in\N_\n$. Then for the entry $K_{jk}(\cdot)$ of the matrix $K$, defined by \er{MaK}, we have
$|K_{jk}(\cdot)|\neq\const$. This, \er{vv1} and the irreducibility
of the matrix $K(\vt)$ according to Proposition \ref{MP}.i--ii
imply that  the spectral radius $\r\big(K(\vt)\big)$ satisfies
$$
\r\big(K(\vt)\big)<\r\big(K(0)\big)=\z_\n(0)=\a-\l_1(0)
$$
for almost all $\vt\in\T^d$. This yields
$$
\a-\l_1(\vt)=\z_\n(\vt)\leq\r\big(K(\vt)\big)<
\a-\l_1(0),
$$
i.e., $\l_1(0)<\l_1(\vt)$ for almost all $\vt\in\T^d$. Thus,  the
first spectral band of $H$ is non-degenerated.

ii)  Define the diagonal operator $B(\vt)$ acting on $\C^\n$ by
\[
\lb{eq.2'}
B(\vt)=\diag(B_1,\ldots,B_\n)(\vt),\qqq
B_j(\vt )=\sum_{k=1}^\n |\wt\D_{jk}(\vt)|,\qq \vt\in\T^d.
\]
From \er{tl2.15} we deduce that
\[
\lb{wvv1} |\wt\D_{jk}(\vt )|\leq|\wt\D_{jk}(0)|=\b_{jk}, \qqq
\forall \ (j,k,\vt)\in\N_\nu^2\ts\T^d,
\]
where $\b_{jk}$ is the number of bridges $(v_j,v_k)$ on $\G_f$. Then \er{wvv1} gives
\[
\lb{B0} B(\vt)\leq B(0),\qqq \forall\,\vt\in\T^d.
\]
Then estimate \er{B0} and Proposition \ref{TK}.i yield
\[
\lb{eq.2} -B(0) \leq -B(\vt)\le\wt\D(\vt)\le B(\vt)\leq B(0),\qqq
\forall \vt\in\T^d.
\]

We use some arguments from \cite{BKS13}, \cite{Ku10}. Combining
\er{eq.1} and \er{eq.2}, we obtain
\[
H_0-B(0)\le H(\vt)\le H_0+B(0).
\]
Thus, the standard perturbation theory (see Proposition
\ref{MP}.iii) gives
\[
\l_n(H_0-B(0))\leq\l_n^-\le \l_n(\vt)\leq\l_n^+\le \l_n(H_0+B(0)),
\qqq \forall \ (n,\vt) \in \N_\n\ts\T^d,
\]
which implies
\[
\lb{009}
 \big|\s(H)\big|\le \sum_{n=1}^{\nu}(\l_n^+-\l_n^-)\leq
 \sum_{n=1}^\nu\big(\l_n(H_0+B(0))-\l_n(H_0-B(0))\big)=
 2\Tr B(0).
\]
In order to determine $2\Tr B(0)$ we use the relations   \er{wvv1}
and we obtain
\[
\lb{pro} 2\Tr B(0)=2\sum_{j=1}^\n
B_j(0)=2\sum_{j,k=1}^\n|\wt\D_{jk}(0)|= 2\sum_{j,k=1}^\n\b_{jk}=2\b.
\]
From \er{009} and \er{pro} it follows the estimate \er{eq.7}.

Now we will prove \er{GEga}.
Since $\l_1^-$ and $\l_\n^+$ are the lower and upper points of the spectrum, respectively, using the estimate \er{eq.7}, we obtain
\[\lb{GEga1}
\sum_{n=1}^s|\g_n(H)|=
\l^+_\n-\l_1^--\big|\s(H)\big|\geq\l^+_\n-\l_1^--2\b.
\]

We rewrite the sequence $q_1,\ldots,q_\n$, defined by \er{pott}, in nondecreasing order
\[
\lb{wtqn}
q_1^\bullet\le q_2^\bullet \le\ldots \le q_\n^\bullet \qq {\rm and \ let } \qq q_1^\bullet=0.
 \]
Here $q_1^\bullet=q_{n_1}, q_2^\bullet=q_{n_2},\ldots,q_\n^\bullet=q_{n_\n}$
for some distinct integers $n_1, {n_2},\ldots, {n_\n}\in \N_\n$ and without loss of generality we may assume that $q_1^\bullet=0$.

Then Proposition \ref{MP}.iv and the basic fact $\s(\D)\subset[0,2\vk_+]$
give that the eigenvalues of the Floquet matrix $H(\vt)$ for $H=\D+Q$, satisfy
\[
\lb{qq0}
\begin{array}{l}
  q_n^\bullet\le\l_n(\vt)\le q_n^\bullet+2\vk_+ \\[6pt]
  \l_n^0(\vt)\le\l_n(\vt)\le \l_n^0(\vt)+q_\n^\bullet
\end{array},
\qqq \forall\, (\vt,n)\in\T^d\ts\N_\n.
\]
The first inequalities in \er{qq0} give
\[\lb{qq}
\l_\n^+\geq q_\n^\bullet,\qqq \l_1^-\leq 2\vk_+,
\]
and, using the second inequalities in \er{qq0}, we have
\[\lb{qq1}
\l_\n^{0+}=\max_{\vt\in\T^d}\l_\n^{0}(\vt)=\l_\n^{0}(\vt_+)\leq\l_\n(\vt_+)\leq \l_\n^+,
\]
\[\lb{qq2}
0=\l_1^{0-}=\min_{\vt\in\T^d}\l_1^{0}(\vt)=\l_1^{0}(\vt_-)\geq\l_1(\vt_-)-q_\n^\bullet
\geq\l_1^--q_\n^\bullet
\]
for some $\vt_-,\vt_+\in\T^d$.
From \er{qq} -- \er{qq2} it follows that
$$
\l^+_\n-\l_1^-\geq q_\n^\bullet-2\vk_+,\qqq
\l^+_\n-\l_1^-\geq \l_\n^{0+}-q_\n^\bullet,
$$
which yields \er{GEga}.
The last statement of the theorem will be proved in Theorem \ref{T100}.
\qq $\BBox$

\

{\bf Definition of generic potentials.} {\it
 A potential $Q$ is called  "generic" if
the values of the potential $q_j=Q(v_j)$ on the vertex set
$V_f=\{v_1,\ldots,v_\n\}$ of some fundamental graph $\G_f$ are
distinct, i.e., they satisfy}
\[
\lb{dif} q_j\neq q_k \qq \textrm{ \ for all } \qq j,k\in\N_\n, \qq
j\neq k.
\]

%***************************************************
\no {\bf Proof of Theorem \ref{T17}.} The Floquet matrix $H_t(\vt)$,
$\vt\in\T^d$, for the Schr\"odinger operator $H_t=\D+tQ$, where the
potential $Q$ is generic, has the form
$$
H_t(\vt)=\D(\vt)+tq,\qqq q=\diag(q_1,\ldots,q_\n),\qqq
q_1<q_2<\ldots<q_\n.
$$
We define the matrix
$$
\wh{H}_t(\vt)=\textstyle\frac1t\,H_t(\vt)=q+\ve\D(\vt),\qqq
\ve=\frac1t\,.
$$
Each eigenvalue $\hat{\l}_n(\vt,t)$ of the matrix
$\widehat{H}_t(\vt)$  has the following asymptotics:
$$
\hat{\l}_n(\vt,t)=q_n+\ve\,\D_{nn}(\vt)-\ve^2\sum_{j=1 \atop j\neq
n}^\n\frac{|\D_{jn}(\vt)|^2}{q_j-q_n}+O(\ve^3), \qqq n\in\N_\n,
$$
(see pp. 7--8 in \cite{RS78}) uniformly in $\vt\in\T^d$ as
$t\rightarrow\infty$. This yields the asymptotics of the eigenvalues
$\l_n(\vt,t)$ of the matrix $H_t(\vt)$:
\[
\lb{ass} \l_n(\vt,t)=t\,\hat{\l}_n(\vt,t)=tq_n+\D_{nn}(\vt)-
\ve\sum_{j=1 \atop j\neq
n}^\n\frac{|\D_{jn}(\vt)|^2}{q_j-q_n}+O(\ve^2).
\]
Since $\s_n(H_t)=\l_n(\T^d)$, the asymptotics \er{ass} also gives
that \[ \lb{ass2} |\s_n(H_t)|=|\D_{nn}(\T^d)|+O(1/t).
\]

If there exists a bridge-loop $(v_n,v_n)\in \cB_f$,  then
$\D_{nn}(\cdot)\neq\const$ and the function $\l_n(\vt,t)$ is not
constant. Thus, the asymptotics \er{ass2} gives that the spectral band
$\s_n(H_t)$ of the Schr\"odinger operator $H_t$ is non-degenerated as
$t\rightarrow\infty$ and we have \er{Qt1} and $C>0$.

If there are no bridge-loops on $\G_f$, then \er{ass2}  yields
$|\s(H_t)|=O(1/t)$. \qq $\BBox$

\

{\bf Remark.} We do not know an example of a connected periodic graph,
when the spectrum of the Schr\"odinger operator $H=\D+Q$ has a flat
band for the case of generic potentials $Q$. There is an open
problem to show that   the spectrum of $H$ is absolutely continuous
for generic potentials~$Q$.

%***********************************
\section{\lb{Sec7b} Schr\"odinger operators on bipartite regular graphs}\lb{cm}
\setcounter{equation}{0}
We recall some definitions. A graph is called \emph{bipartite} if its vertex set is divided into two disjoint sets (called \emph{parts} of the graph) such that
each edge connects vertices from distinct sets (see p.105 in \cite{Or62}).
Examples of bipartite graphs are the cubic lattice (Fig.\ref{ff.0.1}\emph{a})
and the hexagonal lattice (Fig.\ref{ff.0.3}\emph{a}). The  face-centered cubic
lattice (Fig.\ref{ff.FCC}\emph{a}) is non-bipartite. A graph is called \emph{regular of degree $\vk_+$} if each its vertex $v$ has
the degree $\vk_v=\vk_+$.

\begin{proposition}
\lb{fbg}
For a bipartite $\Z^d$-periodic graph there exists a bipartite
fundamental graph.
\end{proposition}
The proof  repeats the case of $\Z^2$-periodic graph from \cite{BKS13} and is omitted. Note that not every fundamental graph of a bipartite periodic graph is bipartite. Indeed, the cubic lattice (Fig.\ref{ff.0.1}\emph{a}) is bipartite, but its fundamental graph shown in Fig.\ref{ff.0.1}\emph{b} is non-bipartite.

\

Now we discuss  spectral properties of Laplace and Schr\"odinger operators on bipartite regular graphs.

\begin{theorem} \lb{TBG}
Let a graph $\G$ be regular of degree $\vk_+$. Then the following statements hold true.

i) If $\G$ is bipartite and
symmetric with respect to the coordinate origin and $Q(v)=-Q(-v)$
for all $v\in V$, then the spectrum of the Schr\"odinger operator
$H$ is symmetric with respect to the point $\vk_+$.

ii) A fundamental graph $\G_f$ is
bipartite iff the spectrum of the Floquet matrix $\D(\vt)$ is symmetric with
respect to the point $\vk_+$ for each $\vt \in\T^d$.

iii) If a fundamental graph $\G_f$ is bipartite and the number $\nu$ of vertices in $\G_f$ is odd, then $\vk_+$ is a flat band of $\D$.

iv) If $\G$ is bipartite and in the spectrum of the Laplacian $\D$ there exist $s$ spectral gaps $\g_1(\D),\ldots,\g_s(\D)$, then the following estimate holds true:
\[\lb{GEgaB}
\sum_{n=1}^s|\g_n(\D)|\ge
2(\vk^+-\b).
\]
Here $\b$ is the number of fundamental graph bridges.

v) If $\G$ is a bipartite loop graph ($\G_f$ is non-bipartite, since there is a loop on $\G_f$),
then each spectral band of the Laplacian $\D$ on $\G$ has the form
$\s_n^0=[\l_n^{0-},\l_n^{0+}]$, $n\in \N_\n$, where $\l_n^{0-}$ and
$\l_n^{0+}$ are the eigenvalues of the matrices $\D(0)$ and
$2\vk_+\1_\n-\D(0)$, respectively ($\1_\n$ is the identity $\n\ts\n$
matrix).
\end{theorem}

\no {\bf Proof.} i) Let $\G$ be a bipartite periodic graph with the parts $V_1$ and
$V_2$.  We define the unitary operator $\cU$ on $\ell^2(V)$ by
$$
(\cU f)(v)=\left\{\begin{array}{rl}
f(v),  \ &  {\rm if}\  \ v\in V_1 \\[6pt]
 -f(v), \ &  {\rm if}\  \ v\in V_2  \\
\end{array}\right.,\qqq f\in\ell^2(V),\qq  v\in V.
$$
Then we have
\[
\lb{sym1}
\cU^{-1}=\cU,\qqq \cU\bA\cU^{-1}=-\bA,\qqq \cU Q\cU^{-1}=Q.
\]
Let $\cU_S : \ell^2(V)\to \ell^2(V)$ be the reflection operator
given by $ (\cU_S f)(v)=f(-v), \  v\in V. $ Then $\cU_S$ satisfies
\[\lb{sym2}
\cU_S^{-1}=\cU_S,\qqq \cU_S\bA=\bA\cU_S,\qqq \cU_S Q=-Q\cU_S.
\]
From \er{sym1} and \er{sym2} it follows that
$$
\cU_S\cU(-\bA+Q)\cU^{-1}\cU_S^{-1}=\bA-Q,
$$
which yields that $\s(-\bA+Q)$ is symmetric with respect to 0.
Using the formula $H=\vk_+\1-\bA+Q$, we deduce that $\s(H)$ is symmetric with respect to
$\vk_+$.

ii) Let $\G_f$ be a bipartite fundamental graph with the parts
$V_1$ and $V_2$. The same arguments used in the proof of the previous item show that the spectrum $\s(\D(\vt))$ is symmetric with respect to
$\vk_+$ for each $\vt \in\T^d$.

Conversely, since $\D(0)$ is the Laplacian on the regular graph
$\G_f$ and its spectrum is symmetric with respect to $\vk_+$, the
spectrum of the adjacency operator $\bA(0)=\vk_+\1_\n-\D(0)$ on $\G_f$ is symmetric
with respect to 0. It gives that the graph $\G_f$ is bipartite (see
Theorem 3.11 in \cite{CDS95}).

iii) For each $\vt\in\T^d$ the matrix $\D(\vt)$ has $\nu$ eigenvalues,
where $\nu$ is odd. The previous item gives that the spectrum $\s(\D(\vt))$ is
symmetric with respect to $\vk_+$. Then $\vk_+\in\s\big(\D(\vt )\big)$ for
any $\vt\in\T^d$. Therefore, $\vk_+$ is a flat band of $\D$.

iv) For the Laplacian $\D$ on a bipartite regular graph $\G$ we have
$\l_1^0(0)=0$, $\l_\n^{0+}=2\vk^+$ and the first inequality in \er{GEga} gives \er{GEgaB}.

v) Since $\G$ is bipartite and regular of degree $\vk_+$, there
exists a bipartite fundamental graph $\wt\G_f$, which is also
regular of degree $\vk_+$. Then, due to item i), the
spectrum of the Floquet matrix $\D(\vt)$ corresponding to the graph
$\wt\G_f$ is symmetric with respect to the point $\vk_+$ for each
$\vt\in\T^d$. Therefore, the spectrum of the Laplacian on $\G$ is
also symmetric with respect to $\vk_+$. From the formula \er{eq.5}
it follows that $\l_1^0(0)\leq\ldots\leq\l_\n^0(0)$ are the lower
endpoints of the spectral bands. Then, by the symmetry of the
spectrum, $2\vk_+-\l_\n^0(0)\leq\ldots\leq2\vk_+-\l_1^0(0)$ are the
upper endpoints of the spectral bands. Thus, the endpoints of the
spectral bands $\l_n^{0\pm}$, $n\in\N_\n$, are the eigenvalues of
the matrices $\D(0)$ and $2\vk_+\1_\n-\D(0)$. \qq $\BBox$

%*****************************************************************
%\newpage
\section{\lb{Sec5}Stability estimates}
\setcounter{equation}{0}\lb{S4}

%****************************************************************

In the following theorem we obtain stability estimates for  the
Schr\"odinger operators satisfying the following conditions.

{\bf Condition U.} {\it   For the Schr\"odinger operator $H=\D+Q$
there exists a quasimomentum $\vt_+\in \T^d$ such that
$\l_n^{+}=\l_n(\vt_+)$ for all $n=1,\ldots,\n$. This point $\vt_+$
is called the $H$-upper quasimomentum.}

{\bf Condition L.}  {\it For the Schr\"odinger operator $H=\D+Q$
there exists a quasimomentum $\vt_-\in \T^d$ such that
$\l_n^{-}=\l_n(\vt_-)$ for all $n=1,\ldots,\n$. This point $\vt_-$
is called the $H$-lower quasimomentum. }

Below we show that there exists a large class of Schr\"odinger
operators, which satisfy  both Conditions U and L.
 We estimate a global variation of the spectrum and a
global variation of gap-length in terms of perturbations of the
entries of the Floquet matrices of the operators.

\begin{theorem} \lb{TpJ} {\bf(Stability Estimates)}
Let fundamental graphs $\G_f$, $\wt \G_f$ of periodic graphs $\G$,
$\wt \G$, respectively, have the same number of vertices $\n$. Let
the Schr\"odinger operators $H=\D+Q$ and $\wt H=\wt\D+\wt Q$ on the
graphs $\G$ and $\wt\G$, respectively, satisfy both Conditions U and
L. Then the spectral bands $\s_n=[\l_{n}^-,\l_{n}^+]$, $\wt
\s_n=[\wt\l_{n}^-,\wt\l_{n}^+]$, $n\in\N_\n$, and the gaps
$\g_n=[\l_{n}^+,\l_{n+1}^-]$, $\wt\g_n=[\wt\l_{n}^+,\wt\l_{n+1}^-]$,
$n\in\N_{\n-1}$, of the Schr\"odinger operators $H$ and $\wt H$,
respectively, satisfy:
\[
\lb{egap} |\l_1^--\wt\l_1^-|+|\l_\n^+-\wt\l_\n^+|+\sum_{n=1}^{\n-1}
\big||\g_n|-|\wt\g_n|\big|\le 2C,
\]
\[
\lb{esp} \sum_{n=1}^\n \big||\s_n|-|\wt \s_n|\big|\le 2C,
\]
where
\[
\lb{C0}
C=C(\vt_\pm,\wt\vt_\pm)=\|H(\vt_-)-\wt H(\wt\vt_-)\|_1+
\|H(\vt_+)-\wt H(\wt\vt_+)\|_1,
\]
$\|V\|_1=\sum\limits_{j,k=1}^\n|V_{jk}|<\iy$ and
$\vt_\pm,\wt\vt_\pm\in\T^d$ are the $H$- and $\wt H$-upper-lower
quasimomentum, respectively. In particular, we have
\[
\lb{egb} C= \ca 2\sum\limits_{n=1}^\n|q_n-\wt q_n|, & if \  \D(\vt_\pm)=\wt\D(\wt\vt_\pm) \\[10pt]
                               \sum\limits_{\pm}\|\D(\vt_\pm)-\wt \D(\wt\vt_\pm)\|_1, \qq & if \ Q=\wt Q \ac \ .
\]
\end{theorem}

\no {\bf Proof.} We have
\begin{multline}\lb{uu}
\big||\g_n|-|\wt\g_n|\big|=
\big|(\l_{n+1}^--\l_{n}^+)-(\wt\l_{n+1}^--\wt\l_{n}^+)\big|\\=
\big|(\l_{n+1}^--\wt\l_{n+1}^-)-(\l_{n}^+-\wt\l_{n}^+)\big|\le
\big|\l_{n+1}^--\wt\l_{n+1}^-\big|+\big|\l_{n}^+-\wt\l_{n}^+\big|.
\end{multline}
Since $H$, $\wt H$ satisfy Conditions U-L, then for all $n\in \N_\n$ we have
\[\label{es11}
\s_n=[\l_n^-,\l_n^+]=[\l_n(\vt_-),\l_n(\vt_+)], \qq
\wt\s_n=[\wt\l_n^-,\wt\l_n^+]=[\wt\l_n(\wt\vt_-),\wt\l_n(\wt\vt_+)].
%\qq \forall \;,
\]
%where $\vt_\pm, \wt\vt_\pm\in \T^d$ are the $H$- and $\wt
% H$-upper-lower quasimomentum, respectively.

Thus, applying Proposition \ref{TK}.ii and using \er{uu}, \er{es11}, we obtain
\begin{multline*}
\sum_{n=1}^{\n-1}
\big||\g_n|-|\wt\g_n|\big|\leq\sum_{n=1}^{\n-1}
\big|\l_{n+1}^--\wt\l_{n+1}^-\big|+\sum_{n=1}^{\n-1}
\big|\l_{n}^+-\wt\l_{n}^+\big|\\=\sum_{n=1}^{\n-1}
\big|\l_{n+1}(\vt_-)-\wt\l_{n+1}(\wt\vt_-)\big|+\sum_{n=1}^{\n-1}
\big|\l_{n}(\vt_+)-\wt\l_{n}(\wt\vt_+)\big|\leq
2C-\big|\l_1^--\wt\l_1^-\big|-
\big|\l_\n^+-\wt\l_\n^+\big|,
\end{multline*}
where $C=C(\vt_\pm,\wt\vt_\pm)$ is given by \er{C0},
which yields \er{egap}. Similar arguments and the  estimate
\[
\lb{anes}
\big||\s_n|-|\wt\s_n|\big|\le
|\l_n^+-\wt\l_n^+|+|\l_n^--\wt\l_n^-|
\]
yield \er{esp}.

If $\D(\vt_\pm)=\wt\D(\wt\vt_\pm)$, then we obtain
$
H(\vt_\pm)-\wt H(\wt\vt_\pm)=q-\wt q\,,
$
where $q=\diag(q_1,\ldots,q_\n)$, $\wt q=\diag(\wt q_1,\ldots,\wt q_\n)$.
Thus, in this case $C=2\sum\limits_{n=1}^\n|q_n-\wt q_n|$.
The proof for the case $Q=\wt Q$ is similar.  \BBox

\

\no\textbf{Remark.} 1) Theorem  \ref{TpJ}  can be applied for the
Schr\"odinger operators on both precise loop graphs and bipartite
loop graphs.

2)  There  is an old open problem to obtain estimates similar to
\er{egap}, \er{esp} in the case of the Schr\"odinger operators on
$\R^d$.

\

\no {\bf Proof of Proposition \ref{TpLG}.} i) Consider
the triangular lattice $\bT$ (Fig.\ref{f.tl}\emph{a}). The periods of $\bT$ are the vectors $a_1,a_2$. The fundamental graph $\bT_f$ consists of one vertex $v$, three bridge-loops
$\be_1,\be_2,\be_3$ and their inverse edges. The indices of the fundamental graph edges are given by
$$
\t(\be_1)=(1,0),\qqq \t(\be_2)=(0,1),\qqq \t(\be_3)=(1,1),
$$
see Fig.\ref{f.tl}\emph{b}. Since all bridges of $\bT_f$ are loops, the graph $\bT$ is a loop graph.
Since there is no point $\vt_0\in\T^2$ such that $\cos\lan\t
({\bf e_s}),\,\vt_0\ran=-1$ for all $s=1,2,3$, the loop graph $\bT$ is not precise.

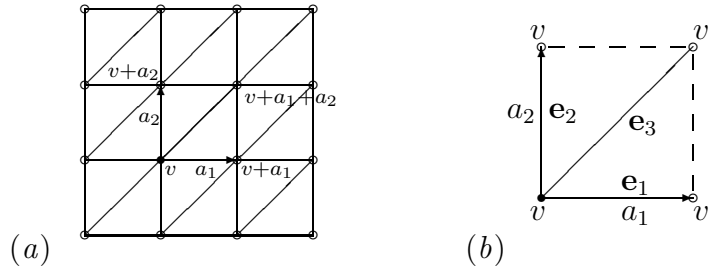
\begin{figure}[h]
\centering
\unitlength 1.0mm % = 2.845pt

\ifx\plotpoint\undefined\newsavebox{\plotpoint}\fi % GNUPLOT compatibility
\begin{picture}(80,45)(0,0)
% “реугольна€ решетка
\put(10,10){\line(1,0){30.00}}
\put(10,20){\line(1,0){30.00}}
\put(10,30){\line(1,0){30.00}}
\put(10,40){\line(1,0){30.00}}
\put(10,10){\line(0,1){30.00}}
\put(20,10){\line(0,1){30.00}}
\put(30,10){\line(0,1){30.00}}
\put(40,10){\line(0,1){30.00}}

\put(10,10){\line(1,1){30.00}}
\put(10,20){\line(1,1){20.00}}
\put(20,10){\line(1,1){20.00}}
\put(10,30){\line(1,1){10.00}}
\put(30,10){\line(1,1){10.00}}

\put(10,10){\circle{1}}
\put(20,10){\circle{1}}
\put(30,10){\circle{1}}
\put(40,10){\circle{1}}

\put(10,20){\circle{1.0}}
\put(20,20){\circle*{1.0}}
\put(30,20){\circle{1.0}}
\put(40,20){\circle{1}}

\put(10,30){\circle{1}}
\put(20,30){\circle{1.0}}
\put(30,30){\circle{1.0}}
\put(40,30){\circle{1}}

\put(10,40){\circle{1}}
\put(20,40){\circle{1}}
\put(30,40){\circle{1}}
\put(40,40){\circle{1}}
\put(20,20){\vector(1,0){10.00}}
\put(20,20){\vector(0,1){10.00}}
\put(20,20){\line(1,1){10.00}}

\put(20.5,18.0){$\scriptstyle v$}
\put(30.5,18.0){$\scriptstyle v+a_1$}
\put(31.0,27.5){$\scriptstyle v+a_1+a_2$}
\put(24.5,18.0){$\scriptstyle a_1$}

\put(17.0,25.0){$\scriptstyle a_2$}

\put(13.0,31.0){$\scriptstyle v+a_2$}

%**********************************

\put(70,15){\vector(1,0){20.00}}
\put(70,15){\vector(0,1){20.00}}
\multiput(71,35)(4,0){5}{\line(1,0){2}}
\multiput(90,16)(0,4){5}{\line(0,1){2}}

\put(80.5,16){$\be_1$}
\put(71.0,25.5){$\be_2$}
\put(81.5,24.0){$\be_3$}
\put(80.5,12.0){$a_1$}
\put(65.5,25.5){$a_2$}
\put(70,15){\line(1,1){20.00}}

\put(68.5,12.0){$v$}
\put(90.0,12.0){$v$}
\put(68.5,36.0){$v$}
\put(90.0,36.0){$v$}
\put(70,15){\circle*{1}}
\put(70,35){\circle{1}}
\put(90,15){\circle{1}}
\put(90,35){\circle{1}}

\put(0,7){(\emph{a})}
%**********************************
\put(60,7){(\emph{b})}
\end{picture}

\vspace{-0.5cm} \caption{\footnotesize \emph{a}) Triangular lattice $\bT$;  \
\emph{b}) the fundamental graph $\bT_f$.}
\label{f.tl}
\end{figure}

%%%%%%%%%%%%%%%%%%%%%%%%%%%%%%%

ii) Since the diameter of the
finite graph $\G'$ is small enough, in some  coordinate system all
bridges of the obtained graph $\G_\ast$ are also the bridges of the
given loop graph $\G$. Then all bridges of $\G_\ast$ are loops. The
last statement follows from the fact that the indices of the bridges
of $\G_\ast$ are the same as the bridge indices of $\G$. Thus, the
precise quasimomentum of $\G$ is also a precise quasimomentum of $\G_\ast$. \qq
\BBox

\

\no \textbf{Remark.} It is known that the spectrum of the Laplacian $\D$ on the triangular lattice $\bT$ is given by
$\s(\D)=\s_{ac}(\D)=[0,9]$.

\

{\bf Example.}  We take the lattice graph $\dL^d$, defined by
\er{dLg}. The lattice $\dL^d$ is a precise loop graph.  We take any connected finite graph $\G'\ss \R^d$ such
that its diameter is small enough. We joint  the graph $\G'$ with
each point from the vertex set $\Z^d$, identifying  any fixed vertex
of $\G'$ with each vertex of $\Z^d$. Due to Proposition \ref{TpLG}
this gives a new precise loop graph $\G_*$, which is
$\Z^d$-periodic.

\

%*****************************************
\no {\bf Proof of Theorem \ref{T100}.} i) Applying \er{eq.1}, \er{tl2.15}
we obtain $H(\vt)=H_0+\wt\D(\vt)$, where the matrix $\wt\D(\vt)$,
$\vt \in \T^d$, is given by
\[
\lb{ele}
\begin{aligned}
\wt\D=\diag\big(\wt\D_{11},\ldots,\wt\D_{\n\n}\big),\qqq
\wt\D_{nn}(\vt )=-\hspace{-5mm}\sum\limits_{{\bf e}=(v_n,\,v_n)
\in{\cB}_f}\cos\lan\t ({\bf e}),\,\vt\ran.
\end{aligned}
\]
From this identity it follows that
$
\wt\D_{nn}(0)\leq\wt\D_{nn}(\vt),\ \forall\,
(\vt,n)\in\T^d\ts\N_\n.
$
Then $\wt\D(0)\leq\wt\D(\vt)$ and we have
\[
H(0)=H_0+\wt\D(0)\leq H_0+\wt\D(\vt)=H(\vt).
\]
Applying Proposition \ref{MP}.iii to the last inequality, we obtain
$
\l_n(0)\leq\l_n(\vt )$ for all $(\vt,n)\in\T^d\ts\N_\n$,
which yields $\l_n^-=\min\limits_{\vt\in\T^d}\l_n(\vt )=\l_n(0)$.

ii)
The identity \er{ele} yields
$
\wt\D_{nn}(\vt)\leq\wt\D_{nn}(\vt_0)$ for all $
(\vt,n)\in\T^d\ts\N_\n
$,
since
$
\cos\lan\t ({\bf e}),\,\vt_0\ran=-1$
for the precise quasimomentum $\vt_0\in\T^d$ and
 for all bridges ${\bf
e}\in\cB_f$.
Then $\wt\D(\vt)\le\wt\D(\vt_0)$ and we have
\[
H(\vt )=H_0+\wt\D(\vt)\leq H_0+\wt\D(\vt_0)=H(\vt_0).
\]
Thus, Proposition \ref{MP}.iii gives
$
\l_n(\vt)\le\l_n(\vt_0)$ for all $(\vt,n)\in\T^d\ts\N_\n $, which yields
 $\l_n^+=\max\limits_{\vt
\in\T^d}\l_n(\vt)=\l_n(\vt_0)$ and due to \er{eq.5}, the band $\s_n(H)$ has the form \er{es}.

Using the formulas \er{es}, \er{ele}, we
obtain
\begin{multline*}
\sum\limits_{n=1}^\n|\s_n|=\sum\limits_{n=1}^\n\big(\l_n(\vt_0)-\l_n(0)\big)
=\Tr\big(H(\vt_0)-H(0)\big)=
\Tr\big(\wt\D(\vt_0)-\wt\D(0)\big)
%\sum\limits_{n=1}^\n\big(\wt\D_{nn}(\vt_0)-\wt\D_{nn}(0)\big)
=
2\b
\end{multline*}
and \er{es.1} has been proved.

Let all loops of $\G_f$ have the form $(v_k,v_k)$ for some vertex
$v_k\in V_f$. We need only prove the first identity in \er{es.111}.
Without loss of generality we may assume that $k=\n$. Then the
Floquet matrix $H(\vt)$ has the form
$$
H(\vt )=\left(
\begin{array}{cc}
  A & y \\
  y^\ast & \D_{\n\n}(\vt)+q_\n
\end{array}\right),
$$
where the entry $y\in\C^{\nu-1}$ is a vector and $A$ is a
self-adjoint $(\nu-1)\ts(\nu-1)$ matrix, do not depending on $\vt$.
The eigenvalues $\m_1\leq\ldots\leq\m_{\n-1}$ of $A$ are constant.
Then Proposition \ref{MP}.v gives that
$$
\l_1(\vt)\leq\mu_1\le\l_2(\vt)\leq\m_2\leq\ldots\leq\mu_{\n-1}\le\l_{\n}(\vt)\qqq \textrm{ for all }\vt\in\T^d.
$$
Then  the spectral bands of $H$ may only touch, but do not overlap,
i.e., $ |\s(H)|=\sum\limits_{n=1}^\n|\s_n|=2\b$. Thus, in this case the estimate \er{eq.7} and hence the first inequality in \er{GEga} become identities.
\qq  \BBox

\

 Schr\"odinger operators on  precise loop graphs satisfy Conditions
L-U. Then Theorem \ref{TpJ} can be applied  for the Schr\"odinger
operators on precise loop graphs. Note that if  $\lan \vt_0,\t({\bf
e})\ran/\pi $ is odd for all bridges $\be\in\cB_f$ and some vector
$\vt_0\in\{0,\pi\}^d$, then $\vt_0$ is a precise quasimomentum of $\G$.

\begin{theorem} \lb{TBi}
Let $\G$, $\wt \G$ be loop graphs with loop fundamental graphs
$\G_f$, $\wt \G_f$ having the same number of vertices $\n$. Let $H$
and $\wt H$ be the Schr\"odinger operators on the graphs $\G$ and
$\wt\G$, respectively. Then the spectral bands
$\s_n=[\l_{n}^{-},\l_{n}^{+}]$,
$\wt\s_n=[\wt\l_{n}^{-},\wt\l_{n}^{+}]$, $n\in\N_\n$, and the gaps
$\g_n=[\l_{n}^{+},\l_{n+1}^{-}]$, $\wt\g_n=
[\wt\l_{n}^{+},\wt\l_{n+1}^{-}]$, $n\in\N_{\n-1}$, of the
Schr\"odinger operators $H$ and $\wt H$, respectively, satisfy the
following estimates:

i) if $\G$ and $\wt \G$ are bipartite and regular of the same degree $\vk_+$ and $H=\D$,  $\wt
H=\wt \D$, then
\[
\lb{egapb}
\sum_{n=1}^{\n-1}
\big||\g_n|-|\wt\g_n|\big|\le 4C,
\]
\[
\lb{espb} \sum_{n=1}^\n \big||\s_n|-|\wt \s_n|\big|\le 4C,\qqq
C=\|\D(0)-\wt \D(0)\|_1;
\]
ii) if $\G$ is precise with a precise quasimomentum $\vt_0$, $\wt \G$ is
bipartite and regular of degree $\vk_+$ and $\wt H=\wt\D$, then
\[
\lb{egapbr} |\l_1^{-}|+|2\vk^+-\l_\n^{+}|+\sum_{n=1}^{\n-1}
\big||\g_n|-|\wt\g_n|\big|\le 2C,
\]
\[
\lb{espbr1}
\sum_{n=1}^\n \big||\s_n|-|\wt \s_n|\big|\le 2C,
\]
\[
\lb{espbr2} \qqq
C=C(\vt_0)=\|H(0)-\wt\D(0)\|_1+\|H(\vt_0)+\wt\D(0)-2\vk_+\1_\n\|_1;
\]
where $\|V\|_1=\sum\limits_{j,k=1}^\n|V_{jk}|<\iy$.
\end{theorem}

\no {\bf Proof.} i) Theorem \ref{TBG}.v gives
\[
\lb{es111}
\s_n(H)=[\l_n(0),2\vk_+-\l_{\n-n+1}(0)], \qq \s_n(\wt H)=[\wt\l_n(0),2\vk_+-\wt\l_{\n-n+1}(0)],\qq
\forall \;n\in \N_\n.
\]
Thus, using \er{uu},
\er{es111} and applying Proposition \ref{TK}.ii,  we obtain
\begin{multline*}
\sum_{n=1}^{\n-1}
\big||\g_n|-|\wt\g_n|\big|\leq\sum_{n=1}^{\n-1}
\big|\l_{n+1}^{-}-\wt\l_{n+1}^{-}\big|+\sum_{n=1}^{\n-1}
\big|\l_{n}^{+}-\wt\l_{n}^{+}\big|=\\\sum_{n=1}^{\n-1}
\big|\l_{n+1}(0)-\wt\l_{n+1}(0)\big|+\sum_{n=1}^{\n-1}
\big|\l_{\n-n+1}(0)-\wt\l_{\n-n+1}(0)\big|\\\leq
4C-|\l_1(0)-\wt\l_1(0)|-|\l_1(0)-\wt\l_1(0)|=4C,
\end{multline*}
where $C$ is defined in \er{espb}. Here we have used that, due to Proposition \ref{LP}, $\l_1(0)$ and $\wt \l_1(0)$ are the lower points of the spectra of the Laplacians $\D$ and $\wt\D$, respectively, i.e., $\l_1(0)=\wt \l_1(0)=0$.
Thus, \er{egapb} has been proved. Similar arguments and the estimate \er{anes}
give \er{espb}.

ii) We have
\[
\lb{est1}
\s_n(H)=[\l_n(0),\l_n(\vt_0)],\qqq \s_n(\wt H)=[\wt\l_n(0),2\vk_+-\wt\l_{\n-n+1}(0)], \qqq \forall \;n\in \N_\n.
\]
Using \er{uu},\er{est1}, the
constant $C$ from \er{espbr2} and applying Proposition \ref{TK}.ii, we obtain
\begin{multline*}
\sum_{n=1}^{\n-1} \big||\g_n|-|\wt\g_n|\big|\leq\sum_{n=1}^{\n-1}
\big|\l_{n+1}^--\wt\l_{n+1}^-\big|+\sum_{n=1}^{\n-1}
\big|\l_{n}^+-\wt\l_{n}^+\big|=\\\sum_{n=1}^{\n-1}
\big|\l_{n+1}(0)-\wt\l_{n+1}(0)\big|+\sum_{n=1}^{\n-1}
\big|2\vk_+-\wt\l_{\n-n+1}(0)-\l_{n}(\vt_0)\big|\\
\leq
2C-|\l_1(0)-\wt\l_1(0)|-|2\vk^+-\wt\l_1(0)-\l_\n(\vt_0)|=
2C-|\l_1^-|-|2\vk^+-\l_\n^+|.
\end{multline*}
This yields \er{egapbr}. Similar arguments and the estimate
\er{anes} yield \er{espbr1}. \qq \BBox

\

%****************************************************************
\section{\lb{Sec6} Various number of flat bands of  Laplacians on specific
graphs}
\setcounter{equation}{0}

We consider the question about the possible number of
non-degenerated  spectral bands and flat bands of
Laplacians on periodic graphs. We need to introduce another numbering $\L_n(\vt)$ of the eigenvalues $\l_n(\vt)$, $n\in\N_\n$, of the Floquet matrix $H(\vt)$, $\vt\in\T^d$. Here we separate all flat bands \er{fb1}, \er{fb2} from other bands of the Schr\"odinger  operator.

We have that  $\l_*$ is an eigenvalue
of $H$ iff $\l_*$ is an eigenvalue of $H(\vt)$ for any $\vt\in\T^d$
(see Proposition 4.2 in \cite{HN09}).
Thus, we can define the multiplicity of a flat band by:
an eigenvalue $\l_*$  of $H$ has the multiplicity $m$ iff
$\l_*=\const$ is an eigenvalue of $H(\vt)$ for each $\vt\in\T^d$ with the multiplicity $m$ (except maybe for a finite number of $\vt\in\T^d$).
Thus, if the operator $H$ has $r\ge 0$ flat bands, then we denote them by
\[\lb{fb1}
\m_j=\L_{\n-r+j}(\vt)=\const, \qqq \forall \ (j,\vt)\in \N_r\ts
\T^d,
\]
and they are labeled by
\[\lb{fb2}
\m_1\le \m_2\le\ldots\le  \m_r,
\]
counting multiplicities. Thus, all other eigenvalues
$\L_n(\vt)$, $n\in\N_{\n-r}$, are not constant. They can be enumerated
in increasing order (counting multiplicities) by
\[
\label{eq.3.1} \L_1(\vt )\leq\L_2(\vt )\leq\ldots\leq\L_{\nu-r}(\vt),
\qqq \forall\,\vt\in\T^d =\R^d/(2\pi\Z)^d.
\]
Define the \emph{spectral bands} $\gS_n=\gS_n(H), n\in\N_{\n-r}$, by
\[
\lb{ban} \gS_n=[\L_n^-,\L_n^+]=\L_n(\T^d).
\]
Each spectral band $\gS_n, n=1,\ldots, \n-r$, is open (non-degenerate), i.e., $\L_n^-<\L_n^+$. Thus, the number of open spectral bands of the
operator $H$ is $\n-r$. Some of them may overlap. Then the number of
gaps is at most $\n-r-1$.

Similarly, the eigenvalues of the Floquet matrix  $\D(\vt)$ for the Laplacian $\D$ will be denoted by
$\L^0_n(\vt)$, $n\in\N_\n$. The \emph{spectral bands} $\gS_n^0=\gS_n(\D)$ for the Laplacian have the form
\[
\lb{ban0.1} \gS_n^0=[\L_n^{0-},\L_n^{0+}]=\L_n^0(\T^d),\qq n\in\N_{\n-r}.
\]

We present  sufficient conditions for the existence of a flat band
for the Schr\"odinger operators.

\begin{lemma}\label{t8.1'}
Let $H=\D+Q$ on some graph $\G$, where its Floquet matrix $H(\vt)$
has the form
\[
\lb{ww.10}
H(\vt )=\left(
\begin{array}{cc}
  A & y \\
  y^\ast & a
\end{array}\right)(\vt),\qqq \vt\in\T^d,
\]
where $y(\vt)\in\C^{\nu-1}$ is a vector-valued function,
$a(\vt)$ is a function,  $A(\vt)$ is a self-adjoint
$(\nu-1)\ts(\nu-1)$ matrix-valued function. If $A(\vt)$ has an
eigenvalue $\eta(\cdot)=\eta=\const$ with multiplicity $m\ge 2$, then
$\eta$ is a flat band with multiplicity  $m-1$ of the operator
$H$.

\end{lemma}

\no {\bf Proof.} Due to Proposition \ref{MP}.v, there exist
eigenvalues $\m_1(\vt),\ldots,\m_{m-1}(\vt)$ of $H(\vt)$
satisfying
$$
\eta\leq\m_1(\vt)\leq\eta\le\ldots\eta\le\m_{m-1}(\vt)\leq\eta\qqq
 \textrm{ for all }\vt\in\T^d,
$$
which yields
$\m_1(\cdot)=\ldots=\m_{m-1}(\cdot)=\eta=\mathrm{const}$, i.e., $\eta$
is a flat band with multiplicity $m-1$ of the Schr\"odinger operator
$H$ on $\G$. \quad
$\BBox$
%*****************************************

%*************************************************************
\setlength{\unitlength}{1.0mm}
\begin{figure}[h]
\centering
\unitlength 1mm % = 2.845pt
\linethickness{0.4pt}
\ifx\plotpoint\undefined\newsavebox{\plotpoint}\fi % GNUPLOT compatibility
\begin{picture}(120,50)(0,0)

\put(10,10){\line(1,0){40.00}}
\put(10,30){\line(1,0){40.00}}
\put(10,50){\line(1,0){40.00}}
\put(10,10){\line(0,1){40.00}}
\put(30,10){\line(0,1){40.00}}
\put(50,10){\line(0,1){40.00}}

\put(10,10){\circle{1}}
\put(30,10){\circle{1}}
\put(50,10){\circle{1}}

\put(10,30){\circle{1}}
\put(30,30){\circle{1}}
\put(50,30){\circle{1}}

\put(10,50){\circle{1}}
\put(30,50){\circle{1}}
\put(50,50){\circle{1}}

\put(10,10){\vector(1,0){20.00}}
\put(10,10){\vector(0,1){20.00}}

\put(7.5,7.5){$\scriptstyle v_\n$}
\put(1,29){$\scriptstyle v_\n+a_2$}
\put(29,7){$\scriptstyle v_\n+a_1$}
\put(20,8){$\scriptstyle a_1$}
\put(6,20){$\scriptstyle a_2$}
\put(14,26.5){$\scriptstyle v_1$}
\put(23.8,17){$\scriptstyle v_{\nu-2}$}

\put(24.0,12){$\scriptstyle v_{\nu-1}$}
\put(20,26.5){$\scriptstyle v_2$}

\put(10,10){\line(2,3){10.00}}
\put(10,10){\line(3,2){15.00}}
\put(10,10){\line(1,3){5.00}}
\put(10,10){\line(3,1){15.00}}
\put(20.5,21){$\ddots$}
\put(15,25){\circle{1}}
\put(20,25){\circle{1}}
\put(25,15){\circle{1}}
\put(25,20){\circle{1}}
%**************************

\put(30,10){\line(2,3){10.00}}
\put(30,10){\line(3,2){15.00}}
\put(30,10){\line(1,3){5.00}}
\put(30,10){\line(3,1){15.00}}
\put(40.5,21){$\ddots$}
\put(35,25){\circle{1}}
\put(40,25){\circle{1}}
\put(45,15){\circle{1}}
\put(45,20){\circle{1}}
%******************************

\put(10,30){\line(2,3){10.00}}
\put(10,30){\line(3,2){15.00}}
\put(10,30){\line(1,3){5.00}}
\put(10,30){\line(3,1){15.00}}
\put(20.5,41){$\ddots$}
\put(15,45){\circle{1}}
\put(20,45){\circle{1}}
\put(25,35){\circle{1}}
\put(25,40){\circle{1}}
%**************************

\put(30,30){\line(2,3){10.00}}
\put(30,30){\line(3,2){15.00}}
\put(30,30){\line(1,3){5.00}}
\put(30,30){\line(3,1){15.00}}
\put(40.5,41){$\ddots$}
\put(35,45){\circle{1}}
\put(40,45){\circle{1}}
\put(45,35){\circle{1}}
\put(45,40){\circle{1}}
\put(-4,10){(\emph{a})}
%**********************************
\put(64,25){(\emph{b})}

\multiput(80,45)(4,0){5}{\line(1,0){2}}
\multiput(100,25)(0,4){5}{\line(0,1){2}}
\put(80,25){\line(1,0){20.00}}
\put(80,25){\line(0,1){20.00}}
\put(80,25){\circle{1}}
\put(100,25){\circle{1}}

\put(80,45){\circle{1}}
\put(100,45){\circle{1}}

\put(80,25){\vector(1,0){20.00}}
\put(80,25){\vector(0,1){20.00}}

\put(77.0,23.0){$\scriptstyle v_\n$}
\put(76.5,45){$\scriptstyle v_\n$}
\put(101,23){$\scriptstyle v_\n$}
\put(101,45){$\scriptstyle v_\n$}
\put(90,23){$\scriptstyle a_1$}
\put(76,35){$\scriptstyle a_2$}

\put(84,41.5){$\scriptstyle v_1$}
\put(93.8,32.5){$\scriptstyle v_{\nu-2}$}

\put(94.0,27){$\scriptstyle v_{\nu-1}$}
\put(90,41.5){$\scriptstyle v_2$}

\put(80,25){\line(2,3){10.00}}
\put(80,25){\line(3,2){15.00}}
\put(80,25){\line(1,3){5.00}}
\put(80,25){\line(3,1){15.00}}
\put(90.5,36){$\ddots$}
\put(85,40){\circle{1}}
\put(90,40){\circle{1}}
\put(95,30){\circle{1}}
\put(95,35){\circle{1}}
%*********************************
\put(64,10){(\emph{c})}

\put(80,10){\line(1,0){50.00}}
\put(80,9){\line(0,1){2.00}}
\put(95,9){\line(0,1){2.00}}
\put(84,9){\line(0,1){2.00}}
\put(130,9){\line(0,1){2.00}}
\put(85,10){\circle*{1}}

\put(80,9.8){\line(1,0){4.00}}
\put(80,10.2){\line(1,0){4.00}}

\put(95,9.8){\line(1,0){35.00}}
\put(95,10.2){\line(1,0){35.00}}

\put(110,12){$\gS_2^0$}
\put(80,11.5){$\scriptstyle\gS_1^0$}
\put(85,11.5){$\scriptstyle\m=1$}
\put(79,6){$\scriptstyle0$}
\put(95,6){$\scriptstyle3$}
\put(123.0,5.5){$\scriptstyle\frac{11+\sqrt{89}}2$}
\put(82.0,5.5){$\scriptstyle\frac{11-\sqrt{89}}2$}
\end{picture}
\vspace{-0.5cm} \caption{\footnotesize  \emph{a}) $\Z^2$-periodic graph $\G$;\quad
\emph{b}) the fundamental graph $\G_f$; only 2 unoriented loops in
the vertex $v_\n$ are bridges;\quad \emph{c}) the spectrum of the
Laplacian ($\n=3$).} \label{ff.10}
\end{figure}
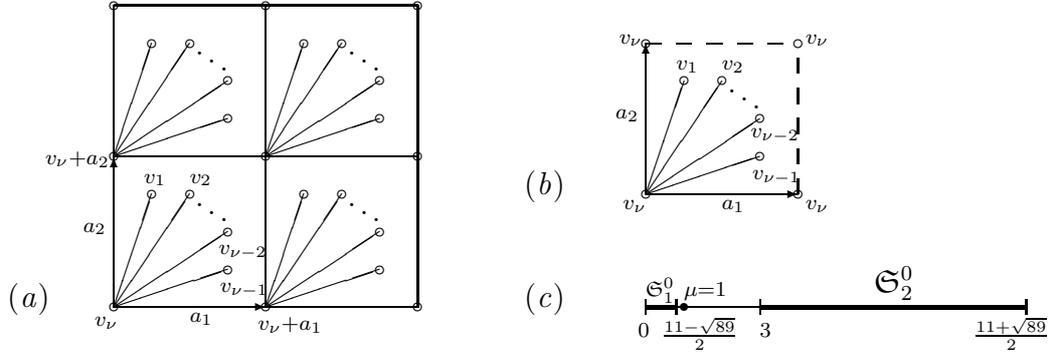

\

We describe the spectrum of the Laplace and Schr\"odinger operators
on the specific precise loop graph $\G$, shown in
Fig.\ref{ff.10}\emph{a}.

\begin{proposition}\lb{TG1}
Let $\G_f$ be obtained from the fundamental graph $\dL^d_f$ of the
$d$-dimension lattice $\dL^d$ by adding $\n-1\geq1$ vertices
$v_1,\ldots,v_{\n-1}$ and $\n-1$ unoriented edges
$(v_1,v_\n),\ldots,(v_{\n-1},v_\n)$  with zero indices (see
Fig.\ref{ff.10}b), $v_\n$ is the unique vertex of $\dL^d_f$. Then
$\G$ is a precise loop graph with the precise quasimomentum
$\vt_\pi=(\pi,\ldots,\pi)\in\T^d$ and satisfies

i) The spectrum of the Laplacian on $\G$  has the form
\[
\lb{acs1}
\s(\D)=\s_{ac}(\D)\cup\s_{fb}(\D),\qqq \s_{fb}(\D)=\{1\},
\]
where the  flat band $1$ has multiplicity $\n-2$ and  $\s_{ac}(\D)$
has only two bands $\gS_1^0$ and $\gS_2^0$ given by
\[
\lb{acs2}
\begin{array}{c}
\s_{ac}(\D)=\gS_1^0\cup\gS_2^0, \\ [6 pt]
 \gS_1^0=\big[\,0,\textstyle x-\sqrt{x^2-4d}\, \big],\qq
\gS_2^0=\big[\,\n,\textstyle x+\sqrt{x^2-4d}\, \big],\qqq
x={\n+4d\/2}\,.
\end{array}
\]

ii) The spectrum of Schr\"odinger operators $H=\D+Q$ on $\G$ has the form
\[
\lb{acs333}
\s(H)=\bigcup_{n=1}^\n\big[\l_n(0),\l_n(\vt_\pi)\big].
\]

iii) Let $q_\n=0$ and let all other values of  the potential $q_1,\ldots,q_{\n-1}$ at the
vertices of the fundamental graph $\G_f$ be distinct. Then
$\s(H)=\s_{ac}(H)$, i.e., $\s_{fb}(H)=\varnothing$.

iv) Let among the numbers $q_1,\ldots,q_{\n-1}$ there exist a value
$q_\ast$ of multiplicity $m$. Then the spectrum of the Schr\"odinger
operator $H$ on $\G$ has the flat band $q_\ast+1$ of multiplicity
$m-1$.

v) The Lebesgue measure of the spectrum of Schr\"odinger operators $H$ on $\G$
satisfies
\[
\lb{sp2}
\textstyle|\s(H)|=4d.
\]
\end{proposition}
\no {\bf Proof.} i) -- ii) The fundamental graph $\G_f$ consists of
$\n\geq2$  vertices $v_1,v_2,\ldots,v_\n$; ${\n-1}$ unoriented edges
${(v_1,v_\n),\ldots,(v_{\n-1},v_\n)}$ with zero indices and $d$
unoriented loops in the vertex $v_\n$ with the indices
$(\pm1,\ldots,0),\ldots,(0,\ldots,\pm1)$. Since all bridges of
$\G_f$ are loops and $\cos\lan\t ({\bf e}),\,\vt_\pi\ran=-1$ for all
bridges $\be\in\cB_f$, the graph $\G$ is a precise loop graph with
the precise quasimomentum $\vt_\pi$. Then, by Theorem \ref{T100}.ii, the
spectral bands of $H$ are given by
$
\s_n(H)=[\l_n^-,\l_n^+]=[\l_n(0),\l_n(\vt_\pi)]$, for all $
n\in\N_\n,
$
and item ii) has been proved.

According to (\ref{l2.15}) we have
\[
\lb{fc0}
\D(\vt)=\left(
\begin{array}{cccc}
  1 & 0  & \ldots & -1  \\[2pt]
  0 & 1 & \ldots & -1 \\[2pt]
  \ldots & \ldots & \ldots & \ldots  \\[2pt]
   -1 & -1 & \ldots & \n-1+\xi  \\
\end{array}\right),\qqq
\ca c_0=c_1+\ldots+c_d\\
c_j=\cos\vt_j, \qq j\in\N_d\\
\x=2(d-c_0)\ac.
\]
Using the formula \er{det}, we obtain
$$
\textstyle\det\big(\D(0)-\l \1_\n\big)=(1-\l)^{\n-2}\,\l\,(\l-\n),
$$
$$
\textstyle\det\big(\D(\vt_\pi)-\l \1_\n\big)=(1-\l)^{\n-2}
\big(\l^2-2\l x+4d\big), \qq x={\n+4d\/2}.
$$
Then the eigenvalues of $\D(0)$ and $\D(\vt_\pi)$ have the form
$$
\textstyle\l_1^0(0)=0,\qqq \l_2^0(0)=\ldots=\l_{\n-1}^0(0)=1,\qqq
\l_\n^0(0)=\n,
$$
$$
\textstyle\l_1^0(\vt_\pi)=x-\sqrt{x^2-4d},\qq
\l_2^0(\vt_\pi)=\ldots=\l_{\n-1}^0(\vt_\pi)=1,\qq
\l_\n^0(\vt_\pi)=x+\sqrt{x^2-4d}.
$$
Thus, the spectrum of the Laplacian on $\G$  has the form \er{acs1}, \er{acs2}.

iii) According to \er{fc0}, we have
\[
\lb{pr32}
H(\vt)=\left(
\begin{array}{cccc}
  1+q_1 & 0  & \ldots & -1  \\[2pt]
  0 & 1+q_2 & \ldots & -1 \\[2pt]
  \ldots & \ldots & \ldots & \ldots  \\[2pt]
  -1 & -1 & \ldots & \n-1+\xi  \\
\end{array}\right).
\]
 Using the formula \er{det},
we write the characteristic polynomial of the matrix $H(\vt)$ in the
form
\[
\lb{lco1} \textstyle  \det(H(\vt)-\l
\1_\n)=(\n-1+\x-\l)W(\l)-W'(\l),
\]
where
\[
\lb{WWW}
W(\l)=(1+q_1-\l)\ldots(1+q_{\n-1}-\l)\,.
\]

We show that $\s(H)=\s_{ac}(H)$ by the contradiction. Let a point  $\l$
be a flat band of the Schr\"odinger operator $H$. Then we obtain $
\det(H(\vt)-\l\1_\n)=0$ for all $\vt\in\T^d$. Since the linear
combination \er{lco1} of the linearly independent functions
$c_1,\ldots,c_d,1$ is equal to 0, then
\[
\lb{sys}
W(\l)=0,\qqq  W'(\l)=0.
\]
All values of the potential $q_1,\ldots,q_{\n-1}$ are distinct. Then
each  zero of the function $W$ is simple. This contradicts the
identities \er{sys}. Thus, $\s(H)=\s_{ac}(H)$.

iv) Without loss of generality we assume that
$q_1=q_2=\ldots=q_m=q^\ast$. Then $\l=q^\ast+1$ is a zero of
multiplicity $m$ of the function $W$, defined  by \er{WWW}, and is a
zero of multiplicity $m-1$ of the function $W'$. Then the identity
\er{sys} that defined all flat bands has the solution $\l=q^\ast+1$
of multiplicity $m-1$. Thus, $\l=q^\ast+1$ is a flat band of $H$ with multiplicity $m-1$.

v) All bridges of $\G_f$ are the loops in the vertex $v_\n$. Thus,
using that the number of fundamental graph bridges $\b=2d$, the
identity (\ref{es.111}) has the form \er{sp2}.\qq \BBox

\

%*****************************************
\setlength{\unitlength}{0.8mm}
\begin{figure}[h]
\centering

\unitlength 1mm % = 2.845pt
\linethickness{0.4pt}
\ifx\plotpoint\undefined\newsavebox{\plotpoint}\fi % GNUPLOT compatibility
\begin{picture}(110,40)(0,0)
\put(-5,10){(\emph{a})}
\put(10,10){\line(1,0){30.00}}
\put(10,25){\line(1,0){30.00}}
\put(10,40){\line(1,0){30.00}}
\put(10,10){\line(0,1){30.00}}
\put(25,10){\line(0,1){30.00}}
\put(40,10){\line(0,1){30.00}}

\put(10,10){\circle*{1}}
\put(15,10){\circle*{1}}
\put(20,10){\circle*{1}}
\put(25,10){\circle{1}}
\put(30,10){\circle{1}}
\put(35,10){\circle{1}}
\put(40,10){\circle{1}}

\put(10,25){\circle{1}}
\put(15,25){\circle{1}}
\put(20,25){\circle{1}}
\put(25,25){\circle{1}}
\put(30,25){\circle{1}}
\put(35,25){\circle{1}}
\put(40,25){\circle{1}}

\put(10,40){\circle{1}}
\put(15,40){\circle{1}}
\put(20,40){\circle{1}}
\put(25,40){\circle{1}}
\put(30,40){\circle{1}}
\put(35,40){\circle{1}}
\put(40,40){\circle{1}}

\put(10,15){\circle*{1}}
\put(10,20){\circle*{1}}
\put(10,30){\circle{1}}
\put(10,35){\circle{1}}

\put(25,15){\circle{1}}
\put(25,20){\circle{1}}
\put(25,30){\circle{1}}
\put(25,35){\circle{1}}

\put(40,15){\circle{1}}
\put(40,20){\circle{1}}
\put(40,30){\circle{1}}
\put(40,35){\circle{1}}
\put(8.0,7.0){$v_5$}
\put(13.0,7.0){$v_1$}
\put(18.0,7.0){$v_2$}
\put(6.0,14.0){$v_3$}
\put(6.0,19.0){$v_4$}

%*********************************
\put(55,10){(\emph{b})}
\put(70,10){\line(1,0){55.00}}
\put(70,9.8){\line(1,0){10.00}}
\put(70,10.2){\line(1,0){10.00}}

\put(84,9.8){\line(1,0){16.00}}
\put(84,10.2){\line(1,0){16.00}}

\put(120,9.8){\line(1,0){5.00}}
\put(120,10.2){\line(1,0){5.00}}

\put(70,9){\line(0,1){2.00}}
\put(84,9){\line(0,1){2.00}}
\put(120,9){\line(0,1){2.00}}
\put(125,9){\line(0,1){2.00}}

\put(80,10){\circle*{1.5}}
\put(100,10){\circle*{1.5}}

\put(120.5,12){$\scriptstyle\gS_3^0$}
\put(88,12){$\scriptstyle\gS_2^0$}

\put(72,12){$\scriptstyle\gS_1^0$}
\put(77,12){$\scriptstyle\m_1=1$}
\put(97,12){$\scriptstyle\m_2=3$}
\put(69,6){$\scriptstyle0$}

\put(81.3,5){$\scriptstyle\frac{\sqrt{17}-1}2$}
\put(122.3,5){$\scriptstyle\frac{\sqrt{17}+7}2$}
\put(119.3,6){$\scriptstyle5$}
\end{picture}
\vspace{-0.5cm}
\caption{\footnotesize  \emph{a}) Graph $\G$ obtained by adding two vertices on each edge of
the lattice $\dL^2$;\quad \emph{b})~the spectrum of the Laplacian.} \label{f.10}
\end{figure}
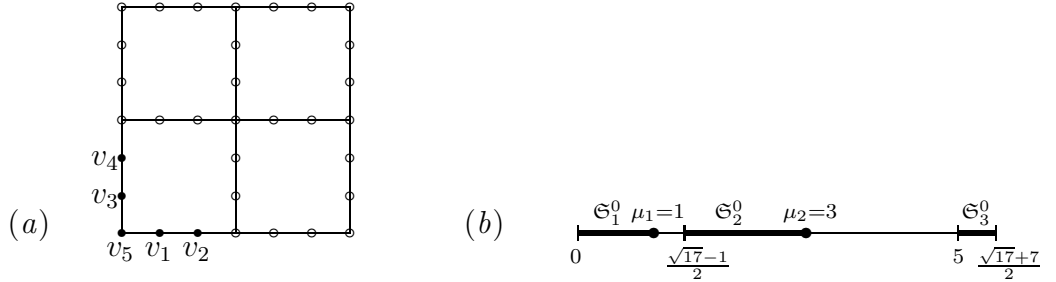

We describe the spectrum of the Laplace operator on the graph $\G$,
shown in Fig.\ref{f.10}\emph{a}.

\begin{proposition}
\lb{TG2} Let $\G$ be the graph obtained from the lattice graph
$\dL^d$ by adding $N$ vertices on each edge of $\dL^d$ (for $N=2$
see Fig.\ref{f.10}a). Then the fundamental graph $\G_f$ has
$\n=dN+1$ vertices and the spectrum of the Laplacian on $\G$ has the
form
\[
\lb{acs3}
\s(\D)=\s_{ac}(\D)\cup\s_{fb}(\D),
\]
\[
\lb{acs4} \textstyle \s_{fb}(\D)=\big\{2+2\cos{\pi n\/ N+1}\;:
\;n=1,\ldots,N\big\},
\]
where the absolutely continuous part $\s_{ac}(\D)$ consists of $N+1$
open spectral bands $\gS_1^0,\ldots,\gS_{N+1}^0$ and each flat band
has multiplicity $d-1$ and is an endpoint of a spectral band.
\end{proposition}
\no {\bf Proof.} Let $N=1$. The fundamental graph $\G_f$ has
$\n=d+1$ vertices. For each  $\vt\in \T^d$ the matrix $\D(\vt)$
defined by (\ref{l2.15}) has the form
\begin{equation}
\label{z2}
\D(\vt)=\left(
\begin{array}{cccc}
   2 & 0 & \ldots & -1-e^{i\vt_1}\\
   0 & 2 & \ldots & -1-e^{i\vt_2}\\
   \ldots & \ldots & \ldots & \ldots\\
   -1-e^{-i\vt_1} & -1-e^{-i\vt_2} & \ldots & 2d\\
\end{array}\right).
\end{equation}
Identity \er{det} yields
\[
\lb{sl1'} \textstyle\det\big(\D(\vt)-\l
\1_\n\big)=(2-\l)^d\big(2d-\l-{2(
d+c_0)\/2-\l}\,\big)=(2-\l)^{d-1}\big(\l^2-2(d+1)\l+2(d-c_0)\big),
\]
where $c_0$ is defined by \er{fc0}. Then the eigenvalues of the
matrix  $\D(\vt)$ have the form
\[
\lb{ev11} \textstyle \L_s^0(\vt)=d+1+(-1)^s \sqrt{d^2+1+2c_0}\,,\qq
s=1,2;
\]
\[
\lb{ev2} \L_3^0(\cdot)=\ldots=\L_{d+1}^0(\cdot)=2.
\]
The identities \er{ev2} give that the point 2 is a flat band of $\D$
of         multiplicity $d-1$. From the identity \er{ev11} we obtain
that
$$
\L_1^{0-}=\min_{\vt\in\T^d}\L_1^0(\vt)=\L_1^0(0)=0,\qqq \L_1^{0+}=
\max_{\vt\in\T^d}\L_1^0(\vt)=\L_1^0(\vt_\pi)=2,\qqq
\vt_\pi=(\pi,\ldots,\pi),
$$
$$
\L_2^{0-}=\min_{\vt\in\T^d}\L_2^0(\vt)=\L_2^0(\vt_\pi)=2d,\qqq
\L_2^{0+}=\max_{\vt\in\T^d}\L_2^0(\vt)=\L_2^0(0)=2d+2,
$$
and thus, $\s_{ac}(\D)=\gS_1^0\cup\gS_2^0$, where $\gS_1^0=[0,2]$, $\gS_2^0=[2d,2d+2]$.

Let $N\geq2$. The fundamental graph $\G_f$ has $\n=dN+1$ vertices.
It is more convenient in this case instead of the Laplacian $\D$
consider the  operator $J=2\1-\D$, where $\1$ is the identity
operator. Then $\s(\D)=2-\s(J)$. The identity holds as an identity
between sets.

The Floquet matrix $J(\vt)$ of the operator $J$ has the form
$$
J(\vt)=2\1_\nu-\D(\vt)=\left(
\begin{array}{cc}
  J_{dN} & y(\vt) \\
  y^\ast(\vt) & 2-2d
\end{array}\right),
$$
where the $dN\ts dN$ matrix $J_{dN}=\diag(J_N,\ldots,J_N)$ and
the vector $y(\vt)$ are given by
\[
\lb{sl0}
\begin{aligned}
J_N=\left(
\begin{array}{cccc}
  0 & 1 & 0 &\ldots \\
  1 & 0 & 1 &\ldots\\
  0 & 1 & 0 &\ldots \\
   \ldots & \ldots & \ldots &\ldots \\
\end{array}\right),\qq
y(\vt)=\ma y_1(\vt_1)\\
\vdots \\
 y_d(\vt_d) \am,
 \qqq
y_s(\vt_s)=\ma 1\\
0\\
\vdots \\
0\\
e^{i\vt_s} \am \in \C^N,
\end{aligned}
\]
$s\in\N_d.$ It is known that all eigenvalues of the $N\ts N$ Jacobi
matrix $J_N$ have the form
\begin{equation}
\label{ww.17}
\textstyle\mu_n=2\cos\frac{\pi n}{N+1}\in(-2,2),\qquad n=1,\ldots, N,
\end{equation}
and they are distinct. Then the matrix $J_{dN}$ has $N$ distinct
eigenvalues $\m_n$ of multiplicity $d$. Thus, due to Lemma
\ref{t8.1'}, the operator $J$ on $\G$ has at least $N=\frac{\n-1}d$
flat bands $\m_n, n\in \N_N$, each of which has multiplicity $d-1$.

We describe $\s_{ac}(J)$. Identity \er{det} yields
\[
\lb{sl1} \det\big(\l \1_\n-J(\vt)\big)=\big(\l-2+2d-y^\ast(\vt) (\l
\1_{dN}-J_{dN})^{-1}y(\vt)\big)\det\big(\l \1_{dN}-J_{dN}\big).
\]
From the form of the matrix $J_{dN}$ it follows that
\[
\lb{sl2}
\det\big(\l \1_{dN}-J_{dN}\big)=\cD^d_N(\l),\qqq \textrm{where }\ \cD_N(\l)=
\det\big(\l \1_N-J_N\big).
\]
A direct calculation gives
\[
\lb{sl3}
(\l \1_{dN}-J_{dN})^{-1}=\diag(B_N,\ldots,B_N), \qqq B_N=(\l \1_N-J_N)^{-1}\,,
\]
\[
\lb{sl5} y^\ast(\vt)(\l \1_{dN}-J_{dN})^{-1}y(\vt)=y_1^\ast(\vt_1)
B_Ny_1(\vt_1)+\ldots+y_d^\ast(\vt_d) B_Ny_d(\vt_d),
\]
\[
\lb{sl16}
\textstyle y_s^\ast(\vt_s)
B_Ny_s(\vt_s)=\frac2{\cD_N(\l)}\big(\cD_{N-1}(\l)+\cos\vt_s\,\big),\qqq s\in\N_d.
\]
Substituting \er{sl2}, \er{sl5} and \er{sl16}  into \er{sl1}, we obtain
\[
\lb{sl4}
\det\big(\l\1_\n-J(\vt)\big)=\cD^{d-1}_N(\l)\big((\l-2+2d)\cD_N(\l)-
2d\cD_{N-1}(\l)-2c_0
\big).
\]
The determinant $\cD_n(\l)$ satisfies the following recurrence relations
 (the Jacobi equation)
\[
\lb{rel}
\begin{array}{c}
\cD_{n+1}(\l)=\l\,\cD_n(\l)-\cD_{n-1}(\l),\qqq \forall \ n\in\N_{\n-1},  \\[6pt]
\textrm{ where } \ \cD_0(\l)=1,\qq \cD_1(\l)=\l.
\end{array}
\]
Thus, the determinants $\cD_n(\l)$, $n=1,2,\ldots,$ are the
Chebyshev polynomials of the second  kind and the following
identities hold true
\[
\lb{sl6}
\textstyle\cD_n(\l)=\frac{\sin(n+1)\vp}{\sin\vp}\,,\qq \textrm{where } \lambda=2\cos\vp\in(-2,2), \qqq
n\in\N_\n.
\]

The eigenvalues of the matrix $J(\vt)$ are zeros of the equations
\[
\lb{spe}
\cD^{d-1}_N(\l)=0,\qqq (\l-2+2d)\cD_N(\l)-
2d\cD_{N-1}(\l)=2c_0\,.
\]
The first identity gives all flat bands $\m_n$, $n=1,\ldots,N$, of $J$ defined by \er{ww.17}. Then the set of all flat bands of the Laplacian has the form \er{acs4}.
From the second identity we obtain that all flat bands $\m_n$, $n=1,\ldots,N$,
are endpoints of the spectral bands. Indeed, since
$\cD^{d-1}_N(\m_n)=0$, $n=1,\ldots,N$, the second identity in \er{spe} gives
\[
\lb{spe1}
\textstyle -\cD_{N-1}(\m_n)=\frac{c_0}d\,.
\]
But, using \er{sl6}, we have
$$
\textstyle \cD_{N-1}(\m_n)=\cD_{N-1}\big(2\cos\frac{\pi n}{N+1}\big)=
\frac{\sin \frac{N\pi n}{N+1}}{\sin\frac{\pi n}{N+1}}=(-1)^{n+1}.
$$
Thus, the identity \er{spe1} has the form $ (-1)^{n}=\frac{c_0}d $
and it holds true when $\vt=0$ (if $n$ is even) or
$\vt=(\pi,\ldots,\pi)$  (if $n$ is odd). Therefore, all flat bands
$\m_n$, $n=1,\ldots,N$, are contained in the absolutely continuous
spectrum. On the other hand, by Proposition \ref{MP}.v, all spectral
bands are separated by the flat bands $\m_n$, $n=1,\ldots,N$. Then
all flat bands $\m_n$ are endpoints of the spectral bands of the operator $J$, which yields that the flat bands \er{acs4} are endpoints of the spectral bands of   $\D$. \qq
\BBox

\

%***************************************************************
\no {\bf Proof of Theorem \ref{Tfb}.}
This statement is a direct consequence of Proposition
\ref{TG1}.
\qq \BBox

%***********************************
\section{\lb{Sec7} Crystal models}
\setcounter{equation}{0}

%***********************************
\subsection{Hexagonal lattice}
We consider the hexagonal lattice $\bG=(V,\cE)$, shown in Fig.\ref{ff.0.3}\emph{a}.
The periods of $\bG$ are the
vectors $a_1=(3/2,\sqrt{3}/2)$, $a_2=(0,\sqrt{3}\,)$ (the
coordinates of $a_1,a_2$ are taken in the orthonormal basis
$e_1,e_2$). The vertex set and the edge set are given by
$$
\textstyle V=\Z^2\cup\big(\Z^2+\big(\frac13\,,\frac13\big)\big),
$$
$$
\textstyle \cE=\big\{\big(\mm,\mm+\big(\frac13\,,\frac13\big)\big),
\big(\mm,\mm+\big(-\frac23\,,\frac13\big)\big),\\
\big(\mm,\mm+\big(\frac13\,,-\frac23\big)\big)\quad\forall\,\mm\in\Z^2\big\}.
$$
Coordinates of all vertices are taken in the basis
$a_1,a_2$.
The fundamental graph $\bG_f=(V_f, \cE_f)$, where $V_f=\{v_1, v_2\}$ consists of two vertices, \emph{multiple} edges $\be_1=\be_2=\be_3=(v_1,v_2)$
(Fig.\ref{ff.0.3}\emph{b}) and their inverse edges $\bar\be_1=\bar\be_2=\bar\be_3$ with
the indices
$\t(\be_1)=(0,0)$, $\t(\be_2)=(1,0)$, $\t(\be_3)=(0,1)$.

We consider the Schr\"odinger operators $H=\D+Q$ on the hexagonal lattice $\bG$.
Without loss of generality we may assume that
\[
\lb{q12}
Q(v_1)=q_1,\qqq Q(v_2)=q_2=-q_1.
\]

\setlength{\unitlength}{1.0mm}
\begin{figure}[h]
\centering

\unitlength 1mm % = 2.845pt
\linethickness{0.4pt}
\ifx\plotpoint\undefined\newsavebox{\plotpoint}\fi % GNUPLOT compatibility
\begin{picture}(100,45)(0,0)

\put(5,10){(\emph{a})}

% √ексагональна€ решетка
\put(14,10){\circle{1}}
\put(28,10){\circle{1}}
\put(34,10){\circle{1}}
\put(48,10){\circle{1}}

\put(18,16){\circle{1}}
\put(24,16){\circle*{1}}
\put(38,16){\circle{1}}
\put(44,16){\circle{1}}

\put(14,22){\circle{1}}
\put(28,22){\circle*{1}}
\put(34,22){\circle{1}}
\put(48,22){\circle{1}}

\put(18,28){\circle{1}}
\put(24,28){\circle{1}}
\put(38,28){\circle{1}}
\put(44,28){\circle{1}}

\put(14,34){\circle{1}}
\put(28,34){\circle{1}}
\put(34,34){\circle{1}}
\put(48,34){\circle{1}}

\put(18,40){\circle{1}}
\put(24,40){\circle{1}}
\put(38,40){\circle{1}}
\put(44,40){\circle{1}}

\put(28,10){\line(1,0){6.00}}
\put(18,16){\line(1,0){6.00}}
\put(38,16){\line(1,0){6.00}}

\put(28,22){\line(1,0){6.00}}
\put(18,28){\line(1,0){6.00}}
\put(38,28){\line(1,0){6.00}}

\put(28,34){\line(1,0){6.00}}
\put(18,40){\line(1,0){6.00}}
\put(38,40){\line(1,0){6.00}}

\put(14,10){\line(2,3){4.00}}
\put(34,10){\line(2,3){4.00}}
\put(24,16){\line(2,3){4.00}}
\put(44,16){\line(2,3){4.00}}

\put(14,22){\line(2,3){4.00}}
\put(34,22){\line(2,3){4.00}}
\put(24,28){\line(2,3){4.00}}
\put(44,28){\line(2,3){4.00}}

\put(14,34){\line(2,3){4.00}}
\put(34,34){\line(2,3){4.00}}

\put(28,10){\line(-2,3){4.00}}
\put(48,10){\line(-2,3){4.00}}
\put(38,16){\line(-2,3){4.00}}
\put(18,16){\line(-2,3){4.00}}

\put(28,22){\line(-2,3){4.00}}
\put(48,22){\line(-2,3){4.00}}
\put(38,28){\line(-2,3){4.00}}
\put(18,28){\line(-2,3){4.00}}

\put(28,34){\line(-2,3){4.00}}
\put(48,34){\line(-2,3){4.00}}

\put(30,18){$\scriptstyle a_1$}
\put(20.5,22){$\scriptstyle a_2$}

\put(24,16){\vector(0,1){12.0}}
\put(33,21.3){\vector(3,2){0.5}}

%\drawline(24,16)(34,22)
\qbezier(24,16)(29,19)(34,22)

\put(24.8,21.5){$\scriptstyle v_1$}
\put(35,21.5){$\scriptstyle v_2+a_1$}
\put(25.0,28.0){$\scriptstyle v_2+a_2$}
\put(17.8,13.5){$\scriptstyle O=v_2$}

%***************************
\put(75,10){\circle*{1}}
\put(83,21){\circle*{1}}

\put(75,30){\circle{1}}
\put(95,40){\circle{1}}
\put(95,20){\circle{1}}

\put(75,10){\vector(0,1){20.0}}
\put(75,10){\vector(2,1){20.0}}

\multiput(95,20)(0,7){3}{\line(0,1){4}}
\put(75,30){\line(2,1){4.0}}
\put(82,33.5){\line(2,1){4.0}}
\put(89,37){\line(2,1){4.0}}

\qbezier(83,21)(89,20.5)(95,20)
\qbezier(83,21)(79,15.5)(75,10)
\qbezier(83,21)(79,25.5)(75,30)

\put(71,8.0){$v_2$}
\put(83,22){$v_1$}
\put(96,19.0){$v_2$}
\put(71.0,31.0){$v_2$}
\put(93.5,42.0){$v_2$}
\put(85,13){$a_1$}
\put(70,20.0){$a_2$}

\put(76,17.2){$\mathbf{e}_1$}
\put(89,21.2){$\mathbf{e}_2$}
\put(78,26.6){$\mathbf{e}_3$}

\put(64,10){(\emph{b})}
\end{picture}

\vspace{-0.5cm} \caption{\footnotesize  \emph{a}) Graphene $\bG$; \quad \emph{b}) the fundamental graph $\bG_f$ of the graphene.} \label{ff.0.3}
\end{figure}
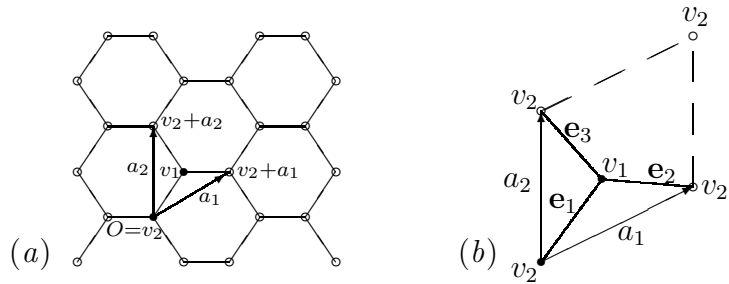

In order to present our result we define
the so-called  symbol $D_m(t)$, $t\in \R^2$, of the 2D Dirac operator  by
\[
\label{DDi}
\begin{aligned}
D_m(t)=\s_1 t_1+\s_2   t_2+m\s_3,\qqq t =(t_1,t_2)\in \R^2,
\end{aligned}
\]
where $m$ is the mass and  $\s_1,\s_2,\s_3$ are the Pauli matrices given by
\[
\label{Dm}
\begin{aligned}
 \s_1=\ma 0 & 1 \\
         1 & 0 \am, \qqq \s_2=\ma 0 & -i \\
         i & 0 \am, \qqq \s_3=\ma 1 & 0 \\
         0 & -1 \am.
\end{aligned}
\]
The 2D Dirac operator has the form $D_m(\pa)$, $\pa =(\pa_1,\pa_2)$, $\pa_j=-i{\pa\/\pa x_j}$\,, $j=1,2$.
Recall that the spectrum of the Dirac operator $D_m(\pa)$ is absolutely continuous
and has the form $\R\sm (-m,m)$.

\begin{proposition}\lb{Gra}
Let the Schr\"odinger operator $H=\D+Q$ act on the hexagonal lattice $\bG$, where
$Q$ satisfies \er{q12}. Then

i) The spectrum of the Schr\"odinger operator $H$ has the form
\[
\lb{SH}
\begin{aligned}
\s(H)=\s_{ac}(H)=[\l^-,\l^+]\sm \g,\\
\l^\pm=3\pm \sqrt{9+m^2}\,,\qqq \g=(3-m,3+m),\qqq m=q_1\,.
\end{aligned}
\]
ii) The Floquet $2\ts 2$  matrix  $H(\vt)$  satisfies
\[
\lb{Di}
H(\vt)=3\1_2+D_m(t)+O(|t|^2) \qqq
\textrm{as}\qq |t|\to 0,
\]
\[
\lb{Di1}
\textstyle t=(t_1,t_2)\in \R^2,\qqq
t_1=\frac12\,(\vt_1+\vt_2),\qqq
t_2=\frac{\sqrt{3}}2\,\big(\vt_1-\vt_2-\frac{4\pi}3\big),
\]
where $\1_2$ is the  $2\ts 2$ identity matrix  and  the matrix $D_m(t)$ is defined by \er{DDi}.
\end{proposition}

\no {\bf Proof.} i) The matrix $H(\vt)$ for
the hexagonal lattice $\bG$ has the form
\[
\lb{ell}
H(\vt)=\left(
\begin{array}{cc}
3+q_1  & -\D_{12} \\[6pt]
-\bar\D_{12} & 3-q_1
\end{array}\right)(\vt),\qqq \textstyle \D_{12}(\vt )=1+e^{i\vt _1}+
e^{i\vt_2},\qq
\forall \ \vt=(\vt_1,\vt_2)\in\T^2.
\]
This gives that the eigenvalues of each matrix $H(\vt)$ are given by
\[
\lb{evg}\textstyle \l_n(\vt)=3+(-1)^{n+1}\sqrt{q_1^2+|\D_{12}(\vt
)|^2}\;,  \qqq n=1,2,
\]
and the spectrum of the Schr\"odinger operator $H$ on the hexagonal
lattice $\bG$ has the form
$$
\s(H)=\s_{ac}(H)=[\l_1^-,\l_1^+]\cup[\l_2^-,\l_2^+],
$$
where
$$
\textstyle
\l_1^+=\max\limits_{\vt\in\T^2}\l_1(\vt)=\l_1({2\pi\/3},-{2\pi\/3})=3-q_1,\qqq
\l_1^-=\min\limits_{\vt\in\T^2}\l_1(\vt)
=\l_1(0)=3-\sqrt{9+q_1^2}\;,
$$
$$
\textstyle \l_2^+=\max\limits_{\vt\in\T^2}\l_2(\vt)=\l_2(0)=
3+\sqrt{9+q_1^2}\;,\qqq
\l_2^-=\min\limits_{\vt\in\T^2}\l_2(\vt)=\l_2({2\pi\/3},-{2\pi\/3})=3+q_1.
$$

ii) It is easy to show that
$$
\textstyle \D_{12}(\vt)=0 \qq
\Leftrightarrow \qq \vt=\pm\vt^0, \qqq \vt^0=(\vt^0_1,\vt^0_2)=
\big(\frac{2\pi}3\,,-\frac{2\pi}3\big)\in\T^2.
$$
The Taylor expansion for $-\D_{12}(\vt_0+s)$ in terms of
$s=(s_1,s_2)=\vt-\vt^0$ is given by
\[
\lb{ex}
\begin{aligned}
\textstyle-\D_{12}(\vt)=-(1+e^{i\vt _1}+e^{i\vt _2})=
-1-e^{i\vt^0_1}(1+s_1)-e^{i\vt^0_2}(1+s_2)+O(|s|^2)\\
\textstyle=-1-\frac12\,\big[(-1+ i\sqrt{3}\;)(1+s_1)- (1+
i\sqrt{3}\;)(1+s_2)\big]+O(|s|^2)= t_1-it_2+O(|t|^2),
\end{aligned}
\]
where the identities $ e^{\,\pm i\frac{2\pi}3}={1\/2}(-1\pm
i\sqrt{3})$ and $t=(t_1,t_2)\in \R^2$, $t_1={s_1+s_2\/2}$,
$t_2=\frac{\sqrt{3}}2\,(s_1-s_2)$ have been used. Thus, we obtain
$$
H(\vt)=\ma 3+q_1 & t_1-it_2 \\
         t_1+it_2 & 3-q_1 \am+O(|t|^2),
$$
which yields \er{Di}. Finally, we note that the Taylor expansion for
$-\D_{12}(\vt)$ about the point $-\vt^0$ is given by the same
asymptotics \er{ex}, but $t_2$ is defined by
$t_2=\vt_1-\vt_2+{4\pi\/3}$\,. \qq
$\BBox$

\

\no \textbf{Remark.} 1) The point $t=0$ is the so-called Dirac point for the Schr\"odinger operator $H=\D+Q$ on the graphene. In the nice survey \cite{CGPNG09} there is a discussion about the tight-banding approximations, the energy dispersion, different defects in graphene etc.

2) Let $q_1=0$, i.e., $H=\D$. From item i) it follows that the spectrum of the Laplacian $\D$ on $\bG$ is given by
$\s(\D)=\s_{ac}(\D)=[0,6]$.

3) Due to \er{SH} any perturbation $Q$ creates a gap in the spectrum of $H$.

%***********************************
\subsection{Body-centered cubic lattice.}
We consider the cubic lattice $\dL^3=(V,\cE)$, defined by \er{dLg}, see Fig.\ref{ff.0.1}\emph{a}. The spectrum of the Laplacian $\D$ on $\dL^3$ has the form
$\s(\D)=\s_{ac}(\D)=[0,12]$.
\setlength{\unitlength}{1.0mm}
\begin{figure}[h]
\centering

\unitlength 1.0mm % = 2.845pt
\linethickness{0.4pt}
\ifx\plotpoint\undefined\newsavebox{\plotpoint}\fi % GNUPLOT compatibility
\begin{picture}(140,60)(0,0)

%  убическа€ решетка

\put(10,10){\line(1,0){40.00}}
\put(10,30){\line(1,0){40.00}}
\put(10,50){\line(1,0){40.00}}

\bezier{60}(17,33.5)(37,33.5)(57,33.5)
\bezier{60}(17,13.5)(37,13.5)(57,13.5)

\put(17,53.5){\line(1,0){40.00}}
\put(10,10){\line(0,1){40.00}}
\put(30,10){\line(0,1){40.00}}
\put(50,10){\line(0,1){40.00}}
\put(57,13.5){\line(0,1){40.00}}
\bezier{12}(10,10)(13.5,11.75)(17,13.5)

\bezier{60}(17,13.5)(17,33.5)(17,53.5)
\bezier{60}(37,13.5)(37,33.5)(37,53.5)
\bezier{12}(30,10)(33.5,11.75)(37,13.5)
\put(50,10){\line(2,1){7.00}}
\bezier{12}(10,30)(13.5,31.75)(17,33.5)
\bezier{12}(30,30)(33.5,31.75)(37,33.5)

\put(50,30){\line(2,1){7.00}}
\put(10,50){\line(2,1){7.00}}
\put(30,50){\line(2,1){7.00}}
\put(50,50){\line(2,1){7.00}}
\put(10,10){\circle{0.8}}
\put(30,10){\circle{0.8}}
\put(50,10){\circle{0.8}}
\put(10,30){\circle{0.8}}
\put(30,30){\circle{0.8}}
\put(50,30){\circle{0.8}}
\put(10,50){\circle{0.8}}
\put(30,50){\circle{0.8}}
\put(50,50){\circle{0.8}}

\put(17,53.5){\circle{0.8}}
\put(37,53.5){\circle{0.8}}
\put(57,53.5){\circle{0.8}}
\put(17,33.5){\circle{0.8}}
\put(37,33.5){\circle{0.8}}
\put(57,33.5){\circle{0.8}}
\put(17,13.5){\circle{0.8}}
\put(37,13.5){\circle{0.8}}
\put(57,13.5){\circle{0.8}}

\put(50,10){\vector(0,1){20.00}}
\put(50,10){\vector(-1,0){20.00}}
\put(50,10){\vector(2,1){7.00}}

\put(49.0,7.0){$v$}
\put(53,10.0){$\scriptstyle a_1$}
\put(35.0,7.5){$\scriptstyle a_2$}
\put(46.5,26){$\scriptstyle a_3$}

\put(-5,5){(\emph{a})}
%*********************************
\bezier{30}(87,13.5)(97,13.5)(107,13.5)
\bezier{30}(87,13.5)(87,23.5)(87,33.5)
\bezier{12}(87,13.5)(83.5,11.75)(80,10)

\put(100,10){\vector(2,1){7.00}}
\bezier{12}(100,30)(103.5,31.75)(107,33.5)
\bezier{12}(80,30)(83.5,31.75)(87,33.5)
\bezier{30}(80,10)(80,20)(80,30)
\bezier{30}(107,13.5)(107,23.5)(107,33.5)
\put(100,10){\vector(-1,0){20.00}}
\bezier{30}(87,33.5)(97,33.5)(107,33.5)

\put(100,10){\vector(0,1){20.00}}
\bezier{30}(80,30)(90,30)(100,30)
\put(80,10){\circle{1}}
\put(80,30){\circle{1}}
\put(100,30){\circle{1}}
\put(87,13.5){\circle{1}}
\put(107,13.5){\circle{1}}
\put(107,33.5){\circle{1}}
\put(87,33.5){\circle{1}}

\put(100,10){\circle*{1}}

\put(100.0,7.5){$\scriptstyle v$}
\put(79.0,7.5){$\scriptstyle v$}
\put(101.0,28.5){$\scriptstyle v$}
\put(78.0,28.5){$\scriptstyle v$}

\put(85.0,33.5){$\scriptstyle v$}
\put(108.0,33.5){$\scriptstyle v$}
\put(108.0,13.5){$\scriptstyle v$}

\put(103.0,9.5){$\scriptstyle {\bf e}_1$}
\put(89.0,8.0){$\scriptstyle {\bf e}_2$}
\put(97.0,21.0){$\scriptstyle {\bf e}_3$}

\put(65,5){(\emph{b})}
\end{picture}

\vspace{-0.5cm} \caption{\footnotesize \emph{a}) Lattice
$\dL^3$;\quad \emph{b})  the fundamental graph $\dL^3_f$.}
\label{ff.0.1}
\end{figure}
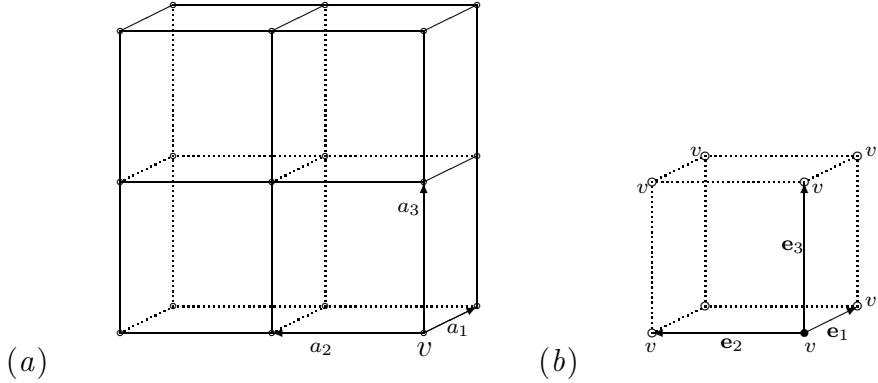

The body-centered cubic lattice $\bB_2$ is obtained from the cubic
lattice $\dL^3$ by adding one vertex in the center of each cube.
This vertex is connected with each corner vertex of the cube (see
Fig.\ref{ff.11}\emph{a}).  We consider the Schr\"odinger operators on
the body-centered cubic lattice $\bB_2$, where the potential $Q(v_j)=q_j$, $j\in
\N_2$. Without loss of generality we assume that $q_2=0$.

\setlength{\unitlength}{1.0mm}
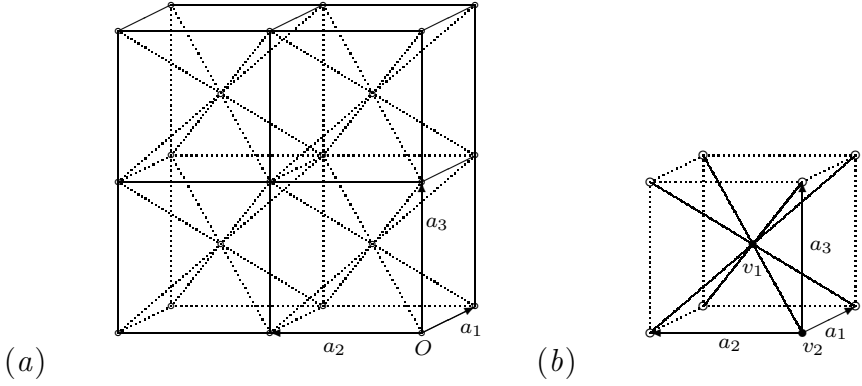
\begin{figure}[h]
\centering

\unitlength 1.0mm % = 2.845pt
\linethickness{0.4pt}
\ifx\plotpoint\undefined\newsavebox{\plotpoint}\fi % GNUPLOT compatibility
\begin{picture}(140,60)(0,0)

\put(10,10){\line(1,0){40.00}}
\put(10,30){\line(1,0){40.00}}
\put(10,50){\line(1,0){40.00}}

\bezier{60}(10,10)(23.5,21.75)(37,33.5)
\bezier{60}(10,30)(23.5,21.75)(37,13.5)
\bezier{40}(17,13.5)(23.5,21.75)(30,30.0)
\bezier{40}(17,33.5)(23.5,21.75)(30,10.0)
\put(23.5,21.75){\circle{0.8}}

\bezier{60}(10,30)(23.5,41.75)(37,53.5)
\bezier{60}(10,50)(23.5,41.75)(37,33.5)
\bezier{40}(17,33.5)(23.5,41.75)(30,50.0)
\bezier{40}(17,53.5)(23.5,41.75)(30,30.0)
\put(23.5,41.75){\circle{0.8}}

\bezier{60}(30,10)(43.5,21.75)(57,33.5)
\bezier{60}(30,30)(43.5,21.75)(57,13.5)
\bezier{40}(37,13.5)(43.5,21.75)(50,30.0)
\bezier{40}(37,33.5)(43.5,21.75)(50,10.0)
\put(43.5,21.75){\circle{0.8}}

\bezier{60}(30,30)(43.5,41.75)(57,53.5)
\bezier{60}(30,50)(43.5,41.75)(57,33.5)
\bezier{40}(37,33.5)(43.5,41.75)(50,50.0)
\bezier{40}(37,53.5)(43.5,41.75)(50,30.0)
\put(43.5,41.75){\circle{0.8}}

\bezier{60}(17,33.5)(37,33.5)(57,33.5)
\bezier{60}(17,13.5)(37,13.5)(57,13.5)

\put(17,53.5){\line(1,0){40.00}}
\put(10,10){\line(0,1){40.00}}
\put(30,10){\line(0,1){40.00}}
\put(50,10){\line(0,1){40.00}}
\put(57,13.5){\line(0,1){40.00}}
\bezier{12}(10,10)(13.5,11.75)(17,13.5)

\bezier{60}(17,13.5)(17,33.5)(17,53.5)
\bezier{60}(37,13.5)(37,33.5)(37,53.5)

\bezier{12}(30,10)(33.5,11.75)(37,13.5)
\put(50,10){\line(2,1){7.00}}
\bezier{12}(10,30)(13.5,31.75)(17,33.5)
\bezier{12}(30,30)(33.5,31.75)(37,33.5)

\put(50,30){\line(2,1){7.00}}
\put(10,50){\line(2,1){7.00}}
\put(30,50){\line(2,1){7.00}}
\put(50,50){\line(2,1){7.00}}
\put(10,10){\circle{0.8}}
\put(30,10){\circle{0.8}}
\put(50,10){\circle{0.8}}
\put(10,30){\circle{0.8}}
\put(30,30){\circle{0.8}}
\put(50,30){\circle{0.8}}
\put(10,50){\circle{0.8}}
\put(30,50){\circle{0.8}}
\put(50,50){\circle{0.8}}

\put(17,53.5){\circle{0.8}}
\put(37,53.5){\circle{0.8}}
\put(57,53.5){\circle{0.8}}
\put(17,33.5){\circle{0.8}}
\put(37,33.5){\circle{0.8}}
\put(57,33.5){\circle{0.8}}
\put(17,13.5){\circle{0.8}}
\put(37,13.5){\circle{0.8}}
\put(57,13.5){\circle{0.8}}

\put(50,10){\vector(0,1){20.00}}
\put(50,10){\vector(-1,0){20.00}}
\put(50,10){\vector(2,1){7.00}}

\put(49,7.0){$\scriptstyle O$}
\put(55,10.0){$\scriptstyle a_1$}
\put(37.0,7.5){$\scriptstyle a_2$}
\put(50.5,24){$\scriptstyle a_3$}

\put(-5,5){(\emph{a})}
%*********************************
\bezier{160}(80,10)(93.5,21.75)(107,33.5)
\bezier{160}(80,30)(93.5,21.75)(107,13.5)
\bezier{140}(87,13.5)(93.5,21.75)(100,30.0)
\bezier{140}(87,33.5)(93.5,21.75)(100,10.0)
\put(93.5,21.75){\circle*{1.0}}

\bezier{30}(87,13.5)(97,13.5)(107,13.5)
\bezier{30}(87,13.5)(87,23.5)(87,33.5)
\bezier{12}(87,13.5)(83.5,11.75)(80,10)

\put(100,10){\vector(2,1){7.00}}
\bezier{12}(100,30)(103.5,31.75)(107,33.5)
\bezier{12}(80,30)(83.5,31.75)(87,33.5)
\bezier{30}(80,10)(80,20)(80,30)
\bezier{30}(107,13.5)(107,23.5)(107,33.5)
\put(100,10){\vector(-1,0){20.00}}
\bezier{30}(87,33.5)(97,33.5)(107,33.5)

\put(100,10){\vector(0,1){20.00}}
\bezier{30}(80,30)(90,30)(100,30)
\put(80,10){\circle{1}}
\put(80,30){\circle{1}}
\put(100,30){\circle{1}}
\put(87,13.5){\circle{1}}
\put(107,13.5){\circle{1}}
\put(107,33.5){\circle{1}}
\put(87,33.5){\circle{1}}

\put(100,10){\circle*{1}}

\put(100.0,7.5){$\scriptstyle v_2$}
\put(92.1,18.0){$\scriptstyle v_1$}

\put(103.0,9.5){$\scriptstyle a_1$}
\put(89.0,8.0){$\scriptstyle a_2$}
\put(101.0,21.0){$\scriptstyle a_3$}

\put(65,5){(\emph{b})}
\end{picture}

\vspace{-0.5cm} \caption{\footnotesize \emph{a}) Body-centered cubic lattice $\bB_2$ ;\quad \emph{b})  the fundamental graph.}
\label{ff.11}
\end{figure}

\begin{proposition}\lb{BCC}
i) The spectrum of $\D$ on the body-centered cubic lattice
$\bB_2$ has the form
\[
\lb{sBCCx}\textstyle \s(\D)=\s_{ac}(\D)=\s_1^0\cup\s_2^0,\qqq
\s_1^0=[0,8], \qqq \s_2^0=[12,20].
\]
ii) The spectrum of Schr\"odinger operators $H=\D+Q$ on  $\bB_2$ has the form
\[
\lb{SB}
\s(H)=\s_{ac}(H)=\s_1\cup\s_2=[\l^-_1,\l^+_2]\sm \g_1,
\]
where
$$
\textstyle\l_1^-=8+p-\sqrt{p^2+64}\;,\qqq
\l_2^+=8+p+\sqrt{p^2+64}\;,\qqq p={q_1\/2}\,,
$$
and the gap $\g_1$ satisfies (without loss of generality we may assume
 that $q_2=0$)
$$
\g_1=(\l^+_1,\l^-_2)=\ca (20,q_1+8), & if \ 12<q_1\\
                         (q_1,q_1)=\varnothing, & if \ 4\leq q_1\leq12\\
                         (q_1+8,12),& if \ q_1<4  \ac.
$$
\end{proposition}

\no {\bf Proof.} i) The fundamental graph $\G_f$ of the
body-centered cubic lattice $\bB_2$  consists of 2 vertices
$v_1,v_2$ with the degrees $\vk_1=8$, $\vk_2=14$ and 11 oriented
edges
$$
\be_1=\be_2=\be_3=(v_2,v_2),\qqq \be_4=\ldots=\be_{11}=(v_1,v_2)
$$
and their inverse edges. The indices of the fundamental graph edges in the coordinate system with the origin $O$ (Fig.\ref{ff.11}\emph{a}) are given by
$$
\begin{aligned}
\t(\be_1)=\t(\be_5)=(1,0,0),\quad \t(\be_2)=\t(\be_6)=(0,1,0),\quad
\t(\be_3)= \t(\be_7)=(0,0,1), \\[4pt]  \t(\be_4)=(0,0,0),\qq
\t(\be_8)=(1,1,0),\qq \t(\be_9)=(1,0,1),\quad \t(\be_{10})=
(0,1,1),\quad \t(\be_{11})=(1,1,1).
\end{aligned}
$$
For each $\vt\in\T^3$ the matrix $\D(\vt)$ has the form
\[
\lb{el}
\D(\vt)=\left(
\begin{array}{cc}
8  & \D_{12} \\[6pt]
\bar\D_{12} & \D_{22}
\end{array}\right)(\vt)\qqq    \ca
\D_{22}(\vt)=14-2c_0,\\
\D_{12}(\vt)=-(1+e^{i\vt_1})(1+e^{i\vt_2})(1+e^{i\vt_3})\ac,
\]
 where
\[
\lb{pol}
c_0=c_1+c_2+c_3,\qqq c_j=\cos\vt_j,\qqq j=1,2,3.
\]
The direct calculation gives
$$
\det\big(\D(\vt)-\l \1_2\big)=\l^2-(22-2c_0)\l+
8\big(14-2c_0-(1+c_1)(1+c_2)(1+c_3)\big).
$$
The eigenvalues of the matrix $\D(\vt)$ are given by
\[
\lb{ev}\textstyle
\l_n^0(\vt)=11-c_0+(-1)^{n}\sqrt{(3-c_0)^2+
8\,(1+c_1)(1+c_2)(1+c_3)}\;,  \qq n=1,2.
\]
Thus, the spectrum of the Laplacian on the body-centered cubic
lattice  $\bB_2$ has the form
$$
\s(\D)=\s_{ac}(\D)=[\l_1^{0-},\l_1^{0+}]\cup[\l_2^{0-},\l_2^{0+}],\qq
[\l_n^{0-},\l_n^{0+}]=\l_n^0(\T^3),\qq n=1,2.
$$
From the Proposition \ref{MP}.v it follows that
\[
\lb{razz}
\l_1^0(\vt)\le8\le\l_2^0(\vt),\qqq \forall\,\vt\in\T^3.
\]
Then, investigating the functions $\l_1^0,\l_2^0$ defined by \er{ev} on the extremes and using \er{razz}, we obtain
$$
\l_1^{0-}=\l_1(0)=0,\qq
\l_1^{0+}=\l_1(\pi,\pi,\pi)=8,\qq
%$$
%$$
\textstyle
\l_2^{0-}=\l_2(\pi,0,0)=12,\qq
\l_2^{0+}=\l_2(\pi,\pi,\pi)=20.
$$

ii) The matrix $H(\vt)$ has the form
$$
H(\vt)=\left(
\begin{array}{cc}
8+q_1  & \D_{12} \\[6pt]
\bar\D_{12} & \D_{22}
\end{array}\right)(\vt),
$$
where $\D_{12},\D_{22}$ are given in \er{el}.
By a direct calculation, we get
$$
\textstyle
\det(H(\vt)-\l\1_2)=\l^2-2(11-c_0+p)\l
-8(1+c_1)(1+c_2)(1+c_3)+112-16c_0+28p-4c_0p,
$$
where $c_0$ and $c_j$ are defined by \er{pol} and $p=q_1/2$.
The eigenvalues of the matrix $H(\vt)$ are given by
\[
\lb{ev1}
\l_n(\vt)=11-c_0+p+(-1)^n
\sqrt{\big(c_0+p-3\big)^2+8(1+c_1)(1+c_2)(1+c_3)}\;,  \qq n=1,2.
\]
Due to Propositions \ref{LP} and \ref{MP}.v, we have
$$
\l_1^-=\l_1(0)=8+p-\sqrt{p^2+64}\;,
$$
$$
\l_1(\vt)\leq8+2p\leq\l_2(\vt),\qqq \forall\,\vt\in\T^3.
$$
The direct calculation gives
$$
\l_1(\vt_0)=8+2p \ \textrm{ for some } \vt_0 \qq \Leftrightarrow \qq p\leq6,
$$
$$
\l_2(\vt_0)=8+2p \ \textrm{ for some } \vt_0 \qq \Leftrightarrow \qq p\geq2.
$$
Then, an extremes properties of the functions $\l_1,\l_2$  gives
$$
\l_1^+=\l_1(\pi,\pi,\pi)=
\left\{\begin{array}{rl}
         8+2p, & \textrm{if} \qq p\leq6\\[4pt]
         20, & \textrm{if} \qq  p>6
       \end{array}\right.,
$$
$$
\l_2^+=\l_2(0)=8+p+\sqrt{p^2+64},\qqq
\l_2^-=\l_2(\pi,0,0)=
\left\{\begin{array}{rl}
         8+2q, & \textrm{if} \qq  2\leq p\\[4pt]
         12, & \textrm{if} \qq  p<2
       \end{array}\right.,
$$
which proves ii). \qq
$\BBox$

\

%********************************************************

\subsection{Face-centered cubic lattice}
The face-centered cubic lattice
$\bB_4$ is obtained from the cubic lattice $\dL^3$ by adding one
vertex in the center of each cube face. This vertex is connected
with each corner vertex of the cube face (Fig.\ref{ff.FCC}\emph{a}).
The fundamental graph $\G_f$ of $\bB_4$ consists of 4 vertices
$v_1,v_2,v_3,v_4$ and 15 edges (Fig.\ref{ff.FCC}\emph{b}).

We consider the Schr\"odinger operators $H=\D+Q$ on the face-centered cubic lattice $\bB_4$.
%,
%where the potential $Q(v_j)=q_j$, $j\in \N_4$. Without loss of generality we %assume that $q_4=0$.

\setlength{\unitlength}{1.0mm}
\begin{figure}[h]
\centering

\unitlength 1.0mm % = 2.845pt
\linethickness{0.4pt}
\ifx\plotpoint\undefined\newsavebox{\plotpoint}\fi % GNUPLOT compatibility
\begin{picture}(140,60)(0,0)

% face-centered решетка
\bezier{80}(17,13.5)(37,33.5)(57,53.5)
\bezier{40}(17,33.5)(27,43.5)(37,53.5)
\bezier{40}(37,13.5)(47,23.5)(57,33.5)

\bezier{80}(57,13.5)(37,33.5)(17,53.5)
\bezier{40}(17,33.5)(27,23.5)(37,13.5)
\bezier{40}(37,53.5)(47,43.5)(57,33.5)

\bezier{40}(10,10)(23.5,11.75)(37,13.5)
\bezier{40}(30,10)(43.5,11.75)(57,13.5)
\bezier{20}(30,10)(23.5,11.75)(17,13.5)
\bezier{20}(50,10)(43.5,11.75)(37,13.5)

\bezier{40}(10,30)(23.5,31.75)(37,33.5)
\bezier{40}(30,30)(43.5,31.75)(57,33.5)
\bezier{20}(30,30)(23.5,31.75)(17,33.5)
\bezier{20}(50,30)(43.5,31.75)(37,33.5)

\bezier{180}(10,50)(23.5,51.75)(37,53.5)
\bezier{180}(30,50)(43.5,51.75)(57,53.5)
\bezier{180}(30,50)(23.5,51.75)(17,53.5)
\bezier{180}(50,50)(43.5,51.75)(37,53.5)

\bezier{40}(10,10)(13.5,21.75)(17,33.5)
\bezier{40}(10,30)(13.5,21.75)(17,13.5)
\bezier{40}(10,30)(13.5,41.75)(17,53.5)
\bezier{40}(10,50)(13.5,41.75)(17,33.5)

\bezier{40}(30,10)(33.5,21.75)(37,33.5)
\bezier{40}(30,30)(33.5,21.75)(37,13.5)
\bezier{40}(30,30)(33.5,41.75)(37,53.5)
\bezier{40}(30,50)(33.5,41.75)(37,33.5)

\bezier{140}(50,10)(53.5,21.75)(57,33.5)
\bezier{140}(50,30)(53.5,21.75)(57,13.5)
\bezier{140}(50,30)(53.5,41.75)(57,53.5)
\bezier{140}(50,50)(53.5,41.75)(57,33.5)

\put(10,10){\line(1,1){40.00}}
\put(10,30){\line(1,-1){20.00}}
\put(30,50){\line(1,-1){20.00}}
\put(10,50){\line(1,-1){40.00}}
\put(10,30){\line(1,1){20.00}}
\put(30,10){\line(1,1){20.00}}

\put(10,10){\line(1,0){40.00}}
\put(10,30){\line(1,0){40.00}}
\put(10,50){\line(1,0){40.00}}
\bezier{60}(17,33.5)(37,33.5)(57,33.5)
\bezier{60}(17,13.5)(37,13.5)(57,13.5)
\put(17,53.5){\line(1,0){40.00}}
\put(10,10){\line(0,1){40.00}}
\put(30,10){\line(0,1){40.00}}
\put(50,10){\line(0,1){40.00}}
\put(57,13.5){\line(0,1){40.00}}
\bezier{12}(10,10)(13.5,11.75)(17,13.5)
\bezier{60}(17,13.5)(17,33.5)(17,53.5)
\bezier{60}(37,13.5)(37,33.5)(37,53.5)

\bezier{12}(30,10)(33.5,11.75)(37,13.5)
\put(50,10){\line(2,1){7.00}}
\bezier{12}(10,30)(13.5,31.75)(17,33.5)
\bezier{12}(30,30)(33.5,31.75)(37,33.5)

\put(50,30){\line(2,1){7.00}}
\put(10,50){\line(2,1){7.00}}
\put(30,50){\line(2,1){7.00}}
\put(50,50){\line(2,1){7.00}}
\put(10,10){\circle{0.8}}
\put(30,10){\circle{0.8}}
\put(50,10){\circle{0.8}}
\put(10,30){\circle{0.8}}
\put(30,30){\circle{0.8}}
\put(50,30){\circle{0.8}}
\put(10,50){\circle{0.8}}
\put(30,50){\circle{0.8}}
\put(50,50){\circle{0.8}}
\put(20,20){\circle{0.8}}
\put(40,20){\circle{0.8}}

\put(20,40){\circle{0.8}}
\put(40,40){\circle{0.8}}

\put(27,23.5){\circle{0.8}}
\put(47,23.5){\circle{0.8}}

\put(27,43.5){\circle{0.8}}
\put(47,43.5){\circle{0.8}}

\put(23.5,11.75){\circle{0.8}}
\put(43.5,11.75){\circle{0.8}}
\put(23.5,31.75){\circle{0.8}}
\put(43.5,31.75){\circle{0.8}}
\put(23.5,51.75){\circle{0.8}}
\put(43.5,51.75){\circle{0.8}}

\put(13.5,21.75){\circle{0.8}}
\put(13.5,41.75){\circle{0.8}}
\put(33.5,21.75){\circle{0.8}}
\put(33.5,41.75){\circle{0.8}}
\put(53.5,21.75){\circle{0.8}}
\put(53.5,41.75){\circle{0.8}}

\put(17,53.5){\circle{0.8}}
\put(37,53.5){\circle{0.8}}
\put(57,53.5){\circle{0.8}}
\put(17,33.5){\circle{0.8}}
\put(37,33.5){\circle{0.8}}
\put(57,33.5){\circle{0.8}}
\put(17,13.5){\circle{0.8}}
\put(37,13.5){\circle{0.8}}
\put(57,13.5){\circle{0.8}}

\put(50,10){\vector(0,1){20.00}}
\put(50,10){\vector(-1,0){20.00}}
\put(50,10){\vector(2,1){7.00}}

\put(49,7.0){$\scriptstyle O$}
\put(53,10.0){$\scriptstyle a_1$}
\put(35.0,7.5){$\scriptstyle a_2$}
\put(46.5,19){$\scriptstyle a_3$}

\put(-5,5){(\emph{a})}
%*********************************
\put(80,10){\line(1,1){20.00}}
\put(100,10){\line(-1,1){20.00}}
\bezier{180}(80,10)(93.5,11.75)(107,13.5)
\bezier{180}(87,13.5)(93.5,11.75)(100,10)
\bezier{30}(87,13.5)(97,13.5)(107,13.5)
\bezier{30}(87,13.5)(87,23.5)(87,33.5)
\bezier{12}(87,13.5)(83.5,11.75)(80,10)

\put(100,10){\vector(2,1){7.00}}
\bezier{12}(100,30)(103.5,31.75)(107,33.5)
\bezier{12}(80,30)(83.5,31.75)(87,33.5)
\bezier{30}(80,10)(80,20)(80,30)
\bezier{30}(107,13.5)(107,23.5)(107,33.5)
\put(100,10){\vector(-1,0){20.00}}
\bezier{30}(87,33.5)(97,33.5)(107,33.5)

\bezier{180}(100,10)(103.5,21.75)(107,33.5)
\bezier{180}(100,30)(103.5,21.75)(107,13.5)
\put(100,10){\vector(0,1){20.00}}
\bezier{30}(80,30)(90,30)(100,30)
\put(80,10){\circle{1}}
\put(80,30){\circle{1}}
\put(100,30){\circle{1}}
\put(87,13.5){\circle{1}}
\put(107,13.5){\circle{1}}
\put(107,33.5){\circle{1}}
\put(87,33.5){\circle{1}}
\put(100,10){\circle*{1}}
\put(90,20){\circle*{1}}
\put(100.0,8.0){$\scriptstyle v_4$}
\put(92.5,13.0){$\scriptstyle v_1$}
\put(88.7,21.7){$\scriptstyle v_3$}
\put(104.2,21.3){$\scriptstyle v_2$}
\put(104.0,10.5){$\scriptstyle a_1$}
\put(85.0,8.0){$\scriptstyle a_2$}
\put(97.0,25.0){$\scriptstyle a_3$}
\put(93.4,11.8){\circle*{1}}
\put(103.5,21.8){\circle*{1}}

\put(65,5){(\emph{b})}
\end{picture}

\vspace{-0.5cm} \caption{\footnotesize  \emph{a}) Face-centered cubic lattice;\quad \emph{b})  the fundamental graph.}
\label{ff.FCC}
\end{figure}
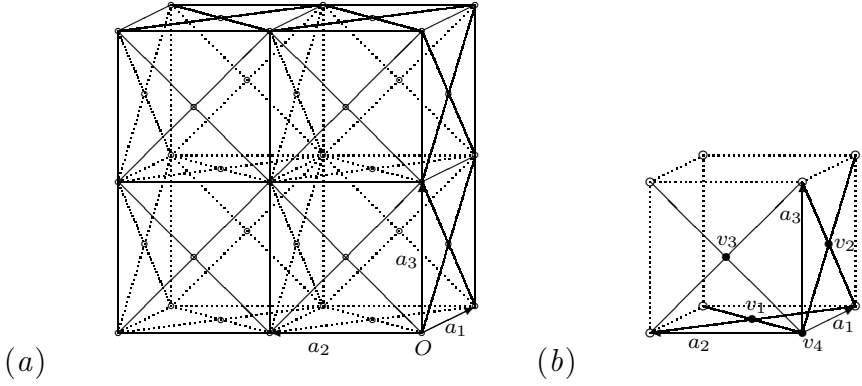

 \begin{proposition}\lb{FCC}
i) The spectrum of $\D$ on the face-centered cubic lattice $\bB_4$ has the form
\[
\lb{sBCC}\textstyle
\s(\D)=\s_{ac}(\D)\cup\s_{fb}(\D),\qqq
\s_{ac}(\D)=[0,4]\cup[16,24],\qqq
\s_{fb}(\D)=\{4\},
\]
where the flat band 4 has multiplicity 2.

ii) Let the Schr\"odinger operator $H=\D+Q$ act
on $\bB_4$, where the potential $Q$ satisfies
\[
Q(v_j)=q_j, \qqq j\in \N_4, \qqq q_4=0.
\]
Then the  spectrum of the Schr\"odinger operator $H$  has the form
\[
\lb{sBCC1}\textstyle
\s(H)=\s_{ac}(H)\cup\s_{fb}(H),
\]
where
\[\lb{fla}
\s_{fb}(H)=\left\{
  \begin{array}{cl}
    \varnothing, & \textrm{ if all values $q_1,q_2,q_3$ are distinct
    (generic potential) }\\[4pt]
    \{q_{(1)}+4\}, & \textrm{ if $q_j=q_k=q_{(1)}\neq q_n$ for some $j,k,n=1,2,3$, $n\ne j\neq k$, $k\neq n$ }\\[4pt]
    \{q_{(2)}+4\} , & \textrm{ if $q_1=q_2=q_3=q_{(2)}$  }\\
  \end{array}\right..
\]
Moreover, the flat band $q_{(1)}+4$ has multiplicity 1, the flat band $q_{(2)}+4$ has multiplicity 2.
\end{proposition}

\textbf{Remark.} The spectrum of the Laplacian $\D$ on $\bB_4$ has the flat band $4$ of multiplicity 2.  If $Q$ is a generic potential, then the spectrum of $H=\D+Q$ is absolutely continuous. In the case of non-generic potential the operator $H$ has a flat band.

\

\no {\bf Proof of Proposition \ref{FCC}.} i) The fundamental graph
$\G_f$ consists of 4 vertices $v_1,v_2,v_3,v_4$ with the degrees
$\vk_1=\vk_2=\vk_3=4$, $\vk_4=18$ and 15 oriented edges
$$
\begin{aligned}
\be_1=\be_2=\be_3=(v_4,v_4),\qq \be_4=\ldots=\be_7=(v_1,v_4),\\
\be_8=\ldots=\be_{11}=(v_2,v_4),\qq\be_{12}=\ldots=\be_{15}=(v_3,v_4)
\end{aligned}
$$
and their inverse edges. The indices of the fundamental graph edges in the
coordinate system with the origin $O$ (Fig.\ref{ff.FCC}\emph{a}) are given by
$$
\begin{aligned}
\t(\be_1)=\t(\be_5)=\t(\be_9)=(1,0,0),\quad \t(\be_2)=\t(\be_6)=
\t(\be_{13})=(0,1,0),\\[4pt]
\t(\be_3)=\t(\be_{10})=\t(\be_{14})=(0,0,1),\qq
\t(\be_4)=\t(\be_8)=\t(\be_{12})=(0,0,0),\\[4pt]
\t(\be_7)=(1,1,0),\quad \t(\be_{11})=(1,0,1),\qq
  \t(\be_{15})=(0,1,1).
\end{aligned}
$$
For each $\vt\in\T^3$ the matrix $\D(\vt)$ has the form
\[
\lb{DDD}
\D(\vt)=\left(
\begin{array}{cccc}
4 & 0  & 0 & \D_{14} \\[6pt]
0 & 4  & 0 & \D_{24} \\[6pt]
0 & 0  & 4 & \D_{34} \\[6pt]
\bar\D_{14} & \bar\D_{24} & \bar\D_{34} & \D_{44} \\[6pt]
\end{array}\right)(\vt),\qqq
\ca \D_{14}(\vt)=-(1+e^{i\vt_1})(1+e^{i\vt_2})\\
\D_{24}(\vt)=-(1+e^{i\vt_1})(1+e^{i\vt_3})\\
\D_{34}(\vt)=-(1+e^{i\vt_2})(1+e^{i\vt_3})\\
\D_{44}(\vt)=18-2c_0\ac,
\]
where recall that $c_0=c_1+c_2+c_3$, $c_j=\cos\vt_j$, $j=1,2,3$. By
a direct calculation, we get
$$
\det\big(\D(\vt)-\l\1_4\big)=(4-\l)^2\big(\l^2-2(11-c_0)\l+4(15-\e-4c_0)\big),
$$
where
$$
\e=\e(\vt)=c_1c_2+c_1c_3+c_2c_3.
$$

The eigenvalues of each matrix $\D(\vt)$ are given by
\[
\lb{evf}\textstyle
\L_s^0(\vt)=(11-c_0)+(-1)^s\sqrt{(3-c_0)^2+52+4\e}\;,  \qq s=1,2,\qqq \L_3^0(\vt)=\L_4^0(\vt)=4.
\]
Thus, the spectrum of the Laplacian on the face-centered cubic lattice $\bB_4$ has the form
$$
\s(\D)=\s_{ac}(\D)\cup\s_{fb}(\D),\qqq
\s_{ac}(\D)=[\L_1^{0-},\L_1^{0+}]\cup[\L_2^{0-},\L_2^{0+}],\qqq
\s_{fb}(\D)=\{4\},
$$
where the flat band 4 has multiplicity 2 and $[\L_s^{0-},\L_s^{0+}]=\L_s(\T^3), s=1,2$.
From Proposition \ref{MP}.v it follows that
\[
\lb{raz2}
\L_1^0(\vt)\le4\le\L_2^0(\vt),\qqq \forall\,\vt\in\T^3.
\]
Then, using Propositions \ref{LP}  and \er{raz2}, we have
$$
\L_1^{0-}=\L_1^0(0)=0,\qqq
\L_1^{0+}=\L_1^0(\pi,\pi,\pi)=4.
$$
A direct calculation yields
$$\textstyle
\L_2^{0-}=\L_2^0(0)=16,\qqq
\L_2^{0+}=\L_2^0(\pi,\pi,\pi)=24.
$$

ii) The matrix $H(\vt)$ has the form
$$
H(\vt)=\D(\vt)+q,\qqq q=\diag(q_1,q_2,q_3,0),
$$
where $\D(\vt)$ is defined by \er{DDD}. We write the characteristic
polynomial of the matrix $H(\vt)$ in the form of the linear
combination of the linearly independent functions
\begin{multline}
\lb{lco} \det(H(\vt)-\l\1_4)=\e_1+2\e_2c_1+2\e_3c_2+
2\e_4c_3+4\e_5c_1c_2+ 4\e_6c_1c_3+4\e_7c_2c_3,
\end{multline}
where  $c_j=\cos\vt_j$, $j=1,2,3$, and $\e_s=\e_s(\l,q)$, $s\in\N_7$, are
given by
$$
\e_1=\l^4-(\z_1+30)\l^3+
(26\z_1+\z_2+252)\l^2-(152\z_1+22\z_2+\z_3+832)\l+256\z_1+68\z_2+18\z_3+960,
$$
$$
\e_2=-\l^3+(\z_1+8)\l^2+
(-6\z_1-\z_2+2q_3-16)\l+8q_1+8q_2+2\z_2+2q_1q_2+\z_3,
$$
$$
\e_3=-\l^3+(\z_1+8)\l^2+
(-6\z_1-\z_2+2q_2-16)\l+8q_1+8q_3+2\z_2+2q_1q_3+\z_3,
$$
$$
\e_4=-\l^3+(\z_1+8)\l^2+
(-6\z_1-\z_2+2q_1-16)\l+8q_2+8q_3+2\z_2+2q_2q_3+\z_3,
$$
$$
\e_5=-\l^2+(q_2+q_3+8)\l-4(q_2+q_3)-q_2q_3-16,
$$
$$
\e_6=-\l^2+(q_1+q_3+8)\l-4(q_1+q_3)-q_1q_3-16,
$$
$$
\e_7=-\l^2+(q_1+q_2+8)\l-4(q_1+q_2)-q_1q_2-16,
$$
$$
\z_1=q_1+q_2+q_3,\qqq \z_2=q_1q_2+q_1q_3+q_2q_3,\qqq \z_3=q_1q_2q_3.
$$
A point $\l$ is a flat band of the Schr\"odinger operator $H$ iff
$\det(H(\vt)-\l\1_4)=0$ for all $\vt\in\T^3$.
Since the linear combination \er{lco} of the linearly independent functions
is equal to 0, then we obtain the system of equations $\e_s(\l,q)=0$,  $s\in\N_7$.
All solutions of this system of equations have the form
$$
\l=q_1+4,\; q_1=q_2, \; q_3\in\R;  \qq
\l=q_1+4,\;q_1=q_3, \; q_2\in\R; \qq
\l=q_2+4,\;q_2=q_3, \; q_1\in\R,
$$
which yields \er{fla}.\qq
$\BBox$

%******************************************************************
\section{\lb{Sec8} Appendix, well-known  properties of matrices}
\setcounter{equation}{0}

\subsection{Properties of  matrices.}
We recall some well-known  properties of matrices (see e.g., \cite{HJ85}).
Denote by
$\l_1(A)\leq\ldots\leq\l_\n(A)$ the eigenvalues of a self-adjoint $\n\ts\n$ matrix $A$, arranged in
increasing order, counting multiplicities. Let $\r(A)$ be the spectral radius of~$A$.

\begin{proposition}\lb{MP}
i) Let $A=\{A_{jk}\}$, $A_+=\big\{|A_{jk}|\big\}$ and $B=\{B_{jk}\}$
be $\nu\ts\nu$ matrices. If ${|A_{jk}|\leq B_{jk}}$ for all
$j,k\in\N_\n$, then $\r(A)\leq\r(A_+)\leq\r(B)$ (see Theorem 8.1.18
in \cite{HJ85}).

ii) Let $A=\{A_{jk}\}$ and $B=\{B_{jk}\}$ be $\nu\ts\nu$ matrices with nonnegative entries and $A$ be irreducible.
If $B\neq0$, then $\r(A+B)>\r(A)$
(see p.515 in \cite{HJ85}).

iii) Let $A,B$ be self-adjoint $\nu\ts\nu$ matrices and let $B\geq0$. Then the eigenvalues $\l_n(A)\leq\l_n(A+B)$ for all $n\in\N_\n$ (see Corollary 4.3.3 in \cite{HJ85}).

iv) Let $A,B$ be self-adjoint $\nu\ts\nu$ matrices. Then for each $n\in\N_\n$ we have
$$
\l_n(A)+\l_1(B)\leq\l_n(A+B)\leq\l_n(A)+\l_\n(B)
$$
(see Theorem 4.3.1 in \cite{HJ85}).

v) Let $B=\ma   A & y \\
  y^\ast & a \am$ be a self-adjoint $(\nu+1)\ts(\nu+1)$ matrix
for some self-adjoint $\nu\ts \nu$ matrix $A$, some  real number $a$
and some vector $y\in\C^{\nu}$.
Then
$$
\l_1(B)\leq\l_1(A)\leq\l_2(B)\leq\ldots\leq\l_\n(B)\leq\l_\n(A)\leq\l_{\nu+1}(B)
$$
(see Theorem 4.3.8 in \cite{HJ85}).

vi) Let $A=\{A_{jk}\}$ be a $\nu\ts\nu$ matrix and let $\G(A)$ be a graph
on $\n$ vertices $v_1,\ldots,v_\n$ such that there is an edge
$(v_j,v_k)$ in $\G(A)$ iff $A_{jk}\neq0$. Then $A$ is irreducible iff
$\G(A)$ is connected (see Theorem 6.2.24 in \cite{HJ85}).

vii) Let $A$ be a irreducible $\nu\ts\nu$ matrix with nonnegative entries.
Then the spectral radius $\r(A)$ is a simple eigenvalue of $A$
(see Theorem 8.4.4 in \cite{HJ85}).

viii) Let $M=\ma   A & B \\
         C & D       \am$ be a $\n\ts\n$ matrix
for some square matrices $A,D$ and some matrices $B,C$. Then
\[
\lb{det}
\det M=\det A\cdot\det\big(D-CA^{-1}B\big)
\]
(see pp.21--22 in \cite{HJ85}).

\end{proposition}

\begin{proposition} (see \cite{K13})\lb{TK}
i) Let $V=\{V_{jk}\}$ be a self-adjoint $\n\ts \n$  matrix, for some
$\n<\iy$ and let $ B=\diag \{B_1,\ldots, B_\n\},\  B_j=\sum\limits_{k=1}^\n
|V_{jk}|$. Then the following estimates hold true:
\[
\begin{aligned}
\lb{V} &-B\le V\le B.
\end{aligned}
\]

ii) Let, in addition, $H_0$, $H=H_0+V$ be self-adjoint $\n\ts\n$
matrices. Then the following estimate holds true:
\[
\begin{aligned}
\lb{1} \sum_{n=1}^\n |\l_n(H)-\l_n(H_0)|\le 2
\sum\limits_{j,k=1}^\n|V_{jk}|.
\end{aligned}
\]
\end{proposition}

%*****************************************************

\medskip

\footnotesize \no\textbf{Acknowledgments.} \footnotesize Various
parts of this paper were written during Evgeny Korotyaev's stay  in
the Mathematical Institute of Tsukuba University, Japan  and
Mittag-Leffler Institute, Sweden and Centre for Quantum Geometry of
Moduli spaces (QGM),  Aarhus University, Denmark. He is grateful to
the institutes for the hospitality. His study was partly supported
by The Ministry of education and science of Russian Federation,
project 07.09.2012 No 8501 and the RFFI grant "Spectral and
asymptotic methods for studying of the differential operators" No
11-01-00458 and the Danish National Research Foundation grant DNRF95
(Centre for Quantum Geometry of Moduli Spaces - QGM)Ф.

\end{document}